\documentclass[10pt, openany]{amsbook}
\usepackage{mathpazo}

\usepackage{amssymb, enumerate, xspace, graphicx, url}
\usepackage[all]{xy}
\SelectTips{cm}{}
\makeindex

\numberwithin{equation}{section}

2

\numberwithin{subsection}{section}

\numberwithin{section}{chapter}

\allowdisplaybreaks[1]

\newenvironment{enumeratea}
{\begin{enumerate}[\upshape (a)]}
{\end{enumerate}}

\newenvironment{enumeratei}
{\begin{enumerate}[\upshape (i)]}
{\end{enumerate}}

\newtheorem*{namedtheorem}{\theoremname}
\newcommand{\theoremname}{testing}
\newenvironment{named}[1]{\renewcommand\theoremname{#1}
\begin{namedtheorem}}
{\end{namedtheorem}}

\newtheorem{theorem}{Theorem}[chapter]
\newtheorem{proposition}[theorem]{Proposition}
\newtheorem{proposition-definition}[theorem]
{Proposition-Definition}
\newtheorem{corollary}[theorem]{Corollary}
\newtheorem{lemma}[theorem]{Lemma}

\theoremstyle{definition}
\newtheorem{definition}[theorem]{Definition}

\newtheorem{example}[theorem]{Example}

\newtheorem{remark}[theorem]{Remark}

\newtheorem{namedtheoremr}[theorem]{\theoremnamer}
\newcommand{\theoremnamer}{testing}
\newenvironment{namedr}[1]{\renewcommand\theoremnamer{#1}
\begin{namedtheoremr}}
{\end{namedtheoremr}}

\theoremstyle{remark}

\newcounter{steps}
\newcommand\steps{\setcounter{steps}0}
\newcommand\step[1][\relax]{\addtocounter{steps}1
\subsubsection*{Step \thesteps #1}}

\newcommand\nome{testing}
\newcommand\call[1]{\label{#1}\renewcommand\nome{#1}}
\newcommand\itemref[1]{\item\label{\nome;#1}}
\newcommand\refall[2]{\ref{#1}~(\ref{#1;#2})}
\newcommand\refpart[2]{(\ref{#1;#2})}

\newcommand\cA{\mathcal{A}}
\newcommand\cB{\mathcal{B}}
\newcommand\cC{\mathcal{C}}
\newcommand\cD{\mathcal{D}}
\newcommand\cE{\mathcal{E}}
\newcommand\cF{\mathcal{F}}
\newcommand\cG{\mathcal{G}}

\newcommand\cL{\mathcal{L}}
\newcommand\cM{\mathcal{M}}
\newcommand\cN{\mathcal{N}}
\newcommand\cO{\mathcal{O}}
\newcommand\cP{\mathcal{P}}

\newcommand\cS{\mathcal{S}}
\newcommand\cT{\mathcal{T}}
\newcommand\cU{\mathcal{U}}
\newcommand\cV{\mathcal{V}}
\newcommand\cW{\mathcal{W}}

\renewcommand\AA{\mathbb{A}}

\newcommand\GG{\mathbb{G}}

\newcommand\PP{\mathbb{P}}

\newcommand\ZZ{\mathbb{Z}}

\newcommand\rC{\mathrm{C}}

\newcommand\rF{\mathrm{F}}

\newcommand\rK{\mathrm{K}}

\newcommand\rO{\mathrm{O}}
\newcommand\rP{\mathrm{P}}

\newcommand\rma{\mathrm{a}}

\newcommand\arr{\ifinner\to\else\longrightarrow\fi}

\newcommand\xarr{\xrightarrow}

\newcommand\into{\hookrightarrow}

\newcommand\larr{\longrightarrow}

\renewcommand\H{\operatorname{H}}

\newcommand\op{^{\mathrm{op}}}

\newcommand\eqdef{\overset{\mathrm{\scriptscriptstyle def}} =}

\newcommand\im{\operatorname{im}}

\renewcommand\th{^\text{th}}

\def\displaytimes_#1{\mathrel{\mathop{\times}\limits_{#1}}}

\def\displayotimes_#1{\mathrel{\mathop{\bigotimes}\limits_{#1}}}

\renewcommand\hom{\operatorname{Hom}}

\newcommand\curshom
{\mathop{\mathcal{H}\hskip -.8pt om}\nolimits}

\newcommand\End{\operatorname{End}}

\newcommand\aut{\operatorname{Aut}}

\newcommand\spec{\operatorname{Spec}}

\newcommand{\proj}{\operatorname{Proj}}

\newcommand\generate[1]{\langle #1 \rangle}

\newcommand{\GL}{\operatorname{GL}}

\newcommand{\SL}{\operatorname{SL}}

\newcommand\M{\mathrm{M}}

\newcommand\id{\mathrm{id}}

\newcommand\pr{\operatorname{pr}}

\newcommand\indlim{\varinjlim}

\newlength{\ignora}
\newcommand{\hsmash}[1]{\settowidth{\ignora}{#1}#1\hspace{-\ignora}}

\newcommand\dash{\nobreakdash-\hspace{0pt}}

\newdir{ >}{{}*!/-5pt/@{>}}

\newcommand\cat[1]{(\mathrm{#1})}

\newcommand\catset{\cat{Set}}

\newcommand\cattop{\cat{Top}}

\newcommand\catover[2]{({\mathrm{#1}/#2})}

\newcommand\catsch[1]{\catover{Sch}{#1}}

\newcommand\catgrpover[1]{\mathrm{Grp}(#1)}

\newcommand\aff{\operatorname{Aff}}

\newcommand\cataff[1]{\catover{Aff}{#1}}

\newcommand\catsh[1]{\catover{Sh}{#1}}

\newcommand\qcoh[1]{\operatorname{QCoh}(#1)}

\newcommand\catmod[1]{\mathrm{Mod}_{#1}}

\newcommand\catqc[1]{\catover{QCoh}{#1}}

\newcommand\qcohalg{\operatorname{QCohComm}}

\newcommand\catqcalg[1]{\catover{QCohComm}{#1}}

\newcommand\catgrp{\cat{Grp}}

\newcommand\sh{\operatorname{Sh}}

\newcommand\h{\mathrm{h}}

\newcommand\func[1]{\mathrm{Hom}\bigl(#1\op, \catset\bigr)}

\newcommand\ass{^\mathrm{a}}

\newcommand\sep{^\mathrm{s}}

\newcommand\cursspec{\mathop{\mathcal{S}pec}}

\newcommand\sym{\operatorname{Sym}}

\newcommand\p{\mathrm{p}}

\newcommand\iden{\mathrm{e}}

\newcommand\mul{\mathrm{m}}

\newcommand\inv{\mathrm{i}}

\newcommand\cl{_\mathrm{cl}}

\newcommand\pt{\mathrm{pt}}

\newcommand\cart{_\mathrm{cart}}

\newcommand{\underhom}
{\mathop{\underline{\mathrm{Hom}}}\nolimits}

\newcommand{\underaut}
{\mathop{\underline{\mathrm{Aut}}}\nolimits}

\newcommand\double{\mathbin{\rightrightarrows}}

\newcommand\doublelong[2]{\mathbin{\xymatrix{{}\ar@<3pt>[r]^{#1}
\ar@<-3pt>[r]_{#2}&}}}

\newcommand\isoiden{\epsilon}

\newcommand\isoass{\alpha}

\newcommand\Arr{\operatorname{Arr}}

\newcommand\gm[1][S]{\mathbb{G}_{\mathrm{m}, #1}}

\newcommand\qc{quasi-coherent\xspace}

\newcommand\comp{_{\mathrm{comp}}}

\newcommand\desc{_{\mathrm{desc}}}

\begin{document}

\title{Notes on Grothendieck topologies,\\
fibered categories\\
and descent theory\\\ \\{\rm\small Version of May 17, 2007}}

\author[Angelo Vistoli]{Angelo Vistoli}

\address{Scuola Normale Superiore\\Piazza dei Cavalieri 7\\
56126\\Pisa\\Italy}
\email{angelo.vistoli@sns.it}


\maketitle

\tableofcontents

\chapter*{Introduction}

Descent theory has a somewhat formidable reputation among algebraic geometers. In fact, it simply says that under certain conditions homomorphisms between \qc sheaves can be constructed locally and then glued together if they satisfy a compatibility condition, while \qc sheaves themselves can be constructed locally and then glued together via isomorphisms that satisfy a cocycle condition.

Of course, if ``locally'' were to mean ``locally in the Zariski topology'' this would be a formal statement, certainly useful, but hardly deserving the name of a theory. The point is that ``locally'' here means locally in the flat topology; and the flat topology is something that is not a topology, but what is called a \emph{Grothendieck topology}. Here the coverings are, essentially, flat surjective maps satisfying a finiteness condition. So there are many more coverings in this topology than in the Zariski topology, and the proof becomes highly nontrivial.

Still, the statement is very simple and natural, provided that one resorts to the usual abuse of identifying the pullback $(gf)^{*}F$ of a sheaf $F$ along the composite of two maps $f$ and $g$ with $f^{*}g^{*}F$. If one wants to be fully rigorous, then one has to take into account the fact that $(gf)^{*}F$ and $f^{*}g^{*}F$ are not identical, but there is a canonical isomorphism between them, satisfying some compatibility conditions, and has to develop a theory of such compatibilities. The resulting complications are, in my opinion, the origin of the distaste with which many algebraic geometers look at descent theory (when they look at all).

There is also an abstract notion of ``category in which descent theory works''; the category of pairs consisting of a scheme and a \qc sheaf on it is an example. These categories are known as \emph{stacks}. The general formalism is quite useful, even outside of moduli theory, where the theory of algebraic stacks has become absolutely central (see for example \cite{deligne-mumford}, \cite{artin74} and \cite{laumon-moretbailly}).

These notes were born to accompany my ten lectures on \emph{Grothendieck topologies and descent theory} in the \emph{Advanced School in Basic Algebraic Geometry} that took place at I.C.T.P., 7--18 July 2003. They form the first part of the book \emph{Fundamental Algebraic Geometry: Grothendieck's FGA Explained}, by Barbara Fantechi (SISSA), Lothar G\"ottsche (ICTP), Luc Illusie (Universit\'e Paris-Sud), Steven L. Kleiman (MIT), Nitin Nitsure (Tata Institute of Fundamental Research), and Angelo Vistoli (Universit\`a di Bologna), published by A.M.S..\footnote{The online version, posted at the address \url{http://homepage.sns.it/vistoli/descent.pdf}, will continue to evolve; at the very least, I will correct the errors that are pointed out to me. I hope that it will also grow with the addition of new material.}

Their purpose is to provide an exposition of descent theory, more complete than the original (still very readable, and highly recommended) article of Grothendieck (\cite{grothendieck-descent}), or than \cite{sga1}. I also use the language of Grothendieck topologies, which is the natural one in this context, but had not been introduced at the time when the two standard sources were written. 

The treatment here is slanted toward the general theory of fibered categories and stacks: so the algebraic geometer searching for immediate gratification will probably be frustrated. On the other hand, I find the general theory both interesting and applicable, and hope that at at least some of my readers will agree.

Also, in the discussion of descent theory for \qc sheaves and for schemes, which forms the real reason of being of these notes, I never use the convention of identifying objects when there is a canonical isomorphism between them, but  I always specify the isomorphism, and write down explicitly the necessary compatibility conditions. This makes the treatment rigorous, but also rather heavy (for a particularly unpleasant example, see \S\ref{subsec:descent-ample}). One may question the wisdom of this choice; but I wanted to convince myself that a fully rigorous treatment was indeed possible. And the unhappy reader may be assured that this has cost more suffering to me than to her.

All of the ideas and the results contained in these notes are due to Grothendieck. There is nothing in here that is not, in some form, either in \cite{sga1} or in \cite{sga4}, so I do not claim any originality at all.

There are modern developments of descent theory, particularly in category theory (see for example \cite{joyal-tierney}) and in non\dash commutative algebra and non\dash commutative geometry (\cite{kontsevich-rosenberg-descent} and \cite{kontsevich-rosenberg-stacks}). One of the most exciting ones, for topologists as well as algebraic geometers, is the idea of ``higher descent'', strictly linked with the important topic of higher category theory (see for example \cite{hirschowitz-simpson} and \cite{street-descent}). We will not discuss any of these very interesting subjects.

\subsection*{Contents}

In Chapter~\ref{ch:recall} I recall some basic notions in algebraic geometry and category theory.

The real action starts in Chapter~\ref{ch:functors}. Here first I discuss Grothendieck's philosophy of representable functors, and give one of the main illustrative examples, by showing how this makes the notion of group scheme, and action of a group scheme on a scheme, very natural and easy. All of algebraic geometry can be systematically developed from this point of view, making it very clean and beautiful, and incomprehensible for the beginner (see \cite{demazure-gabriel}).

In Section~\ref{sec:sheaves} I define and discuss Grothendieck topologies and sheaves on them. I use the naive point of view of pretopologies, which I find much more intuitive. However, the more sophisticated point of view using sieves has advantages, so I try to have my cake and eat it too (the Italian expression, more vivid, is ``have my barrel full and my wife drunk'') by defining sieves and characterizing sheaves in terms of them, thus showing, implicitly, that the sheaves only depend on the topology and not on the pretopology. In this section I also introduce the four main topologies on the category of schemes, Zariski, \'etale, fppf and fpqc, and prove Grothendieck's theorem that a representable functor is a sheaf in all of them.

There are two possible formal setups for descent theory, fibered categories and pseudo-functors. The first one seems less cumbersome, so Chapter~\ref{ch:fibered} is dedicated to the theory of fibered categories. However, here I also define pseudo-functors, and relate the two points of view, because several examples, for example \qc sheaves, are more naturally expressed in this language. I prove some important results (foremost is Yoneda's lemma for fibered categories), and conclude with a discussion of equivariant objects in a fibered category (I hope that some of the readers will find that this throws light on the rather complicated notion of equivariant sheaf).

The heart of these notes is Chapter~\ref{ch:stacks}. After a thorough discussion of descent data (I give several definitions of them, and prove their equivalence) I define the central concept, that of \emph{stack}: a stack is a fibered category over a category with a Grothendieck topology, in which descent theory works (thus we see all the three notions appearing in the title in action). Then I proceed to proving the main theorem, stating that the fibered category of \qc sheaves is a stack in the fpqc topology. This is then applied to two of the main examples where descent theory for schemes works, that of affine morphisms, and morphisms endowed with a canonical ample line bundle. I also discuss a particularly interesting example, that of descent along principal bundles (torsors, in Grothendieck's terminology).

In the last section I give an example to show that \'etale descent does not always work for schemes, and end by mentioning that there is an extension of the concept of scheme, that of \emph{algebraic space}, due to Michael Artin. Its usefulness is that on one hand algebraic spaces are, in a sense, very close to schemes, and one can define for them most of the concepts of scheme theory, and on the other hand fppf descent always works for them. It would have been a natural topic to include in the notes, but this would have further delayed their completion.

\subsection*{Prerequisites}

I assume that the reader is acquainted with the language of schemes, at least at the level of Hartshorne's book (\cite{hartshorne}). I use some concepts that are not contained in \cite{hartshorne}, such as that of a morphism locally of finite presentation; but I recall their main properties, with references to the appropriate parts of \emph{\'{E}l\'ements de g\'eom\'etrie alg\'ebrique}, in Chapter~\ref{ch:recall}. 

I make heavy use of the categorical language: I assume that the reader is acquainted with the notions of category, functor and natural transformation, equivalence of categories. On the other hand, I do not use any advanced concepts, nor do I use any real results in category theory, with one exception: the reader should know that a fully faithful essentially surjective functor is an equivalence.

\subsection*{Acknowledgments} Teaching my course at the \emph{Advanced School in Basic Algebraic Geometry} has been a very pleasant experience, thanks to the camaraderie of my fellow lecturers (Lothar G\"ottsche, Luc Illusie, Steve Kleiman and Nitin Nitsure) and the positive and enthusiastic attitude of the participants. I am also in debt with Lothar, Luc, Steve and Nitin because they never once complained about the delay with which these notes were being produced.

I am grateful to Steve Kleiman for useful discussions and suggestions, particularly involving the fpqc topology, and to Pino Rosolini, who, during several hikes on the Alps, tried to enlighten me on some categorical constructions.

 I have had some interesting conversations with Behrang Noohi concerning the definition of a stack: I thank him warmly.
 
I learned about the counterexample in \cite[XII 3.2]{raynaudample} from Andrew Kresch.

I also thank the many participants to the school who showed interest in my lecture series, and particularly those who pointed out mistakes in the first version of the notes. I am especially in debt with Zoran Skoda, who sent me several helpful comments, and also for his help with the bibliography.

Joachim Kock read carefully most of this, and sent me a long list of comments and corrections, which were very useful. More corrections where provided by the referees, by Ms. Elaine Becker, from the A.M.S., by Luigi Previdi and by Henning Ulfarsson, who also found mistakes in the wording of Definition~\ref{def:saturation} and in the statement of Proposition~\ref{prop:characterization-quasifunctors}. Alon Shapira found a serious error in the proof of Lemma~\ref{lem:criterion-sheaf}. I am grateful to them.

Finally, I would like to dedicate these notes to the memory of my father-in-law, Amleto Rosolini, who passed away at the age of 86 as they were being completed. He would not have been interested in learning descent theory, but he was a kind and remarkable man, and his enthusiasm about mathematics, which lasted until his very last day, will always be an inspiration to me.

\chapter{Preliminary notions} \label{ch:recall}

\section{Algebraic geometry}

In this chapter we recall, without proof, some basic notions of scheme theory that are used in the notes. All rings and algebras will be commutative.

We will follow the terminology of \emph{\'{E}l\'ements de g\'eom\'etrie alg\'ebrique}, with the customary exception of calling a ``scheme'' what is called there a ``prescheme'' (in \emph{\'{E}l\'ements de g\'eom\'etrie alg\'ebrique}, a scheme is assumed to be separated).

We start with some finiteness conditions. Recall if $B$ is an algebra over the ring $A$, we say that $B$ is \emph{finitely presented} if it is the quotient of a polynomial ring $A[x_1, \dots, x_n]$ over $A$ by a finitely generated ideal. If $A$ is noetherian, every finitely generated algebra is finitely presented.

If $B$ is finitely presented over $A$, whenever we write $B = A[x_1, \dots, x_n]/I$, $I$ is always finitely generated in $A[x_1, \dots, x_n]$ (\cite[Proposition~1.4.4]{ega4-1}).

\begin{definition}[See \mbox{\cite[1.4.2]{ega4-1}}]
A morphism of schemes $f \colon X \arr Y$ is \emph{locally of finite presentation}
\index{morphism of schemes!locally of finite presentation}
\index{locally of finite presentation!morphism of schemes}
if for any $x \in X$ there are affine neighborhoods $U$ of $x$ in $X$ and $V$ of $f(x)$ in $Y$ such that $f(U) \subseteq V$ and $\cO(U)$ is finitely presented over $\cO(V)$.
\end{definition}

Clearly, if $Y$ is locally noetherian, then $f$ is locally of finite presentation if and only if it is locally of finite type.

\begin{proposition}[\mbox{\cite[1.4]{ega4-1}}]\hfil
\begin{enumeratei}

\item If $f \colon X \arr Y$ is locally of finite presentation, $U$ and $V$ are open affine subsets of $X$ and $Y$ respectively, and $f(U) \subseteq V$, then $\cO(U)$ is finitely presented over $\cO(V)$.

\item The composite of morphisms locally of finite presentation is locally of finite presentation.

\item Given a cartesian diagram
   \[
   \xymatrix{
   X'\ar[r]\ar[d] & X\ar[d] \\
   Y' \ar[r] & Y
   }
   \]
if $X \arr Y$ is locally of finite presentation, so is $X' \arr Y'$.
\end{enumeratei}
\end{proposition}

\begin{definition}[See \mbox{\cite[6.6.1]{ega1}}]
A morphism of schemes $X \arr Y$ is \emph{quasi-compact}\index{morphism of schemes!quasi-compact}\index{quasi-compact morphism of schemes} if the inverse image in $X$ of a quasi-compact open subset of $Y$ is quasi-compact.

An affine scheme is quasi-compact, hence a scheme is quasi-compact if and only if it is the finite union of open affine subschemes; using this, it is easy to prove the following.

\end{definition}

\begin{proposition}[\mbox{\cite[Proposition~6.6.4]{ega1}}]\hfil

Let $f \colon X \arr Y$ be a morphism of schemes. The following are equivalent.

\begin{enumeratei}

\item $f$ is quasi-compact.

\item The inverse image of an open affine subscheme of $Y$ is quasi-compact.

\item There exists a covering $Y = \cup_i V_i$ by open affine subschemes, such that the inverse image in $X$ of each $V_i$ is quasi-compact.
\end{enumeratei}

In particular, a morphism from a quasi-compact scheme to an affine scheme is quasi-compact.
\end{proposition}

\begin{remark}
It is not enough to suppose that there is a covering of $Y$ by open quasi-compact subschemes $V_i$, such that the inverse image of each $V_i$ is quasi-compact in $X$, without additional hypotheses. For example, consider a ring $A$ that does not satisfy the ascending chain condition on radical ideals (for example, a polynomial ring in infinitely many variables), and set $X = \spec A$. In $X$ there will be an open subset $U$ that is not quasi-compact; denote by $Y$ the scheme obtained by gluing two copies of $X$ together along $U$, and by $f \colon X \arr Y$ the inclusion of one of the copies. Then $Y$ and $X$ are both quasi-compact; on the other hand there is an affine open subset of $Y$ (the other copy of $X$) whose inverse image in $X$ is $U$, so $f$ is not quasi-compact.
\end{remark}

\begin{proposition}[\mbox{\cite[6.6]{ega1}}]\hfil

\begin{enumeratei}

\item The composite of quasi-compact morphisms is quasi-compact.

\item Given a cartesian diagram
   \[
   \xymatrix{
   X'\ar[r]\ar[d] & X\ar[d] \\
   Y' \ar[r] & Y
   }
   \]
if $X \arr Y$ is quasi-compact, so is $X' \arr Y'$.
\end{enumeratei}

\end{proposition}

Let us turn to flat morphisms. 

\begin{definition}
A morphism of schemes $f \colon X \arr Y$ is \emph{flat}\index{morphism of schemes!flat}\index{flat morphism of schemes} if for any $x \in X$, the local ring $\cO_{X,x}$ is flat as a module over
$\cO_{Y, f(x)}$.
\end{definition}

\begin{proposition}[\mbox{\cite[Proposition~2.1.2]{ega4-2}}]
Let $f \colon X \arr Y$ be a morphism of schemes. Then the following are equivalent.
\begin{enumeratei}

\item $f$ is flat.

\item For any $x \in X$, there are affine neighborhoods $U$ of $x$ in $X$ and $V$ of $f(x)$ in $Y$ such that $f(U) \subseteq V$, and $\cO(U)$ is flat over $\cO(V)$.

\item For any open affine subsets $U$ in $X$ and $V$ in $Y$ such that $f(U) \subseteq V$, $\cO(U)$ is flat over $\cO(V)$.

\end{enumeratei}
\end{proposition}

\begin{proposition}[\mbox{\cite[2.1]{ega4-2}}]\hfil

\begin{enumeratei}

\item The composite of flat morphisms is flat.

\item Given a cartesian diagram
   \[
   \xymatrix{
   X'\ar[r]\ar[d] & X\ar[d] \\
   Y' \ar[r] & Y
   }
   \]
if $X \arr Y$ is flat, so is $X' \arr Y'$.
\end{enumeratei}
\end{proposition}

\begin{definition}
A morphism of schemes $f \colon X \arr Y$ is \emph{faithfully flat}\index{morphism of schemes!faithfully flat}\index{faithfully flat morphism of schemes} if it is flat and surjective.
\end{definition}

Let $B$ be an algebra over $A$. We say that $B$ is faithfully flat if the associated morphism of schemes $\spec B \arr \spec A$ is faithfully flat.

\begin{proposition}[\mbox{\cite[Theorems 7.2 and 7.3]{matsumura89}}]

Let $B$ be an algebra over $A$. The following are equivalent.

\begin{enumeratei}

\item $B$ is faithfully flat over $A$.

\item A sequence of $A$-modules $M' \arr M \arr M''$ is exact if and only if the induced sequence of $B$-modules $M' \otimes_A B \arr M \otimes_A B \arr M'' \otimes_A B$  is exact.

\item A homomorphism of $A$-modules $M' \arr M$ is injective if and only if the associated homomorphism of $B$-modules $M' \otimes_A B \arr M \otimes_A B$ is injective.

\item $B$ is flat over $A$, and if $M$ is a module over $A$ with $M
\otimes_A B = 0$, we have $M = 0$.

\item $B$ is flat over $A$, and $\mathfrak{m}B \neq B$ for all maximal ideals $\mathfrak{m}$ of $A$. 

\end{enumeratei}
\end{proposition}

The following fact is very important.

\begin{proposition}[\mbox{\cite[Proposition~2.4.6]{ega4-2}}]
\label{prop:flat->open}
A flat morphism that is locally of finite presentation is open.
\end{proposition}

This is not true in general for flat morphisms that are not locally of finite presentation; however, a weaker version of this fact holds.

\begin{proposition}[\mbox{\cite[Corollaire~2.3.12]{ega4-2}}]
\label{prop:flat->quotient-topology}
If $f \colon X \arr Y$ is a faithfully flat quasi-compact morphism, a subset of $Y$ is open if and only if its inverse image in $X$ is open in $X$.
\end{proposition}

In other words, $Y$ has the topology induced by that of $X$.

\begin{remark}\label{rmk:need-finiteness}
For this we need to assume that $f$ is quasi-compact, it is not enough to assume that it is faithfully flat. For example, let $Y$ be an integral smooth curve over an algebraically closed field, $X$ the disjoint union of the $\spec\cO_{Y,y}$ over all closed points $y \in Y$. The natural projection $f \colon X \arr Y$ is clearly flat. However, if $S$ is a subset of $Y$ containing the generic point, then $f^{-1}S$ is always open in $X$, while $S$ is open in $Y$ if and only if its complement is finite.
\end{remark}

\begin{proposition}[\mbox{\cite[Proposition~2.7.1]{ega4-2}}]
\label{prop:local-qcflat}
Let
   \[
   \xymatrix{
   X' \ar[r] \ar[d] &
   X \ar[d] \\
   Y' \ar[r] &
   Y
   }
   \]
be a cartesian diagram of schemes in which $Y' \arr Y$ is faithfully flat and either quasi-compact or locally of finite presentation. Suppose that $X' \arr Y'$ has one of the following properties:
\begin{enumeratei}

\item is separated,

\item is quasi-compact,

\item is locally of finite presentation,

\item is proper,

\item is affine,

\item is finite,

\item is flat,

\item is smooth, 

\item is unramified,

\item is \'etale,

\item is an embedding,

\item is a closed embedding.

\end{enumeratei}

Then $X \arr Y$ has the same property.
\end{proposition}

In \cite{ega4-2} all these statements are proved when $Y' \arr Y$ is quasi-compact. Using Proposition~\ref{prop:flat->open}, and the fact that all those properties are local in the Zariski topology of $Y$, it is not hard to prove the statement also when $Y' \arr Y$ is locally of finite presentation.

\section{Category theory}

We will assume that the reader is familiar with the concepts of category, functor and natural transformation. The standard reference in category theory, containing a lot more than what we need, is \cite{maclane98}; also very useful are \cite{borceux1}, \cite{borceux2} and \cite{borceux3}.

We will not distinguish between small and large categories. More generally, we will ignore any set-theoretic difficulties. These can be overcome with standard arguments using universes.

If $F \colon \cA \arr \cB$ is a functor, recall that $F$ is called \emph{fully faithful}\index{functor!fully faithful}\index{fully faithful functor} when for any two objects $A$ and $A'$ of $\cA$, the function
   \[
   \hom_{\cA}(A, A') \arr \hom_{\cB}(FA, FA')
   \]
defined by $F$ is a bijection. $F$ is called \emph{essentially surjective}\index{functor!essentially surjective}\index{essentially surjective functor} if every object of $\cB$ is isomorphic to the image of an object of $\cA$.

Recall also that $F$ is called an \emph{equivalence}\index{equivalence!of categories} when there exists a functor $G \colon \cB \arr \cA$, such that the composite $GF \colon \cA \arr \cA$ is isomorphic to $\id_{\cA}$, and $FG \colon \cB \arr \cB$ is isomorphic to $\id_{\cB}$.

The composite of two equivalences is again an equivalence. In particular, ``being equivalent'' is an equivalence relation among categories.

The following well-known fact will be used very frequently.

\begin{proposition}
A functor is an equivalence if and only if it is both fully faithful and essentially surjective.
\end{proposition}

If $\cA$ and $\cB$ are categories, there is a category $\hom(\cA, \cB)$, whose objects are functors $\Phi \colon \cA \arr \cB$, and whose arrows $\alpha \colon \Phi \arr \Psi$ are natural transformations. If $F \colon \cA' \arr \cA$ is a functor, there is an induced functor
   \[
   F^{*}\colon \hom(\cA, \cB) \arr \hom(\cA', \cB)
   \]
defined at the level of objects by the obvious rule
   \[
   F^{*}\Phi \eqdef \Phi \circ F \colon \cA' \arr \cB
   \]
for any functor $\Phi\colon \cA \arr \cB$. At the level of arrows $F^{*}$ is defined by the formula
   \[
   (F^{*}\alpha)_{A'} \eqdef \alpha_{FA'} \colon \Phi FA' \arr \Psi FA'
   \]
for any natural transformation $\alpha \colon \Phi \arr \Psi$.

Also for any functor $F \colon \cB \arr \cB'$ we get an induced functor
   \[
   F_{*}\colon \hom(\cA, \cB) \arr \hom(\cA, \cB')
   \]
obtained by the obvious variants of the definitions above.

The  reader should also recall that a \emph{groupoid}\index{groupoid} is a category in which every arrow is invertible.

We will also make considerable use of the notions of fibered product and cartesian diagram in an arbitrary category. We will manipulate some cartesian diagrams. In particular the reader will encounter diagrams of the type
   \[
   \xymatrix{
   A'\ar[r]\ar[d] & B' \ar[r]\ar[d] & C'\ar[d]\\
   A\ar[r] & B \ar[r] & C\hsmash{;}
   }
   \]
we will say that this is cartesian when both squares are cartesian. This is equivalent to saying that the right hand square and the square
   \[
   \xymatrix{
   A'\ar[r]\ar[d] & C'\ar[d]\\
   A\ar[r] & C\hsmash{,}
   }
   \]
obtained by composing the rows, are cartesian. There will be other statements of the type ``there is a cartesian diagram \dots''. These should all be straightforward to check.

For any category $\cC$ and any object $X$ of $\cC$ we denote by $(\cC/X)$ the \emph{comma category}\index{category!comma}\index{comma!category}, whose objects are arrows $U \arr X$ in $\cC$, and whose arrows are commutative diagrams
   \[
   \xymatrix@-15pt{
   U\ar[rr]\ar[rd] &&V\ar[ld]\\
			&X
   }
   \]

We also denote by $\cC\op$ the \emph{opposite category}\index{category!opposite} \index{opposite category} of $\cC$, in which the objects are the same, and the arrows are also the same, but sources and targets are switched. A \emph{contravariant functor}\index{functor!contravariant}\index{contravariant functor} from $\cC$ to another category $\cD$ is a functor $\cC\op \arr \cD$.

Whenever we have a fibered product $X_1 \times_Y X_2$ in a category, we denote by $\pr_1 \colon X_1 \times_Y X_2 \arr X_1$ and $\pr_2 \colon X_1 \times_Y X_2 \arr X_2$ the two projections. We will also use a similar notation for the product of three or more objects: for example, we will denote by
   \[
   \pr_{i} \colon X_1 \times_Y X_2 \times_{Y} X_{3} \arr X_i
   \]
the $i\th$ projection, and by
   \[
   \pr_{ij} \colon X_1 \times_Y X_2 \times_{Y} X_{3} \arr X_i \times_{Y} X_{j}
   \]
the projection into the product of the $i\th$ and  $j\th$ factor.

Recall that a category has finite products if and only if it has a terminal object (the product of $0$ objects) and products of two objects.

Suppose that $\cC$ and $\cD$ are categories with finite products; denote their terminal objects by $\pt_{\cC}$ and $\pt_{\cD}$. A functor $F \colon \cC \arr \cD$ is said to \emph{preserve finite products}\index{functor!preserving finite products} if the following holds. Suppose that we have objects $U_1$, \dots,~$U_r$ of $\cC$: the projections $U_1 \times \dots \times U_r \arr U_i$ induce arrows $F(U_1 \times \dots \times U_r) \arr F U_i$. Then the corresponding arrow
   \[
   F(U_1 \times \dots \times U_r) \arr F U_1 \times \dots \times F U_r
   \]
is an isomorphism in $\cC$.

If $f_1 \colon U_1 \arr V_1$, \dots,~$f_r \colon U_r \arr V_r$ are arrows in $\cC$, the diagram
   \[
   \xymatrix{
   F(U_1 \times \dots \times U_r) \ar[r]
           \ar[d]^{F(f_1\times \dots \times f_r)} &
   F U_1 \times \dots \times F U_r \ar[d]^{Ff_1 \times \dots \times Ff_r} \\
   F(V_1 \times \dots \times V_r) \ar[r] &
   F V_1 \times \dots \times F V_r
   }
   \]
in which the horizontal arrows are the isomorphism defined above, obviously commutes.

By a simple induction argument, $F$ preserves finite products if and only if $F\pt_{\cC}$ is a terminal object of $\cD$, and for any two objects $U$ and $V$ of $\cC$ the arrow $F(U \times V) \arr FU \times FV$ is an isomorphism (in other words, for $F$ to preserve finite products it is enough that it preserves products of $0$ and $2$ objects).

Finally, we denote by $\catset$ the category of sets, by $\cattop$ the category of topological spaces, $\catgrp$ the category of groups, and by $\catsch{S}$ the category of schemes over a fixed base scheme $S$.

\chapter{Contravariant functors}\label{ch:functors}

\section{Representable functors and the Yoneda Lemma}\label{sec:repfunctors}

\subsection{Representable functors}
Let us start by recalling a few basic notions of category theory.

Let $\cC$ be a category. Consider functors from $\cC\op$ to $\catset$. These are the objects of a category, denoted by
   \[
   \func{\cC},
   \]
in which the arrows are the natural transformations. From now on we will refer to natural transformations of contravariant functors on $\cC$ as \emph{morphisms}.

Let $X$ be an object of $\cC$. There is a functor\index{$\h_{X}$}
   \[
   \h_X\colon \cC\op \arr\catset
   \]
to the category of sets, which sends an object $U$ of $\cC$ to the set
   \[
   \h_X U = \hom_\cC(U, X).
   \] If $\alpha
\colon U' \arr U$ is an arrow in $\cC$, then $\h_X \alpha  \colon \h_X U \arr \h_X U'$ is defined to be composition with $\alpha$. (When $\cC$ is the category of schemes over a fixed base scheme, $\h_{X}$ is often called the \emph{functor of points of $X$}\index{functor!of points})

Now, an arrow $f \colon X \arr Y$ yields a function $\h_f U \colon \h_X U \arr \h_Y U$\index{$\h_{f}$} for each object $U$ of $\cC$, obtained by composition with $f$. This defines a morphism $\h_X \arr \h_Y$, that is, for all arrows $\alpha\colon U' \arr U$ the diagram
   \[
   \xymatrix{
   {}\h_X U\ar[r]^{\h_f U} \ar[d]^{\h_X \alpha}
   & {}\h_Y U \ar[d]^{\h_Y \alpha}\\
   {}\h_X U'\ar[r]^{\h_f U'}&{}\h_Y U'\\
   }
   \]
commutes.

Sending each object $X$ of $\cC$ to $\h_X$, and each arrow $f \colon X \arr Y$ of $\cC$ to $\h_f \colon \h_X \arr \h_Y$ defines a functor $\cC \arr \func{\cC}$.

\begin{named}{Yoneda Lemma (weak version)}
\index{Yoneda Lemma!weak version}
\index{Lemma!Yoneda!weak version}
Let $X$ and $Y$ be
objects of $\cC$. The function
   \[
   \hom_\cC(X, Y) \arr
   \hom (\h_X, \h_Y)
   \]
that sends $f \colon X \arr Y$ to $\h_f \colon \h_X \arr \h_Y$ is bijective. 
\end{named}

In other words, the functor $\cC \arr \func{\cC}$ is fully faithful. It fails to be an equivalence of categories, because in general it will not be essentially surjective. This means that not every functor $\cC\op \arr \catset$ is isomorphic to a functor of the form $\h_X$. However, if we restrict to the full subcategory of $\func{\cC}$ consisting of functors $\cC\op \arr \catset$ which are isomorphic to a functor of the form $\h_X$, we do get a category which is equivalent to $\cC$.

\begin{definition} A \emph{representable functor}\index{functor!representable}\index{representable!functor} on the category $\cC$ is a functor
   \[
   F \colon \cC\op \arr \catset
   \] which is isomorphic to a functor of the form $\h_X$ for some object $X$ of $\cC$. 

If this happens, we say that \emph{$F$ is represented by $X$}.
\end{definition}

Given two isomorphisms $F \simeq \h_X$ and $F \simeq \h_Y$, we have that the resulting isomorphism $\h_X \simeq \h_Y$ comes from a unique isomorphism $X \simeq Y$ in $\cC$, because of the weak form of Yoneda's lemma. Hence \emph{two objects representing the same functor are canonically isomorphic.}
 
\subsection{Yoneda's lemma}

The condition that a functor be representable can be given a new expression with the more general version of Yoneda's lemma. Let $X$ be an object of $\cC$ and $F \colon \cC\op \arr \catset$ a functor. Given a natural transformation $\tau\colon \h_X \arr F$, one gets an element $\xi \in F X$, defined as the image of the identity map $\id_X \in \h_X X$ via the function $\tau_X\colon \h_X X \arr F X$. This construction defines a function $\hom(\h_X, F) \arr F X$.
 
Conversely, given an element $\xi\in F X$, one can define a morphism $\tau\colon \h_X \arr F$ as follows. Given an object $U$ of $\cC$, an element of $\h_X U$ is an arrow $f\colon U \arr X$; this arrow induces a function $F f \colon F X \arr F U$. We define a function $\tau_{U} \colon \h_X U \arr F U$ by sending $f \in \h_X U$ to $F f(\xi) \in F U$. It is straightforward to check that the $\tau$ that we have defined is in fact a morphism. In this way we have defined functions
   \[
   \hom(\h_X,F) \arr F X
   \]
and
   \[
   F X \arr \hom(\h_X,F).
   \]

\begin{named}{Yoneda lemma}
\index{Yoneda Lemma!strong version}
\index{Lemma!Yoneda!strong version}
These two functions are inverse to each other, and therefore establish a bijective correspondence
   \[
   \hom(\h_X,F) \simeq F X.
   \]
\end{named}

The proof is easy and left to the reader. Yoneda's lemma is not a deep fact, but its importance cannot be overestimated.

Let us see how this form of Yoneda's lemma implies the weak form above. Suppose that $F = \h_Y$: the function $\hom(X,Y) = \h_Y X \arr \hom(\h_X,\h_Y)$ constructed here sends each arrow $f \colon X \arr Y$ to
   \[
   \h_Y f(\id_Y) = \id_Y \circ f \colon X \arr Y,
   \]
so it is exactly the function $\hom(X,Y) \arr \hom(\h_X,\h_Y)$ appearing in the weak form of the result.

One way to think about Yoneda's lemma is as follows. The weak form says that the category $\cC$ is embedded in the category $\func{\cC}$. The strong version says that, given a functor $F \colon \cC\op \arr \catset$, this can be extended to the representable functor $\h_F \colon \func{\cC}\op \arr \catset$: thus, every functor becomes representable, when extended appropriately. (In practice, the functor category $\func{\cC}$ is usually much too big, and one has to restrict it appropriately.)

We can use Yoneda's lemma to give a very important characterization of representable functors.

\begin{definition} Let $F \colon \cC\op \arr \catset$ be a functor. A \emph{universal object}\index{object!universal}\index{universal!object} for $F$ is a pair $(X, \xi)$ consisting of an object $X$ of $\cC$, and an element $\xi \in F X$, with the property that for each object $U$ of $\cC$ and each $\sigma\in F U$, there is a unique arrow $f \colon U \arr X$ such that $F f(\xi) = \sigma \in F U$.
\end{definition}

In other words: the pair $(X, \xi)$ is a universal object if the morphism $\h_X \arr F$ defined by $\xi$ is an isomorphism. Since every natural transformation $\h_X \arr F$ is defined by some object $\xi \in F X$, we get the following.

\begin{proposition}
\index{representable!functor!characterization via existence of a universal object}
\index{functor!representable!characterization via existence of a universal object}
A functor $F \colon \cC\op \arr \catset$ is representable if and only if it has a universal object.
\end{proposition}

Also, if $F$ has a universal object $(X, \xi)$, then $F$ is represented by $X$.

Yoneda's lemma ensures that the natural functor $\cC \arr \func{\cC}$ which sends an object $X$ to the functor $\h_X$ is an equivalence of $\cC$ with the category of representable functors. From now on we will not distinguish between an object $X$ and the functor $\h_X$ it represents. So, if $X$ and $U$ are objects of $\cC$, we will write $X(U)$ for the set  $\h_X U = \hom_{\cC}(U,X)$ of arrows $U \arr X$. Furthermore, if $X$ is an object and $F \colon \cC\op \arr \catset$ is a functor, we will also identify the set $\hom(X, F) = \hom(\h_X, F)$ of morphisms from $\h_X$ to $F$ with $F X$.

\subsection{Examples}\label{subsec:examples}
\index{representable!functor!examples of}
\index{functor!representable!examples of}

Here are some examples of representable and non\dash represent\-a\-ble functors.

\begin{enumeratei}

\item Consider the functor $\rP \colon \catset\op \arr \catset$ that sends each set $S$ to the set $\rP(S)$ of subsets of $S$. If $f \colon S \arr T$ is a function, then $\rP(f) \colon \rP(T) \arr \rP(S)$ is defined by $\rP(f)\tau = f^{-1}\tau$ for all $\tau \subseteq T$.

Given a subset $\sigma \subseteq S$, there is a unique function $\chi_\sigma \colon S \arr \{0,1\}$ such that $\chi_\sigma^{-1}(\{1\}) = \sigma$, namely the \emph{characteristic function}, defined by 
   \[
   \chi_\sigma(s) = \begin{cases}1 & \text{if $s \in \sigma$}\\
   0 & \text{if $s \notin \sigma$.}
   \end{cases}
   \]
Hence the pair $(\{0,1\}, \{1\})$ is a universal object, and the functor $\rP$ is represented by $\{0,1\}$.

\item\label{ex:rep-opensubsets} This example is similar to the previous one. Consider the category $\cat{Top}$ of all topological spaces, with the arrows being given by continuous functions. Define a functor $\rF \colon \cat{Top}\op \arr \catset$ sending each topological space $S$ to the collection $\rF(S)$ of all its open subspaces. Endow $\{0,1\}$ with the coarsest topology in which the subset $\{1\} \subseteq \{0,1\}$ is open; the open subsets in this topology are $\emptyset$, $\{1\}$ and $\{0,1\}$. A function $S \arr \{0,1\}$ is continuous if and only if $f^{-1}(\{1\})$ is open in $S$, and so one sees that the pair $(\{0,1\}, \{1\})$ is a universal object for this functor.

The space $\{0,1\}$ is called \emph{the Sierpinski space}.

\item\label{ex:notrep-opensubsets} The next example may look similar, but the conclusion is very different. Let $\cat{HausTop}$ be the category of all Hausdorff topological spaces, and consider the restriction  $\rF\colon \cat{HausTop}\op \arr \catset$ of the functor above. I claim that this functor is not representable. 

In fact, assume that $(X, \xi)$ is a universal object. Let $S$ be any set, considered with the discrete topology; by definition, there is a unique function $f \colon S \arr X$ with $f^{-1} \xi = S$, that is, a unique function $S \arr \xi$. This means that $\xi$ can only have one element. Analogously, there is a unique function $S \arr X \setminus \xi$, so $X \setminus \xi$ also has a unique element. But this means that $X$ is a Hausdorff space with two elements, so it  must have the discrete topology; hence $\xi$ is also closed in $X$. Hence, if $S$ is any topological space with a closed subset $\sigma$ that is not open, there is no continuous function $f\colon S \arr X$ with $f^{-1}{\xi} = \sigma$.

\item Take $\cat{Grp}$ to be the category of groups, and consider the functor 
   \[
   \operatorname{Sgr} \colon \cat{Grp}\op \arr \catset
   \]
that associates with each group $G$ the set of all its subgroups. If $f \colon G \arr H$ is a group homomorphism, we take $\operatorname{Sgr} f \colon \operatorname{Sgr} H \arr \operatorname{Sgr} G$ to be the function associating with each subgroup of $H$ its inverse image in $G$.

This is not representable: there does not exist a group $\Gamma$, together with a subgroup $\Gamma_1 \subseteq \Gamma$, with the property that for all groups $G$ with a subgroup $G_1 \subseteq G$, there is a unique homomorphism $f \colon G \arr \Gamma$ such that $f^{-1} \Gamma_1 = G_1$. This can be checked in several ways; for example, if we take the subgroup $\{0\} \subseteq \ZZ$, there should be a unique homomorphism $f \colon \ZZ \arr \Gamma$ such that $f ^{-1} \Gamma_1 = \{0\}$. But given one such $f$, then the homomorphism $\ZZ \arr \Gamma$ defined by $n \mapsto f(2n)$ also has this property, and is different, so this contradicts uniqueness.

\item Here is a much more sophisticated example. Let $\cat{Hot}$ be the category of CW complexes, with the arrows being given by homotopy classes of continuous functions. If $n$ is a fixed natural number, there is a functor $\H^n \colon \cat{Hot}\op \arr \catset$ that sends a CW complex $S$ to its $n\th$ cohomology group $\H^n(S, \ZZ)$. Then it is a highly nontrivial fact that this functor is represented by a CW complex, known as a Eilenberg--Mac Lane space, usually denoted by $\rK( \ZZ, n)$.

\end{enumeratei}

But we are really interested in algebraic geometry, so let's give some examples in this context. Let $S = \spec R$ (this is only for simplicity of notation, if $S$ is not affine, nothing substantial changes).

\begin{example}\label{ex:affinespace} Consider the affine line $\AA^1_S = \spec R[x]$. We have a functor
   \[
   \cO\colon \catsch{S}\op \arr \catset
   \]
that sends each $S$-scheme $U$ to the ring of global sections $\cO(U)$. If $f\colon U \arr V$ is a morphism of schemes, the corresponding function $\cO(V) \arr \cO(U)$ is that induced by $f^{\sharp}\colon \cO_{V} \arr f_{*}\cO_{U}$.

Then $x \in \cO(\AA^1_S)$, and given a scheme $U$ over $S$, and an element $f \in \cO(U)$, there is a unique morphism $U \arr \AA^1_S$ such that the pullback of $x$ to $U$ is precisely $f$. This means that the functor $\cO$ is represented by $\AA^1_S$, and the pair $(\AA_S, x)$ is a universal object.

More generally, the affine space $\AA^n_S$ represents the functor $\cO^n$ that sends each scheme $S$ to the ring $\cO(S)^n$.
\end{example}

\begin{example}
Now we look at $\gm = \AA^1_S \setminus 0_{S} = \spec R[x, x^{-1}]$. Here by $0_S$ we mean the image of the zero-section $S \arr \AA^1_S$. Now, a morphism of $S$-schemes $U \arr \gm$ is determined by the image of $x \in \cO(\gm)$ in $\cO(S)$; therefore $\gm$ represents the functor $\cO^*\colon \catsch{S}\op \arr \catset$ that sends each scheme $U$ to the group $\cO^*(U)$ of invertible sections of the structure sheaf.
\end{example}

A much more subtle example is given by projective spaces.

\begin{example}
On the projective space $\PP^n_S = \proj R[x_0, \ldots, x_n]$ there is a line bundle $\cO(1)$, with $n+1$ sections $x_0$, \dots,~$x_n$, which generate it.

Suppose that $U$ is a scheme, and consider the set of sequences
   \[
   (\cL, s_0, \ldots, s_n),
   \]
where $\cL$ is an invertible sheaf on $U$, $s_0$, \dots,~$s_n$ sections of $\cL$ that generate it. We say that $(\cL, s_0, \ldots, s_n)$ is equivalent to $(\cL', s'_0, \ldots, s'_n)$ if there exists an isomorphism of invertible sheaves $\phi\colon \cL \simeq \cL'$ carrying each $s_i$ into $s'_i$. Notice that, since the $s_i$ generate $\cL$, if $\phi$ exists then it is unique.

One can consider a function $Q_n \colon \catsch{S} \arr \catset$ that associates with each scheme $U$ the set of sequences $(\cL, s_0, \ldots, s_n)$ as above, modulo equivalence. If $f \colon U \arr V$ is a morphism of $S$-schemes, and $(\cL, s_0, \ldots, s_n)\in Q_{n}(V)$, then there are sections $f^*s_0$, \dots,~$f^*s_n$ of $f^* \cL$ that generate it; this makes $Q_n$ into a functor $\catsch{S}\op \arr \catset$.

Another description of the functor $Q_n$ is as follows. Given a scheme $U$ and a sequence $(\cL, s_0, \ldots, s_n)$ as above, the $s_i$ define a homomorphism $\cO_U^{n+1} \arr \cL$, and the fact that the $s_i$ generate is equivalent to the fact that this homomorphism is surjective. Then two sequences are equivalent if and only if the represent the same quotient of $\cO_S^n$.

It is a very well-known fact, and, indeed, one of the cornerstones of algebraic geometry, that for any sequence $(\cL, s_0, \ldots, s_n)$ over an $S$-scheme $U$, there it exists a unique morphism $f \colon U \arr \PP^n_S$ such that $(\cL, s_0, \ldots, s_n)$ is equivalent to $(f^*\cO(1), f^*x_0, \ldots, f^*x_n)$. This means precisely that $\PP^n_S$ represents the functor $Q_n$.
\end{example}

\begin{example}
This example is an important generalization of the previous one.

Here we will let $S$ be an arbitrary scheme, not necessarily affine, $\cM$ a \qc sheaf on $S$. In Grothendieck's notation, $\pi\colon \PP(\cM) \arr S$ is the relative homogeneous spectrum $\proj_{S}\sym_{\cO_{S}} \cM$ of the symmetric sheaf of algebras of $\cM$ over $\cO_{S}$. Then on $\PP(\cM)$ there is an invertible sheaf, denoted by $\cO_{\PP(\cM)}(1)$, which is a quotient of $\pi^{*}\cM$. This is a universal object, in the sense that, given any $S$-scheme $\phi\colon U \arr S$, with an invertible sheaf $\cL$ and a surjection $\alpha\colon \phi^{*}\cM \twoheadrightarrow \cL$, there is unique morphism of $S$-schemes $f\colon U \arr \PP(\cM)$, and an isomorphism of $\cO_{U}$-modules $\alpha\colon \cL \simeq f^{*}\cO_{\PP(\cM)}(1)$, such that the composite
   \[
   f^{*}\pi^{*}\cM \simeq \phi^{*}\cM \twoheadrightarrow \cL
      \xarr{\alpha} f^{*}\cO_{\PP(\cM)}(1)
   \]
is the pullback of the projection $\pi^{*}\cM \twoheadrightarrow \cO_{\PP(\cM)}(1)$ (\cite[Proposition~4.2.3]{ega1}).

This means the following. Consider the functor $Q_{\cM} \colon \catsch{S}\op \arr \catset$ that sends each scheme $\phi \colon U \arr S$ over $S$ to the set of all invertible quotients of the pullback $\phi^*\cM$. If $f\colon V \arr U$ is a morphism of $S$-schemes from $\phi\colon U \arr S$ to $\psi \colon V \arr S$, and $\alpha\colon \phi^* \cM \twoheadrightarrow \cL$ is an object of $Q_{\cM}(U)$, then
   \[
   f^*\alpha \colon \psi^* \cE \simeq f^*\phi^* \cM
   \xymatrix@C-5pt{{}\ar@{->>}[r]&{}} f^* \cL
   \]
is an object of $Q_{\cM}(V)$: this defines the pullback $Q_{\cM}(U) \arr Q_{\cM}(V)$. Then this functor is represented by $\PP(\cM)$.

When $\cM = \cO_{S}^{n+1}$, we recover the functor $Q_{n}$ of the previous example.
\end{example}

\begin{example}\label{ex:grassmanians}
With the same setup  as in the previous example, fix a positive integer $r$. We consider the functor $\catsch{S}\op \arr \catset$ that sends each $\phi \colon U \arr S$ to the set of quotients of $\phi^{*}\cM$ that are locally free of rank $r$. This is also representable by a scheme $\GG(r,\cM) \arr S$.
\end{example}

Finally, let us a give an example of a functor that is not representable.

\begin{example} This is very similar to
Example~(\ref{ex:notrep-opensubsets}) of \S\ref{subsec:examples}. Let
$\kappa$ be a field, $\catsch{\kappa}$ the category of schemes over
$\kappa$. Consider the functor $F \colon \catsch{\kappa}\op \arr \catset$
that associates with each scheme $U$ over $\kappa$ the set of all of its
open subsets; the action of $F$ on arrows is obtained by taking inverse
images. 

I claim that this functor is not representable. In fact, suppose that it is represented by a pair $(X, \xi)$, where $X$ is a scheme over $\kappa$ and $\xi$ is an open subset. We can consider $\xi$ as an open subscheme of $X$. If $U$ is any scheme over $\kappa$, a morphism of $\kappa$-schemes $U \arr \xi$ is a  morphism of $\kappa$-schemes $U \arr X$ whose image is contained in $\xi$; by definition of $X$ there is a unique such morphism, the one corresponding to the open subset $U$, considered as an element of $FU$. Hence the functor represented by the $\kappa$-scheme $\xi$ is the one point functor, sending any $\kappa$-scheme $U$ to a set with one element, and this is represented by $\spec \xi$. Hence $\xi$ is isomorphic to $\spec \kappa$ as a $\kappa$-scheme; this means that $\xi$, viewed as an open subscheme of $X$, consists of a unique $\kappa$-rational point of $X$. But a $\kappa$-rational point of a $\kappa$-scheme is necessarily a closed point (this is immediate for affine schemes, and follows in the general case, because being a closed subset of a topological space is a local property). So $\xi$ is also closed; but this would imply that every open subset of a $\kappa$-scheme is also closed, and this fails, for example, for $\AA^1_{\kappa}
\setminus \{0\} \subseteq \AA^1_\kappa$.
\end{example}

\begin{remark}\label{rmk:dual-Yoneda}
\index{Yoneda Lemma!dual form}
\index{Lemma! Yoneda!dual form}
There is a dual version of Yoneda's lemma, which will be used in  \S\ref{subsec:fibered-quasi-coherent}. Each object $X$ of $\cC$ defines a functor
   \[
   \hom_{\cC}(X,-) \colon \cC \arr \catset.
   \]
This can be viewed as the functor $\h_X \colon (\cC\op)\op \arr \catset$; hence, from the usual form of Yoneda's lemma applied to $\cC\op$ for any two objects $X$ and $Y$ we get a canonical bijective correspondence between $\hom_{\cC}(X,Y)$ and the set of natural transformations $\hom_{\cC}(Y,-) \arr \hom_{\cC}(X,-)$.
\end{remark}

\section{Group objects}\label{sec:group-objects}

In this section the category $\cC$ will have finite products; we will denote a terminal object by $\pt$.

\begin{definition}
A \emph{group object}\index{group object} of $\cC$ is an object $G$ of $\cC$, together with a functor $\cC\op \arr \catgrp$ into the category of groups, whose composite with the forgetful functor $\catgrp \arr \catset$ equals $\h_G$.

A group object in the category of topological spaces is called a \emph{topological group}\index{topological!group}. A group object in the category of schemes over a scheme $S$ is called a \emph{group scheme over $S$}\index{group scheme}.
\end{definition}

Equivalently: a group object is an object $G$, together with a group structure on $G(U)$ for each object $U$ of $\cC$, so that the function $f^* \colon G(V) \arr G(U)$ associated with an arrow $f \colon U \arr V$ in $\cC$ is always a homomorphism of groups.

This can be restated using Yoneda's lemma.

\begin{proposition}\index{group object!characterization via diagrams}
To give a group object structure on an object $G$ of $\cC$ is equivalent to assigning three arrows $\mul_G \colon G \times G \arr G$ (the multiplication), $\inv_G \colon G \arr G$ (the inverse), and $\iden_G \colon \pt \arr G$ (the identity), such that the following diagrams commute.

\begin{enumeratei}

\item The identity is a left and right identity:
   \[
   \xymatrix@C+10pt{
   {}\pt \times G \ar[r]^{\iden_G \times \id_G}\ar@{=}[rd]&
			G \times G\ar[d]^{\mul_G}\\
			&G
   }
\quad\text{and}\quad
   \xymatrix@C+10pt{
   G \times \pt \ar[r]^{\id_G \times\iden_G}\ar@{=}[rd]&
			G \times G\ar[d]^{\mul_G}\\
			&G
   }
   \]

\item Multiplication is associative:
   \[
   \xymatrix@C+20pt{
   G \times G \times G\ar[r]^-{\mul_G \times \id_G}
			\ar[d]^{\id_G \times \mul_G} &
			G \times G\ar[d]^{\mul_G}\\
   G \times G\ar[r]^-{\mul_G} &
			G
   }
   \]

\item The inverse is a left and right inverse:
   \[
   \xymatrix@C+10pt{
   G\ar[r]^-{\generate{\inv_G,\id_G}}\ar[d] &
			G \times G\ar[d]^{\mul_G}\\
			{}\pt \ar[r]^-{\iden_G}&
   G}
   \quad\text{and}\quad
   \xymatrix@C+10pt{
   G\ar[r]^-{\generate{\id_G, \inv_G}}\ar[d] &
			G \times G\ar[d]^{\mul_G}\\
			{}\pt \ar[r]^-{\iden_G}&
   G}
   \]
\end{enumeratei}
\end{proposition}

\begin{proof}
It is immediate to check that, if $\cC$ is the category of sets, the commutativity of the diagrams above gives the usual group axioms. Hence the result follows by evaluating the diagrams above (considered as diagrams of functors) at any object $U$ of $\cC$.
\end{proof}

Thus, for example, a topological group is simply a group, that has a structure of a topological space, such that the multiplication map and the inverse map are continuous (of course the map from a point giving the identity is automatically continuous).

Let us give examples of group schemes.

The first examples are the schemes $\AA^n_S \arr S$; these represent the functor $\cO^n $ sending a scheme $U \arr S$ to the set $\cO(U)^n$, which has an evident additive group structure. 

The group scheme $\AA^1_S$ is often denote by $\GG_{\rma, S}$.

Also, $\gm = \AA_S^1 \setminus 0_S$ represents the functor $\cO^*\colon \catsch{S}\op \arr \catset$, that sends each scheme $U \arr S$ to the group $\cO^*(U)$; this gives $\gm$ an obvious structure of group scheme.

Now consider the functor $\catsch{S}\op \arr \catset$ that sends each scheme $U \arr S$ to the set $\M_n\bigl(\cO(U)\bigr)$ of $n \times n$ matrices with coefficients in the ring $\cO(U)$. This is obviously represented by the scheme $\M_{n,S} \eqdef \AA_S^{n^2}$. Consider the determinant mapping as morphism of schemes $\det \colon \M_{n,S} \arr \AA_S^1$; denote by $\GL_{n,S}$ the inverse image of the open subscheme $\gm \subseteq \AA_{S}^1$. Then $\GL_{n,S}$ is an open subscheme of $\M_{n,S}$; the functor it represents is the functor sending each scheme $U \arr S$ to the set of matrices in $\M_n\bigl(\cO(U)\bigr)$ with invertible determinants. But these are the invertible matrices, and they form a group. This gives $\GL_{n,S}$ the structure of a group scheme on $S$.

There are various subschemes of $\GL_{n,S}$ that are group schemes. For example, $\SL_{n,S}$, the inverse image of the identity section $1_S \colon S \arr \gm$ via the morphism $\det \colon \GL_{n,S} \arr \gm$ represents the functor sending each scheme $U \arr S$ to the group $\SL \bigl(\cO(U)\bigr)$ of $n \times n$ matrices with determinant~$1$.

We leave it to the reader to define the orthogonal group scheme $\rO_{n,S}$ and the symplectic group scheme $\mathrm{Sp}_{n,S}$.

\begin{definition}
If $G$ and $H$ are group objects, we define a \emph{homomorphism}\index{group object!homomorphism of group objects}\index{homomorphism of group objects} of group objects as an arrow $G \arr H$ in $\cC$, such that for each object $U$ of $\cC$ the induced function $G(U) \arr H(U)$ is a group homomorphism.

Equivalently, a homomorphism is an arrow $f \colon G \arr H$ such that the diagram
   \[
   \xymatrix{
   G \times G \ar[r]^-{\mul_G} \ar[d]_{f \times f} & G\ar[d]^f\\
   H \times H \ar[r]^-{\mul_H}                     & H
   }
   \]
commutes.
\end{definition}

The identity is obviously a homomorphism from a group object to itself. Furthermore, the composite of homomorphisms of group objects is still a homomorphism; thus, group objects in a fixed category form a category, which we denote by $\catgrpover{\cC}$.

\begin{remark}\label{rmk:preserve-products->preserve-groups}

Suppose that $\cC$ and $\cD$ are categories with products and terminal objects $\pt_{\cC}$ and $\pt_{\cD}$. Suppose that $F \colon \cC \arr \cD$ is a functor that preserves finite products, and $G$ is a group object in $\cC$. The arrow $\iden_G \colon \pt_{\cC} \arr G$ yields an arrow $F\iden_G \colon F\pt_{\cC} \arr FG$; this can be composed with the inverse of the unique arrow $F\pt_{\cC} \arr \pt_{\cD}$, which is an isomorphism, because $F\pt_{\cC}$ is a terminal object, to get an arrow $\iden_{FG} \colon \pt_{\cD} \arr FG$. Analogously one uses $F\mul_{G} \colon F(G \times G) \arr FG$ and the inverse of the isomorphism $F(G \times G) \simeq FG \times FG$ to define an arrow $\mul_{FG} \colon FG \times FG \arr FG$. Finally we set $\inv_{FG} \eqdef F\inv_{G} \colon FG \arr FG$.

We leave it to the reader to check that this gives $FG$ the structure of a group object, and this induces a functor from the category of group objects on $\cC$ to the category of group objects on $\cD$.
\end{remark}

\subsection{Actions of group objects}\label{subsec:actions}

There is an obvious notion of left action of a functor into groups on a functor into sets.

\begin{definition}\index{action!of a functor into groups on a functor into sets}
A left action $\alpha$ of a functor $G \colon \cC\op \arr \catgrp$ on a functor $F \colon \cC\op \arr \catset$ is a natural transformation $G \times F \arr F$, such that for any object $U$ of $\cC$, the induced function $G(U) \times F(U) \arr F(U)$ is an action of the group $G(U)$ on the set $F(U)$.
\end{definition}

In the definition above, we denote by $G \times F$ the functor that sends an object $U$ of $\cC$ to the product of the set underlying the group $GU$ with the set $FU$. In other words, $G \times F$ is the product $\widetilde{G} \times F$, where $\widetilde{G}$ is the composite of $G$ with the forgetful functor $\catgrp \arr \catset$.

Equivalently, a left action of $G$ on $F$ consists of an action of $G(U)$ on $F(U)$ for all objects $U$ of $\cC$, such that for any arrow $f \colon U \arr V$ in $\cC$, any $g \in G(V)$ and any $x \in F(V)$ we have
   \[
   f^*g \cdot f^*x = f^*(g \cdot x) \in F(U).
   \]

Right actions are defined analogously.

We define an action of a group object $G$ on an object $X$ as an action of the functor $\h_G \colon \cC\op \arr \catgrp$ on $\h_X \colon \cC\op \arr \catset$\index{action!of a group object on an object}.

Again, we can reformulate this definition in terms of diagrams.

\begin{proposition}\label{prop:diagram-action}
\index{action!of a group object on an object!characterization via diagrams}
Giving a left action of a group object $G$ on an object $X$ is equivalent to assigning an arrow $\alpha \colon G \times X \arr X$, such that the following two diagrams commute.

\begin{enumeratei}

\item The identity of $G$ acts like the identity on $X$:
   \[
   \xymatrix@C+20pt{
   {}\pt \times X\ar@{=}[rd]\ar[r]^-{\iden_G \times \id_X}&
   G\times X\ar[d]^{\alpha}\\
   &X
   }
   \]

\item The action is associative with respect to the multiplication
on $G$:
   \[
   \xymatrix@C+20pt{
   G \times G \times X
			\ar[r]^-{\mul_G \times\id_X}\ar[d]^{\id_G \times\alpha}&
   G \times X\ar[d]^{\alpha}\\
   G \times X\ar[r]^{\alpha}&
   X
   }
   \]

\end{enumeratei}
\end{proposition}

\begin{proof}
It is immediate to check that, if $\cC$ is the category of sets, the commutativity of the diagram above gives the usual axioms for a left action. Hence the result follows from Yoneda's lemma by evaluating the diagrams above (considered as diagrams of functors) on any object $U$ of $\cC$.
\end{proof}

\begin{definition}
Let $X$ and $Y$ be objects of $\cC$ with an action of $G$, an arrow $f \colon X \arr Y$ is called \emph{$G$-equivariant}\index{equivariant!arrow}\index{arrow!equivariant} if for all objects $U$ of $\cC$ the  induced function $X(U) \arr Y(U)$ is $G(U)$-equivariant.
\end{definition}

Equivalently, $f$ is $G$-equivariant if the diagram
   \[
   \xymatrix{
   G \times X \ar[r] \ar[d]^{\id_G \times f}
   &
   X \ar[d]^f\\
   G \times Y \ar[r]&
   Y\\
   }
   \]
where the rows are given by the actions, commutes (the equivalence of these two definitions follows from Yoneda's lemma).

There is yet another way to define the action of a functor $G \colon \cC\op \arr \catgrp$ on an object $X$ of $\cC$. Given an object $U$ of $\cC$, we denote by $\End_U(U \times X)$ the set of arrows $U \times X \arr U \times X$ that commute with the projection $\pr_1 \colon U \times X \arr U$; this set has the structure of a monoid, the operation being the composition. In other words, $\End_U(U \times X)$ is the monoid of endomorphisms of  $\pr_1 \colon U \times X \arr U$ considered as an object of the comma category $(\cC/U)$. We denote the group of automorphisms in $\End_U(U \times X)$ by $\aut_U(U \times X)$.

Let us define a functor
   \[
   \underaut_{\cC}(X) \colon \cC\op \arr \catgrp
   \]
sending each object $U$ of $\cC$ to the group $\underaut_{\cC}(X)(U) \eqdef \aut_U(U \times X)$. The group $\underaut_{\cC}(X)(\pt)$ is canonically isomorphic to $\aut_{\cC}(X)$.

Consider an arrow $f \colon U \arr V$ in $\cC$; with this we need to associate a group homomorphism $f^* \colon \aut_V(V \times X) \arr \aut_U(U \times X)$. The diagram
   \[
   \xymatrix@C+10pt{
   U \times X \ar[r]^{f \times \id_X} \ar[d]^{\pr_1} &
   V \times X \ar[d]^{\pr_1}\\
   U \ar[r]^f & V
   }
   \]
is cartesian; hence, given an arrow $\beta \colon V \times X \arr V \times X$ over $V$, there is a unique arrow $\alpha \colon U \times X \arr U \times X$ making the diagram
   \[
   \xymatrix@C+10pt{
   U \times X
   \ar[rr]^{f \times \id_X} \ar@/_7pt/[rdd]_{\pr_1}
   \ar@{-->}[rd]^{\alpha}
   &&
   V \times X \ar[d]^{\beta}\\
   & U \times X \ar[r]^{f \times \id_X} \ar[d]^{\pr_1} &
   V \times X \ar[d]^{\pr_1}\\
   & U \ar[r]^f & V
   }
   \]
commute. This gives a function from the set $\End_{V}(V \times X)$ to $\End_{U}(U \times X)$, which is easily checked to be a homomorphism of monoids (that is, it sends the identity to the identity, and it preserves composition). It follows that it restricts to a homomorphism of groups $f^* \colon \aut_V(V \times X) \arr \aut_U(U \times X)$. This gives $ \underaut_{\cC}(X)$ the structure of a functor.

This construction is a very particular case of that of Section~\ref{sec:func-arrows}.

\begin{proposition}\label{prop:action-homo}\index{action!of a group object on an object!as a group homomorphism}
Let $G \colon \cC\op \arr \catgrp$ a functor, $X$ an object of $\cC$. To give an action of $G$ on $X$ is equivalent to giving a natural transformation $G \arr \underaut_{\cC}(X)$ of functors $\cC\op \arr \catgrp$.
\end{proposition}

\begin{proof}
Suppose that we are given a natural transformation $G \arr \underaut_{\cC}(X)$. Then for each object $U$ of $\cC$ we have a group homomorphism $G(U) \arr \aut_U(U \times X)$. The set $X(U)$ is in bijective correspondence with the set of sections $U \arr U \times X$ to the projection $\pr_1 \colon U \times X \arr U$, and if $s \colon U \arr U \times X$ is a section, $\alpha \in \aut_U(U \times X)$, then $\alpha \circ s \colon U \arr U \times X$ is still a section. This induces an action of $\aut_U(U \times X)$ on $X(U)$, and, via the given homomorphism $G(U) \arr \aut_U(U \times X)$, also an action of $G(U)$ on $X(U)$. It is easy to check that this defines an action of $G$ on $X$.

Conversely, suppose that $G$ acts on $X$, let $U$ be an object of $\cC$, and $g \in G(U)$. We need to associate with $g$ an object of 
   \[
   \underaut_{\cC}(X)(U) = \aut_{U}(U \times X).
   \]
We will use Yoneda's lemma once again, and consider $U \times X$ as a functor $U \times X \colon (\cC/U)\op \arr \catset$. For each arrow $V \arr U$ in $\cC$ there is a bijective correspondence between the set $X(V)$ and the set of arrows $V \arr U \times X$ in $\cC/U$, obtained by composing an arrow $V \arr U \times X$ with the projection $\pr_{2}\colon U \times X \arr X$. Now we are given an action of $G(V)$ on $X(V)$, and this induces an action of $G(V)$ on $\hom_{(\cC/U)}(V, U \times X) = (U \times X)(V)$. The arrow $V \arr U$ induces a group homomorphism $G(U) \arr G(V)$, so the element $g \in G(U)$ induces a permutation of $(U \times X)(V)$. There are several things to check: all of them are straightforward and left to the reader as an exercise.

\begin{enumeratei}

\item This construction associates with each $g$ an automorphism of the functor $U \times X$, hence an automorphism of $U \times X$ in $(\cC/U)$.

\item The resulting function $G(U) \arr \underaut_{\cC}(X)(U)$ is a group homomorphism.

\item This defines a natural transformation $G \arr \underaut_{\cC}(X)$.

\item The resulting functions from the set of actions of $G$ on $X$ and the set of natural transformations $G \arr \underaut_{\cC}(X)$ are inverse to each other.\qedhere

\end{enumeratei}
\end{proof}

\subsection{Discrete groups}\label{subsec:discrete-groups}

There is a standard notion of action of a group on an object of a category: a group $\Gamma$ acts on an object $X$ of $\cC$ when there is given a group homomorphism $\Gamma \arr \aut_{\cC}(X)$. With appropriate hypotheseas, this action can be interpreted as the action of a \emph{discrete group object} on $X$.

In many concrete cases, a category of geometric objects has objects that can be called \emph{discrete}. For example in the category of topological spaces we have discrete spaces: these are spaces with the discrete topology, or, in other words, disjoint unions of points. In the category $\catsch{S}$ of schemes over $S$ an object should be called discrete when it is the disjoint union of copies of $S$. In categorical terms, a disjoint union is a coproduct; thus a discrete object of $\catsch{S}$ is a scheme $U$ over $S$, with the property that the functor $\hom_S(U, -) \colon \catsch{S} \arr \catset$ is the product of copies of $\hom_S(S, -)$.

\begin{definition}\call{def:has-discrete-objects}
Let $\cC$ be a category. We say that $\cC$ \emph{has discrete objects}\index{category!with discrete objects} if it has a terminal object $\pt$, and for any set $I$ the coproduct $\coprod_{i \in I} \pt$ exists.

An object of $\cC$ that is isomorphic to one of the form $\coprod_{i \in I} \pt$ for some set $I$ is called a \emph{discrete object}\index{object!discrete}\index{discrete!object}.
\end{definition}

Suppose that $\cC$ has discrete objects. If $I$ and $J$ are two sets and $\phi \colon I \arr J$ is a function, we get a collection of arrows $\pt \arr \coprod_{j \in J} \pt$ parametrized by $I$: with each $i \in I$ we associate the tautological arrow $\pt \arr \coprod_{j \in J} \pt$ corresponding to the element $\phi(i) \in J$. In this way we have defined an arrow
   \[
   \phi_* \colon \coprod_{i \in I}\pt \arr
   \coprod_{j \in J}\pt.
   \]
It is immediate to check that if $\phi \colon I \arr J$ and $\psi \colon J \arr K$ are functions, we have
   \[
   (\psi \circ \phi)_* = \psi_* \circ \phi_* \colon 
   \coprod_{i \in I}\pt \arr \coprod_{k \in K}\pt.
   \]
In this way we have defined a functor $\Delta \colon \catset \arr
\cC$ that sends a set $I$ to $\coprod_{i \in I} \pt$. This is called \emph{the discrete object functor}\index{discrete object functor}\index{functor!discrete object}\index{discrete!object functor}. By construction, it is a left adjoint to the functor $\hom_{\cC}(\pt, -)$. Recall that this means that for every set $I$ and every object $U$ of $\cC$ one has a bijective correspondence between $\hom_{\cC}(\Delta I, U)$ and the set of functions $I \arr \hom_{\cC}(\pt, U)$; furthermore this bijective correspondence is functorial in $I$ and $U$.

Conversely, if we assume that $\cC$ has a terminal object $\pt$, and that $\Delta \colon \catset \arr \cC$ is a left adjoint to the functor $\hom_{\cC}(\pt, -)$, then it is easy to see for each set $I$ the object $\Delta I$ is a coproduct $\coprod_{i \in I} \pt$.

We are interested in constructing discrete group objects in a category $\cC$; for this, we need to have discrete objects, and, according to Remark~\ref{rmk:preserve-products->preserve-groups}, we need to have that the discrete object functor $\catset \arr \cC$ preserves finite products. Here is a condition to ensure that this happens.

Suppose that $\cC$ is a category with finite products. Assume furthermore that for any object $U$ in $\cC$ and any set $I$ the coproduct $\coprod_{i \in I} U$ exists in $\cC$; in particular, $\cC$ has discrete objects.  If $U$ is an object of $\cC$ and $I$ is a set, we will set $I \times U \eqdef \coprod_{i \in I}U$. By definition, an arrow $I \times U \arr V$ is defined by a collection of arrows $f_i \colon U \arr V$ parametrized by $I$. In particular, $\emptyset \times U$ is an initial object of $\cC$.

Notice the following fact. Let $I$ be a set, $U$ an object of 
$\cC$. If $i \in I$ then $\Delta\{i\}$ is a terminal object of $\cC$, hence there is a canonical isomorphism $U \simeq \Delta\{i\} \times U$ (the inverse of the projection $\Delta \{i\} \times U \arr U$). On the other hand the embedding $\iota_i \colon \{i\} \into I$ induces an arrow $\Delta\iota_i \times \id_U \colon \Delta\{i\} \times U \arr \Delta I \times U$. By composing these $\Delta\iota_i$ with the isomorphisms $U \simeq \Delta\{i\} \times U$ we obtain a set of arrows $U \arr \Delta I \times U$ parametrized by $I$, hence an arrow $I \times U \arr \Delta I \times U$.

\begin{definition}\call{def:has-discrete-groups}
A category $\cC$ \emph{has discrete group objects}\index{category!with discrete group object} when the following conditions are satisfied.

\begin{enumeratei}

\itemref{1} $\cC$ has finite products.

\itemref{2}  For any object $U$ in $\cC$ and any set $I$ the coproduct $I \times U \eqdef \coprod_{i \in I} U$ exists;

\itemref{3} For any object $U$ in $\cC$ and any set $I$, the canonical arrow $I \times U \arr \Delta I \times U$ is an isomorphism.
\end{enumeratei}

\end{definition}

For example, in the category $\cattop$ a terminal object is a point (in other words, a topological space with one element), while the coproducts are disjoint unions. The conditions of the definition are easily checked. This also applies to the category $\catsch S$ of schemes over a fixed base scheme $S$; in this case a terminal object is $S$ itself.

\begin{proposition}\label{prop:has-discrete-groups}
If $\cC$ has discrete group objects, then the discrete object functor $\Delta \colon \catset \arr \cC$ preserves finite products.\index{discrete object functor!preserves finite products}\index{functor!discrete object!preserves finite products}\index{discrete!object functor!preserves finite products}
\end{proposition}

So, by Remark~\ref{rmk:preserve-products->preserve-groups}, when the category $\cC$ has discrete group object the functor $\Delta \colon \catset \arr \cC$ gives a functor, also denoted by $\Delta \colon \catgrp \arr \catgrpover{\cC}$, from the category of groups to the category of group objects in $\cC$. A group object in $\cC$ is called \emph{discrete}\index{group object!discrete}\index{discrete!group object} when it isomorphic to one of the form $\Delta\Gamma$, where $\Gamma$ is a group.

\begin{proof}

Let $\cC$ be a category with discrete group objects. To prove that $\Delta$ preserves finite products, it is enough to check that $\Delta$ sends a terminal object to a terminal object, and that it preserves products of two objects. The  first fact follows  immediately from the definition of $\Delta$.

Let us show that, given two sets $I$ and $J$, the natural arrow $\Delta(I \times J) \arr \Delta I \times \Delta J$ is an isomorphism. By definition, $\Delta(I \times J) = (I \times J) \times \pt$. On the other hand there is a canonical well-known isomorphism of
   \[
   (I \times J) \times \pt = \coprod_{(i,j) \in I \times J} \pt
   \]
with
   \[
   \coprod_{i \in I} \biggl(\,\coprod_{j \in J}
   \pt\biggr) = I \times (J
   \times \pt) = I \times \Delta J.
   \]
If we compose this isomorphism $\Delta(I \times J) \simeq I \times \Delta J$ with the isomorphism $I \times \Delta J \simeq \Delta I \times \Delta J$ discussed above we obtain an isomorphism $\Delta(I \times J) \simeq \Delta I \times \Delta J$. It is easy to check that the projections $\Delta(I \times J) \arr \Delta I$ and $\Delta(I \times J) \arr \Delta J$ are induced by the projections $I \times J \arr I$ and $I \times J \arr J$; this finishes the proof.
\end{proof}

An action of a group is the same as an action of the associated discrete group object.

\begin{proposition}\label{prop:action-discrete-group}\index{action!of a discrete group object}
Suppose that $\cC$ has finite group objects. Let $X$ be an object of $\cC$, $\Gamma$ a group, $\Delta\Gamma$ the associated discrete group object of $\cC$. Then giving an action of $\Gamma$ on $X$, that is, giving a group homomorphism $\Gamma \arr \aut_{\cC}(X)$, is equivalent to giving an action of the group object $\Delta\Gamma$ on $X$.
\end{proposition}

\begin{proof}
A function from $\Gamma$ to the set $\hom_{\cC}(X, X)$ of arrows from $X$ to itself corresponds, by definition, to an arrow $\Gamma \times X \arr X$; the isomorphism $\Gamma \times X \simeq \Delta\Gamma \times X$ above gives a bijective correspondence between functions $\Gamma \arr \hom_{\cC}(X, X)$ and arrows $\Delta\Gamma \times X \arr X$. We have to check that a function $\Gamma \arr \hom_{\cC}(X, X)$ gives an action of $\Gamma$ on $X$ if and only if the corresponding arrow $\Delta\Gamma \times X \arr X$ gives an action of $\Delta\Gamma$ on $X$. This is straightforward and left to the reader. 
\end{proof}

\begin{remark}
The terminology ``$\cC$ has discrete group objects'' is perhaps misleading; for $\cC$ to have discrete group objects would be sufficient to have discrete objects, and that the functor $\Delta$ preserves finite products.

However, for discrete group objects to be well behaved we need more than their existence, we want Proposition~\ref{prop:action-discrete-group} to hold: and for this purpose the conditions of Definition~\ref{def:has-discrete-groups} seem to be optimal (except that one does not need to assume that $\cC$ has all products; but this hypothesis is satisfied in all the examples I have in mind).
\end{remark}

\section{Sheaves in Grothendieck topologies}\label{sec:sheaves}

\subsection{Grothendieck topologies}

The reader is familiar with the notion of sheaf on a topological space. A presheaf on a topological space $X$ can be considered as a functor. Denote by $X\cl$ the category in which the objects are the open subsets of $X$, and the arrows are given by inclusions. Then a presheaf of sets on $X$ is a functor ${X\cl}\op \arr \catset$; and this is a sheaf when it satisfies appropriate gluing conditions. 

There are more general circumstances under which we can ask whether a functor is a sheaf. For example, consider a functor $F \colon \cattop\op \arr \catset$; for each topological space $X$ we can consider the restriction $F_X$ to the subcategory $X\cl$ of $\cattop$. We say that $F$ is a $\emph{sheaf}$ on $\cattop$ if $F_X$ is a sheaf on $X$ for all $X$.

There is a very general notion of  sheaf in a Grothendieck topology; in this Section we review this theory.

In a Grothendieck topology the ``open sets'' of a space are \emph{maps} into this space; instead of intersections we have to look at fibered products, while unions play no role. The axioms do not describe the ``open sets'', but the coverings of a space.

\begin{definition}\call{def:grothtop}
Let $\cC$ be a category. A \emph{Grothendieck topology}\index{topology!Grothendieck}\index{Grothendieck!topology} on $\cC$ is the assignment to each object $U$ of $\cC$ of a collection of sets of arrows $\{U_i \to U\}$, called \emph{coverings of $U$}\index{covering}, so that the following conditions are satisfied.

\begin{enumeratei}

\itemref{1} If $V \arr U$ is an isomorphism, then the set $\{V \arr U\}$ is a covering.

\itemref{2} If $\{U_i \to U\}$ is a covering and $V \arr U$ is any arrow, then the fibered products $\{U_i \times_U V\}$ exist,  and the collection of projections $\{U_i \times_U V \arr V \}$ is a covering.

\itemref{3} If $\{U_i \to U\}$ is a covering, and for each index $i$ we have a covering $\{V_{ij} \arr U_i\}$ (here $j$ varies on a set depending on $i$), the collection of composites $\{V_{ij} \arr U_i \arr U\}$ is a covering of $U$.

\end{enumeratei}

A category with a Grothendieck topology is called a \emph{site}\index{site}.
\end{definition}

Notice that from \refpart{def:grothtop}{2} and \refpart{def:grothtop}{3} it follows that if $\{U_i \to U\}$ and $\{V_j \to U\}$ are two coverings of the same object, then $\{U_i \times_U V_j \arr U\}$ is also a covering.

\begin{remark}
In fact what we have defined here is what is called a \emph{pretopology}\index{pretopology} in \cite{sga4}; a pretopology defines a topology, and very different pretopologies can define the same topology. The point is that the sheaf theory only depends on the topology, and not on the pretopology. Two pretopologies induce the same topology if and only if they are equivalent, in the sense of Definition~\ref{def:same-sheaves}.

Despite its unquestionable technical advantages, I do not find the notion of topology, as defined in \cite{sga4}, very intuitive, so I prefer to avoid its use (just a question of habit, undoubtedly). 

However, \emph{sieves}, the objects that intervene in the definition of a topology, are quite useful, and will be used extensively.
\end{remark}

Here are some examples of Grothendieck topologies. In what follows, a set $\{U_{i} \arr U\}$ of functions, or morphisms of schemes, is called \emph{jointly surjective}\index{jointy surjective} when the set-theoretic union of their images equals $U$.

\begin{example}[The site of a topological space]
\label{ex:classical-topology}\index{classical topology}\index{topology!classical}
Let $X$ be a topological space; denote by $X\cl$ the category in which the objects are the open subsets of $X$, and the arrows are given by inclusions. Then we get a Grothendieck topology on $X\cl$ by associating with each open subset $U \subseteq X$ the set of open coverings of $U$.

In this case, if $U_1 \arr U$ and $U_2 \arr U$ are arrows, the fibered product $U_1 \times_U U_2$ is the intersection $U_1 \cap U_2$.
\end{example}

\begin{example}[The global classical topology]\label{ex:global-classical}\index{global!classical topology}\index{topology!classical!global}
Here $\cC$ is the category $\cattop$ of topological spaces. If $U$ is a topological space, then a covering of $U$ will be a jointly surjective collection of open embeddings $U_i \arr U$.

Notice here we must interpret ``open embedding'' as meaning an open continuous injective map $V \arr U$; if by an open embedding we mean the inclusion of an open subspace, then condition \refpart{def:grothtop}{1} of Definition~\ref{def:grothtop} is not satisfied.
\end{example}

\begin{example}[The global \'etale topology for topological spaces]\label{ex:top-global-etale}
Here $\cC$ is the category $\cattop$ of topological spaces. If $U$ is a topological space, then a covering of $U$ will be a jointly surjective collection of local homeomorphisms $U_i \arr U$.
\end{example}

Here is an extremely important example from algebraic geometry.

\begin{example}[The small \'etale site of a scheme]
\index{topology!\'etale}
\index{\'etale topology}
\index{site!small \'etale}
\index{small \'etale site}
Let $X$ be a scheme. Consider the full subcategory $X_{\textrm{\'et}}$ of $\catsch{X}$, consisting of morphisms $U \arr X$ locally of finite presentation, that are \'etale. If $U \arr X$ and $V \arr X$ are objects of $X_{\textrm{\'et}}$, then an arrow $U \arr V$ over $X$ is necessarily \'etale.

A covering of $U \arr X$ in the small \'etale topology is a jointly surjective collection of morphisms $U_i \arr U$.
\end{example}

Here are topologies that one can put on the category $\catsch{S}$ of schemes over a fixed scheme $S$. Several more have been used in different contexts.

\begin{example}[The global Zariski topology]
\index{topology!global Zariski}
\index{global!Zariski topology}
\index{Zariski topology!global}
Here a covering $\{U_i \to U\}$ is a collection of open embeddings covering $U$. As in the example of the global classical topology, an open embedding must be defined as a morphism $V \arr U$ that gives an isomorphism of $V$ with an open subscheme of $U$, and not simply as the embedding of an open subscheme.
\end{example}

\begin{example}[The global \'etale topology]
\index{topology!global \'etale}
\index{global!\'etale topology}
\index{\'etale topology!global}
A covering $\{U_i \to U\}$ is a jointly surjective collection  of \'etale maps locally of finite presentation.
\end{example}

\begin{example}[The fppf topology]
\index{topology!fppf}
\index{fppf topology}
A covering $\{U_i \to U\}$ is a jointly surjective collection of flat maps locally of finite presentation.

The  abbreviation fppf stands for ``fid\`element plat et de pr\'esentation finie''.
\end{example}

\subsection{The fpqc topology}\label{subsec:fpqc}

It is sometimes useful to consider coverings that are not locally finitely presented. One can define a topology on $\catsch{S}$ simply by taking all collections of morphisms $\{U_i \to U\}$ such that the resulting morphism $\coprod_i U_i \arr U$ is faithfully flat\index{wild flat topology}\index{topology!wild flat}. Unfortunately, this topology is not well behaved (see Remarks \ref{rmk:need-finiteness} and \ref{rmk:really-need-finiteness}). One needs some finiteness condition in order to get a reasonable topology. 

For example, one could define a covering as a collection of morphisms $\{U_i \to U\}$ such that the resulting morphism $\{\coprod_i U_i \arr U\}$ is faithfully flat and quasi-compact, as I did in the first version of these notes; but then Zariski covers would not be included, and the resulting topology would not be comparable with the Zariski topology. The definition of the fpqc topology that follows, suggested by Steve Kleiman, gives the correct sheaf theory.

\begin{proposition}\call{prop:eq-quasicompact}
Let $f \colon X \arr Y$ be a surjective morphism of schemes. Then the following properties are equivalent.

\begin{enumeratei}

\itemref{1} Every quasi-compact open subset of $Y$ is the image of a quasi-compact open subset of $X$.

\itemref{2} There exists a covering $\{V_i\}$ of $Y$ by open affine subschemes, such that each $V_i$ is the image of a quasi-compact open subset of $X$.

\itemref{3} Given a point $x \in X$, there exists an open neighborhood $U$ of $x$ in $X$, such that the image $fU$ is open in $Y$, and the restriction $U \arr fU$ of $f$ is quasi-compact.

\itemref{4} Given a point $x \in X$, there exists a quasi-compact open neighborhood $U$ of $x$ in $X$, such that the image $fU$ is open and affine in $Y$.

\end{enumeratei}
\end{proposition}

\begin{proof}

It is obvious that \refpart{prop:eq-quasicompact}{1} implies \refpart{prop:eq-quasicompact}{2}. The fact that \refpart{prop:eq-quasicompact}{4} implies \refpart{prop:eq-quasicompact}{3} follows from the fact that a morphism from a quasi-compact scheme to an affine scheme is quasi-compact.

It is also easy to show that \refpart{prop:eq-quasicompact}{3} implies \refpart{prop:eq-quasicompact}{4}: if $U'$ is an open subset of $X$ containing $x$, whose image $fU'$ in $Y$ is open, take an affine neighborhood $V$ of $f(x)$ in $fU'$, and set $U = f^{-1}V$.

Since $f$ is surjective, we see that \refpart{prop:eq-quasicompact}{4} implies \refpart{prop:eq-quasicompact}{2}.

Conversely, assuming \refpart{prop:eq-quasicompact}{2}, take a point $x \in X$. Then $f(x)$ will be contained in some $V_i$. Let $U'$ be a quasi-compact open subset of $X$ with image $V_i$, and $U''$ an open neighborhood of $x$ in $f ^{-1} V_i$. Then $U = U' \cup U''$ is quasi-compact, contains $x$ and has image $V_i$.

We only have left to prove that \refpart{prop:eq-quasicompact}{2} implies \refpart{prop:eq-quasicompact}{1}. Let $V$ be a quasi-compact open subset of $Y$. The open affine subsets of $Y$ that are contained in some $V \cap V_i$ form a covering of $V$, so we can choose finitely many of them, call them $W_1$, \dots,~$W_r$. Given one of the $W_j$, choose an index $i$ such that $W_j \subseteq V_i$ and a quasi-compact open subset $U_i$ of $X$ with image $V_i$; the restriction $U_i \arr V_i$ is quasi-compact, so the inverse image $W'_j$ of $W_j$ in $U_i$ is quasi-compact. Then $\bigcup_{j=1}^r W'_j$ is an open quasi-compact subscheme of $X$ with image $\bigcup_{j=1}^r W_j = V$.
\end{proof}

\begin{definition}
An \emph{fpqc morphism of schemes}\index{morphism of schemes!fpqc}\index{fpqc!morphism of schemes} is a faithfully flat morphism that satisfies the equivalent conditions of Proposition~\ref{prop:eq-quasicompact}.
\end{definition}

The  abbreviation fpqc stands for ``fid\`element plat et quasi-compact''.

Here are some properties of fpqc morphisms.

\begin{proposition}\call{prop:properties-fpqc}\hfil

\begin{enumeratei} \index{morphism of schemes!fpqc!properties of}\index{fpqc!morphism of schemes!properties of}

\itemref{1} The composite of fpqc morphisms is fpqc.

\itemref{5} If $f \colon X \arr Y$ is a morphism of schemes, and there is an open covering $V_i$ of $Y$, such that the restriction $f^{-1}V_i \arr V_i$ is fpqc, then $f$ is fpqc.

\itemref{3} An open faithfully flat morphism is fpqc.

\itemref{4} A faithfully flat morphism that is locally of finite
presentation is fpqc.

\itemref{2} A morphism obtained by base change from an fpqc morphism is fpqc.

\itemref{6} If $f \colon X \arr Y$ is an fpqc morphism, a subset of $Y$ is open in $Y$ if and only if its inverse image is open in $X$.
\end{enumeratei}

\end{proposition}

\begin{proof} \refpart{prop:properties-fpqc}{1} follows from the definition, using the characterization \refpart{prop:eq-quasicompact}{1} in Proposition~\ref{prop:eq-quasicompact}. Also \refpart{prop:properties-fpqc}{5} follows easily, using the characterization \refpart{prop:eq-quasicompact}{2}, and \refpart{prop:properties-fpqc}{3} follows from condition~\refpart{prop:eq-quasicompact}{3}. \refpart{prop:properties-fpqc}{4} follows from \refpart{prop:properties-fpqc}{3} and the fact that a faithfully flat morphism that is locally of finite presentation is open (Proposition~\ref{prop:flat->open}).

For \refpart{prop:properties-fpqc}{3}, suppose that we are given a cartesian diagram of schemes
   \[
   \xymatrix{
   X' \ar[r]\ar[d] & X \ar[d] \\
   Y' \ar[r] & Y
   }
   \]
such that $X \arr Y$ is fpqc. Take a covering $V_i$ of $Y$ by open affine subschemes, and for each of them choose an open quasi-compact open subset $U_i$ of $X$ mapping onto $V_i$. If we denote by $V_i'$ its inverse image of $V_i$ in $Y'$ and $U'_i$ the inverse image of $U_i$ in $X$, it is easy to check that $U'_i = V'_i \times_{V_i} U_i$. Since the morphism $U_i \arr V_i$ is quasi-compact, it follows that $U'_i \arr V'_i$ is also quasi-compact. Now take a covering $\{V''_j\}$ by open affine subschemes, such that each $V''_j$ is contained in some $V'_i$; then each $V''_j$ is the image of a quasi-compact open subset of $X'$, its inverse image in some $U_i$.

Let us prove \refpart{prop:properties-fpqc}{6}. Let $A$ be a subset of $Y$ whose inverse image in $X$ is open. Pick a covering $\{V_i\}$ of $Y$ by open affine subsets, each of which is the image of a quasi-compact open subset $U_i$ of $X$. Then the inverse image of $A$ in each $U_i$ will be open, and according to Proposition~\ref{prop:flat->quotient-topology} this implies that each $A \cap V_i$ is open in $V_i$, so $A$ is open in $Y$.
\end{proof}

The \emph{fpqc topology}\index{topology!fpqc}\index{fpqc!topology} on the category $\catsch{S}$ is the topology in which the coverings $\{U_i \to U\}$ are collections of morphisms, such that the induced morphism $\coprod U_i \arr U$ is fpqc.

Let us verify that this is indeed a topology, by checking the three conditions of Definition~\ref{def:grothtop}. Condition~\refpart{def:grothtop}{1} is obvious, because an isomorphism is fpqc.

Condition~\refpart{def:grothtop}{2} follows from Proposition~\refall{prop:properties-fpqc}{2}.

Condition~\refpart{def:grothtop}{3} is easy to prove, from parts
\refpart{prop:properties-fpqc}{1} and \refpart{prop:properties-fpqc}{5} of Proposition~\ref{prop:properties-fpqc}.

The fpqc topology is finer than the fppf topology, which is finer than the \'etale topology, which is in turn finer than the Zariski topology.

Many properties of morphisms are local on the codomain in the fpqc topology.

\begin{proposition}\label{prop:local-fpqc}\index{topology!fpqc!local properties in the}\index{fpqc!topology!local properties in the}
Let $X \arr Y$ be a morphism of schemes, $\{Y_i \arr Y\}$ an fpqc covering. Suppose that for each $i$ the projection $Y_i \times_{Y} X \arr Y_i$ has one of the following properties:

\begin{enumeratei}

\item is separated,

\item is quasi-compact,

\item is locally of finite presentation,

\item is proper,

\item is affine,

\item is finite,

\item is flat,

\item is smooth, 

\item is unramified,

\item is \'etale,

\item is an embedding,

\item is a closed embedding.

\end{enumeratei}

Then $X \arr Y$ has the same property.
\end{proposition}

\begin{proof}
This follows easily from the fact that each of the properties above is local in the Zariski topology in the codomain, from the characterization \refpart{prop:eq-quasicompact}{2} in Proposition~\ref{prop:eq-quasicompact}, and from Proposition~\ref{prop:local-qcflat}.
\end{proof}

\subsection{Sheaves}

If $X$ is a topological space, a presheaf of sets on $X$ is a functor ${X\cl}\op \arr \catset$, where $X\cl$ is the category of open subsets of $X$, as in Example~\ref{ex:classical-topology}. The condition that $F$ be a sheaf can easily be generalized to any site, provided that we substitute intersections, which do not make sense, with fibered products. (Of course, fibered products in $X\cl$ are just intersections.)

\begin{definition}
Let $\cC$ be a site, $F \colon \cC\op \arr
\catset$ a functor.

\begin{enumeratei}

\item $F$ is \emph{separated}\index{functor!separated}\index{separated functor} if, given a covering $\{U_i \to U\}$ and two sections $a$ and $b$ in $F U$ whose pullbacks to each $F U_i$ coincide, it follows that $a = b$.

\item $F$ is a \emph{sheaf}\index{sheaf} if the following condition is satisfied. Suppose that we are given a covering $\{U_i \to U\}$ in $\cC$, and a set of elements $a_i \in F U_i$. Denote by $\pr_1 \colon U_i \times_U U_j \arr U_i$ and $\pr_2 \colon U_i \times_U U_j \arr U_j$ the first and second projection respectively, and assume that $\pr_1^* a_i = \pr_2^* a_j \in F(U_i \times_U U_j)$ for all $i$ and $j$. Then there is a unique section $a \in F U$ whose pullback to $F U_i$ is $a_i$ for all $i$.

If $F$ and $G$ are sheaves on a site $\cC$, a \emph{morphism of sheaves}\index{morphism of sheaves} $F \arr G$ is simply a natural transformation of functors.

\end{enumeratei}
\end{definition}

A sheaf on a site is clearly separated.

Of course one can also define sheaves of groups, rings, and so on, as usual: a functor from $\cC\op$ to the category of groups, or rings, is a sheaf if its composite with the forgetful functor to the category of sets is a sheaf.

The reader might find our definition of sheaf pedantic, and wonder why we did not simply say ``assume that the pullbacks of $a_i$ and $a_j$ to $F(U_i \times_U U_j)$ coincide''. The reason is the following: when $i = j$, in the classical case of a topological space we have $U_i \times_U U_i = U_i \cap U_i = U_i$, so the two possible pullbacks from $U_i \times_U U_i \arr U_i$ coincide; but if the map $U_i \arr U$ is not injective, then the two projections $U_i \times_U U_i \arr U_i$ will be different. So, for example, in the classical case coverings with one subset are not interesting, and the sheaf condition is automatically verified for them, while in the general case this is very far from being true.

An alternative way to state the condition that $F$ is a sheaf is the following.

Let $A$, $B$ and $C$ be sets, and suppose that we are given a diagram
   \[
   \xymatrix{
   A\ar[r]^f & B\ar@<3pt>[r]^g \ar@<-3pt>[r]_h &C.
   }
   \]
(that is, we are given a function $f \colon A \arr B$ and two functions $g, h \colon B \arr C$). We say that the diagram is an \emph{equalizer}\index{equalizer} if $f$ is injective, and maps $A$ surjectively onto the subset $\{b \in B \mid g(b) = h(b)\} \subseteq B$.

Equivalently, the diagram is an equalizer if $g \circ f = h \circ f$, and every function $p \colon D \arr B$ such that $g \circ p = h \circ p$ factors uniquely through $A$.

Now, take a functor $F \colon \cC\op \arr \catset$ and a covering $\{U_i \to U\}$ in $\cC$. There is a diagram
   \begin{equation}\label{eq:equalizer-sheaves}
   FU \arr \prod_i
   FU_i\doublelong{\pr_1^*}{\pr_2^*}
   \prod_{i,j} F(U_i \times_U U_j)
   \end{equation}
where the function $FU \arr \prod_i FU_i$ is induced by the restrictions $FU \arr FU_i$, while
   \[
   \pr_1^* \colon \prod_i FU_i \arr
   \prod_{i,j} F(U_i \times_U U_j)
   \]
sends an element $(a_i) \in \prod_i FU_i$ to the element $\pr_1^*(a_i) \in \prod_{i,j} F(U_i \times_U U_j)$ whose component in $ F(U_i\times_U U_j)$ is the pullback $\pr_1^* a_i$ of $a_i$ along the first projection $U_i \times_U U_j \arr U_i$. The function
   \[
   \pr_2^* \colon \prod_i FU_i \arr
   \prod_{i,j} F(U_i \times_U U_j)
   \]
is defined similarly.

One immediately sees that $F$ is a sheaf if and only if the diagram (\ref{eq:equalizer-sheaves}) is an equalizer for all coverings $\{U_i \to U\}$ in $\cC$.

\subsection{Sieves}

Given an object $U$ in a category $\cC$ and a set of arrows $\cU = \{U_i \to U\}$ in $\cC$, we define a subfunctor $\h_{\cU}  \subseteq \h_U$\index{$\h_{\cU}$}, by taking $\h_{\cU}(T)$ to be the set of arrows $T \arr U$ with the property that for some $i$ there is a factorization $T \arr U_i \arr U$. In technical terms, $\h_{\cU}$ is the \emph{sieve}\index{sieve!associated with a set of arrows} associated with the covering $\cU$. The term is suggestive: think of the $U_{i}$ as holes on $U$. Then an arrow $T \arr U$ is in $\h_{\cU}T$ when it goes through one of the holes. So a sieve is determined by what goes through it.

\begin{definition}
Let $U$ be an object of a category $\cC$. A \emph{sieve}\index{sieve} on $U$ is a subfunctor of $\h_{U} \colon \cC\op \arr \catset$.
\end{definition}

Given a subfunctor $S \subseteq \h_{U}$, we get a collection $\cS$ of arrows $T \arr U$ (consisting of union of the $ST$  with $T$ running through all objects of $\cC$), with the property that every time an arrow $T \arr U$ is in $\cS$, every composite $T' \arr T \arr U$ is in $\cS$. Conversely, from such a collection we get a subfunctor $S \subseteq \h_{U}$, in which $ST$ is the set of all arrows $T \arr U$ that are in $\cS$.

Now, let $\cU = \{U_{i} \to U\}$ be a set of arrows, $F \colon \cC\op \arr \catset$ a functor. We define $F\cU$ to be the set of elements of $\prod_i FU_i$ whose images in $\prod_{i,j} F(U_i \times_U U_j)$ are equal. Then the restrictions $FU \arr FU_i$ induce a function $FU \arr F\cU$; by definition, a sheaf is a functor $F$ such that $FU \arr F\cU$ is a bijection for all coverings $\cU = \{U_{i} \to U\}$.

The set $F\cU$ can be defined in terms of sieves.

\begin{proposition}\label{prop:funnychar-gluing}
There is a canonical bijection 
   \[
   R \colon \hom(\h_{\cU}, F) \simeq F\cU
   \]
such that the diagram
   \[
   \xymatrix{
   {}\hom(\h_{U}, F)   \ar[r]\ar[d]  & FU \ar[d] \\
   {}\hom(\h_{\cU}, F) \ar[r]^-{R} & F\cU
   }
   \]
in which the top row is the Yoneda isomorphism, the left hand column is the restriction function induced by the embedding of $\h_{\cU}$ in $\h_{U}$ and the right hand column is induced by the restriction functions $FU \arr FU_{i}$, commutes.
\end{proposition}

\begin{proof}
Take a natural transformation $\phi\colon \h_{\cU} \arr F$. For each $i$, the arrow $U_{i} \to U$ is an object of $\h_{\cU}U_{i}$; from this we get an element $R\phi \eqdef \bigl(\phi(U_{i} \to U)\bigr) \in \prod_i FU_i$. The pullbacks $\pr_{1}^{*} \phi(U_{i} \to U)$ and $\pr_{2}^{*} \phi(U_{j} \to U)$ to $FU_{ij}$ both coincide with $\phi(U_{ij} \to U)$, hence $R\phi$ is an element of $F\cU$. This defines a function $R \colon \hom(\h_{\cU}, F) \arr F\cU$; the commutativity of the diagram is immediately checked.

We need to show that $R$ is a bijection. For this purpose take two natural transformations $\phi, \psi \colon \h_{\cU} \arr F$ such that $R\phi = R\psi$. Consider an element $T \arr U$ of some $\h_{\cU}T$; by definition, this factors as $T \xarr{f} U_{i} \arr U$ for some arrow $f \colon T \arr U_{i}$. Then by definition of a natural transformation we have
   \[
   \phi(T \to U) = f^{*}\phi(U_{i}\to U) =
   f^{*}\psi(U_{i}\to U) = \psi(T \to U),
   \]
hence $\phi = \psi$. This proves the injectivity of R.

For surjectivity, take an element $(\xi_{i}) \in F\cU$; we need to define a natural transformation $\h_{\cU} \arr F$. If $T \arr U$ is an element of $\h_{\cU}T$, choose a factorization $T \xarr{f} U_{i} \arr U$; this defines an element $f^{*}\xi_{i}$ of $FU$. This element is independent of the factorization: two factorizations $T \xrightarrow{f} U_i \arr U$ and $T \xrightarrow{g} U_j \arr U$ give an arrow $T \arr U_{ij}$, whose composites with $\pr_1 \colon U_{ij} \arr U_i$ and $\pr_2 \colon U_{ij} \arr U_j$ are equal to $f$ and $g$. Since $\pr^*_1 \xi_i = \pr^*_2 \xi_j$, we see that $f^{*}\xi_{i} = g^{*}\xi_{j}$.

This defines a function $\h_{\cU}T \arr FT$ for each $T$. We leave it to the reader to show that this defines a natural transformation $\phi \colon \h_{\cU} \arr F$, and that $R\phi = (\xi_{i})$.
\end{proof}

As an immediate corollary, we get the following characterization of sheaves.

\begin{corollary}\label{cor:funnychar-sheaves}
A functor $F \colon \cC\op \arr \catset$ is a sheaf if and only if for any covering $\cU = \{U_i \to U\}$ in $\cC$, the induced function
   \[
   FU \simeq \hom(\h_U, F) \arr  \hom(\h_\cU, F)
   \]
is bijective. Furthermore, $F$ is separated if and only if this function is always injective.
\end{corollary}

This characterization can be sharpened.

\begin{definition}
Let $\cT$ be a Grothendieck topology on a category $\cC$. A sieve $S \subseteq \h_{U}$ on an object $U$ of $\cC$ is said to \emph{belong to $\cT$}\index{sieve!belonging to a topology} if there exists a covering $\cU$ of $U$ such that $\h_{\cU} \subseteq S$.
\end{definition}

If $\cC$ is a site, we will talk about the \emph{sieves of $\cC$} to mean the sieves belonging to the topology of $\cC$.

The importance of the following characterization will be apparent after the proof of Proposition~\ref{prop:same-sheaves}.

\begin{proposition}\label{prop:funnychar-sheaves}\index{sheaf!characterization via sieves}\index{separated functor!characterization via sieves}\index{functor!separated!characterization via sieves}
A functor $F \colon \cC\op \arr \catset$ is a sheaf in a topology $\cT$ if and only if for any sieve $S$ belonging to $\cT$ the induced function
   \[
   FU \simeq \hom(\h_U, F) \arr  \hom(S, F)
   \]
is bijective. Furthermore, $F$ is separated if and only if this function is always injective.
\end{proposition}

\begin{proof}
The fact that this condition implies that $F$ is a sheaf is an immediate consequence of Corollary~\ref{cor:funnychar-sheaves}.

To show the converse, let $F$ be a sheaf, take a sieve $S \subseteq \h_{U}$ belonging to $\cC$, and choose a covering $\cU$ of $U$ with $\h_{\cU} \subseteq S$. The composite
   \[
   \hom(\h_{U}, F) \arr \hom(S, F) \arr \hom(\h_{\cU}, F)
   \]
is a bijection, again because of Corollary~\ref{cor:funnychar-sheaves}, so the thesis follows from the next Lemma.

\begin{lemma}\label{lem:separated->injective}
If $F$ is separated, the restriction function
   \[
   \hom(S, F) \arr \hom(\h_{\cU}, F)
   \]
is injective.
\end{lemma}

\begin{proof}
Let us take two natural transformations $\phi, \psi \colon S \arr F$ with the same image in $\hom(\h_{\cU}, F)$, an element $T \to U$ of $ST$, and let us show that $\phi(T \to U) = \psi(T \to U) \in FT$.

Set $\cU = \{U_{i} \to U\}$, and consider the fibered products $T \times_{U} U_{i}$ with their projections $p_{i}\colon T \times_{U} U_{i} \arr T$. Since $T \times_{U} U_{i} \arr U$ is in $\h_{\cU}(T \times_{U} U_{i})$ we have
   \[
   p_{i}^{*}\phi(T \to U) = \phi(T \times_{U} U_{i}) = 
   \psi(T \times_{U} U_{i}) = p_{i}^{*}\psi(T \to U).
   \]
Since $\{p_{i}\colon T \times_{U} U_{i} \arr T\}$ is a covering and $F$ is a separated presheaf, we conclude that $\phi(T \to U) = \psi(T \to U) \in FT$, as desired.
\end{proof}

This concludes the proof of Proposition~\ref{prop:funnychar-sheaves}.
\end{proof}

We conclude with a remark. Suppose that $\cU = \{U_{i}\to U\}$ and $\cV = \{V_{j}\to U\}$ are coverings. Then $\cU \times_{U} \cV \eqdef \{U_{i}\times_{U} V_{j}\to U\}$ is a covering. An arrow $T \arr U$ factors through $U_{i}\times_{U} V_{j}$ if and only if it factors through $U_{i}$ and through $V_{j}$. This simple observation is easily seen to imply the following fact.

\begin{proposition}\label{prop:directed-sieves}\hfil
\begin{enumerate}

\item If $\cU = \{U_{i}\to U\}$ and $\cV = \{V_{j}\to U\}$ are coverings, then
   \[
   \h_{\cU \times_{U} \cV} = \h_{\cU} \cap \h_{\cV} \subseteq \h_{U}.
   \]

\item If $S_{1}$ and $S_{2}$ are sieves on $U$ belonging to $\cT$, the intersection $S_{1} \cap S_{2} \subseteq \h_{U}$ also belongs to $\cT$.
\end{enumerate}

\end{proposition}

\subsection{Equivalence of Grothendieck topologies}

Sometimes two different to\-pol\-ogies on the same category define the same sheaves.

\begin{definition}
Let $\cC$ be a category, $\{U_i \to U\}_{i \in I}$ a set of arrows. A \emph{refinement}\index{refinement} $\{V_a \arr U\}_{a \in A}$ is a set of arrows such that for each index $a \in A$ there is some index $i \in I$ such that $V_a \arr U$ factors through $U_i \arr U$.
\end{definition}

Notice that the choice of factorizations $V_a \arr U_i \arr U$ is \emph{not} part of the data, we simply require their existence.

This relation between sets of arrows is most easily expressed in terms of sieves. The following fact is immediate.

\begin{proposition}\label{prop:funnychar-refinement}\index{refinement!characterization via sieves}
Let there be given two sets of arrows $\cU = \{U_i \to U\}$ and $\cV = \{V_a \arr U\}$. Then $\cV$ is a refinement of $\cU$ if and only if $\h_{\cV} \subseteq \h_{\cU}$.
\end{proposition}

A refinement of a refinement is obviously a refinement. Also, any covering is a refinement of itself: thus, the relation of being a refinement is a pre-order on the set of coverings of an object $U$.

\begin{definition}\label{def:same-sheaves}
Let $\cC$ be a category, $\cT$ and $\cT'$ two topologies on $\cC$. We say that $\cT$ is \emph{subordinate}\index{topology!subordinate}\index{subordinate topology} to $\cT'$, and write $\cT \prec \cT'$, if every covering in $\cT$ has a refinement that is a covering in $\cT'$.

If $\cT \prec \cT'$ and $\cT' \prec \cT$, we say that $\cT$ and $\cT'$ are \emph{equivalent}\index{topology!equivalent}\index{equivalent topologies}, and write $\cT \equiv \cT'$.
\end{definition}

Being a refinement is a relation between sets of arrows into $U$ that is transitive and reflexive. Therefore being subordinate is a transitive and reflexive relation between topologies on $\cC$, and being equivalent is an equivalence relation.

This relation between topologies is naturally expressed in terms of sieves.

\begin{proposition}\label{prop:funnychar-subordinate}\index{topology!subordinate!characterization via sieves}\index{subordinate topology!characterization via sieves}
Let $\cT$ and $\cT'$ be topologies on a category $\cC$. Then $\cT \prec \cT'$ if and only if every sieve belonging to $\cT$ also belongs to $\cT'$.

In particular, two topologies are equivalent if and only if they have the same sieves.
\end{proposition}

This is clear from Proposition~\ref{prop:funnychar-refinement}.

\begin{proposition}\label{prop:same-sheaves}
\index{topology!equivalent!has the same sheaves}
\index{equivalent topologies!have the same sheaves}
Let $\mathcal T$ and $\cT'$ be two Grothendieck topologies on the same category $\cC$. If $\cT$ is subordinate to $\cT'$, then every sheaf in $\cT'$ is also a sheaf in $\cT$.

In particular, two equivalent topologies have the same sheaves.
\end{proposition}

The proof is immediate from Propositions \ref{prop:funnychar-sheaves} and \ref{prop:funnychar-subordinate}.

In Grothendieck's language what we have defined would be called a pretopology, and two equivalent pretopologies define the same topology.

\begin{example}
The global classical topology on $\cattop$ (Example~\ref{ex:global-classical}), and the global \'etale topology of Example~\ref{ex:top-global-etale}, are equivalent.
\end{example}

\begin{example}
If $S$ is a base scheme, there is another topology that we can define over the category $\catsch{S}$, the \emph{smooth topology}\index{topology!smooth}\index{smooth!topology}, in which a covering $\{U_{i}\to U\}$ is a jointly surjective set of smooth morphisms locally of finite presentation.

By \cite[Corollaire~17.16.3]{ega4-4}, given a smooth covering $\{U_{i}\to U\}$ we can find an \'etale surjective morphism $V \arr U$ that factors through the disjoint union $\coprod_{i}U_{i} \arr U$; given such a factorization, if $V_{i}$ is the inverse image of $U_{i}$ in $V$, we have that $\{V_{i} \to U\}$ is an \'etale covering that is a refinement of $\{U_{i}\to U\}$. This means that the smooth topology is subordinate to the \'etale topology. Since obviously every \'etale covering is a smooth cover, the two topologies are equivalent.
\end{example}

\begin{definition}\label{def:saturation}
A topology $\cT$ on a category $\cC$ is called \emph{saturated}\index{topology!saturated}\index{saturated topology} if a set of arrows $\{U_{i} \arr U\}$ which has a refinement that is in $\cT$ is also in $\cT$.

If $\cT$ is a topology of $\cC$, the \emph{saturation}\index{topology!saturation of a}\index{saturation of a topology} $\overline{\cT}$ of $\cT$ is the set of sets of arrows which have a refinement in $\cT$.
\end{definition}

\begin{proposition}
\index{topology!saturation of a! properties of}
\index{saturation of a topology!properties of}
Let $\cT$ be a topology on a category $\cC$.

\begin{enumeratei}

\item The saturation $\overline{\cT}$ of $\cT$ is a saturated topology.

\item $\cT \subseteq \overline{\cT}$.

\item $\overline{\cT}$ is equivalent to $\cT$.

\item The topology $\cT$ is saturated if and only if $\cT = \overline{\cT}$.

\item A topology $\cT'$ on $\cC$ is subordinate to $\cT$ if and only if $\cT' \subseteq \overline{\cT}$.

\item A topology $\cT'$ on $\cC$ is equivalent to $\cT$ if and only if $\overline{\cT'} =
\overline{\cT}$.

\item A topology on $\cC$ is equivalent to a unique saturated topology.

\end{enumeratei}
\end{proposition}

We leave the easy proofs to the reader.

\subsection{Sheaf conditions on representable functors}

\begin{proposition}\label{prop:rep->sheaf}
\index{representable!functor!is a sheaf in the global classical topology}
\index{functor!representable!is a sheaf in the global classical topology}
A representable functor $\cattop\op \arr \catset$ is a sheaf in the global classical topology.
\end{proposition}

This amounts to saying that, given two topological spaces $U$ and $X$, an open covering $\{U_i \subseteq U\}$, and continuous functions $f_i \colon U_i \arr X$, with the property that the restriction of $f_i$ and $f_j$ to $U_i \cap U_j$ coincide for all $i$ and $j$, there exists a unique continuous function $U \arr X$ whose restriction $U_{i} \arr X$ is $f_{i}$. This is essentially obvious (it boils down to the fact that, for a function, the property of being continuous is local on the domain). For similar reasons, it is easy to show that a representable functor on the category $\catsch{S}$ over a base scheme $S$ is a sheaf in the Zariski topology.

On the other hand the following is not easy at all: a scheme is a topological space, together with a sheaf of rings in the Zariski topology. A priori, there does not seem to be a reason why we should be able to glue morphisms of schemes in a finer topology than the Zariski topology.

\begin{theorem}[Grothendieck]\label{thm:rep-fppf}
\index{representable!functor!is a sheaf in the fpqc topology}
\index{functor!representable!is a sheaf in the fpqc topology}
\index{theorem!a representable functor is a sheaf in the fpqc topology}
A representable functor on $\catsch{S}$ is a sheaf in the fpqc topology.
\end{theorem}

So, in particular, it is also a sheaf in the \'etale and in the fppf topologies.

Here is another way of expressing this result. Recall that in a category $\cC$ an arrow $f \colon V \arr U$ is called an \emph{epimorphism}\index{epimorphism} if, whenever we have two arrows $U \double X$ with the property that the two composites $V \arr U \double X$ coincide, then the two arrows are equal. In other words, we require that the function $\hom_{\cC}(U,X) \arr \hom_{\cC}(V,X)$ be injective for any object $X$ of $\cC$.

On the other hand, $V \arr U$ is called an \emph{effective epimorphism}\index{epimorphism!effective}\index{effective!epimorphism} if for any object $X$ of $\cC$, any arrow $V \arr X$ with the property that the two composites
   \[
   V \times_U V \doublelong{\pr_1}{\pr_2} V \arr X
   \]
coincide, factors uniquely through $U$. In other words, we require that the diagram
   \[
   \hom_{\cC}(U,X) \arr \hom_{\cC}(V,X)
   \doublelong{}{} \hom_{\cC}(V \times_U V,X)
   \]
be an equalizer.

Then Theorem~\ref{thm:rep-fppf} says that every fpqc morphism of schemes is an effective epimorphism in $\catsch{S}$.

\begin{remark}\label{rmk:really-need-finiteness}
\index{wild flat topology}
\index{topology!wild flat}
As we have already observed at the beginning of  \S\ref{subsec:fpqc}, there is a ``wild'' flat topology in which the coverings are  jointly surjective sets $\{U_i \to U\}$ of flat morphisms. However, this topology is very badly behaved; in particular, not all representable functors are sheaves.

Take an integral smooth curve $U$ over an algebraically closed field, with quotient field $K$ and let $V_p = \spec\cO_{U,p}$ for all closed points $p \in U(k)$, as in Remark~\ref{rmk:need-finiteness}. Then $\{V_p \arr U\}$ is a covering in this wild flat topology.

Each $V_p$ contains the closed point $p$, and $V_p \setminus \{p\} = \spec K$ is the generic point of $V_p$; furthermore $V_p \times_U V_q = V_p$ if $p = q$, otherwise $V_p \times_U V_q = \spec K$.

We can form a (very non-separated) scheme $X$ by gluing together all the $V_p$ along $\spec K$; then the embeddings $V_p \into X$ and $V_q \into X$ agree when restricted to $V_p \times_U V_q$, so the give an element of $\prod_{p}\h_X V_p$ whose two images in $\prod_{p,q} \h_X(V_p \times_U V_q)$ agree. However, there is no morphism $U \arr X$ whose restriction to each $V_p$ is the natural morphism $V_p \arr U$. In fact, such a morphism would have to send each closed point $p \in U$ into $p \in V_p \subseteq X$, and the generic point to the generic point; but the resulting set-theoretic function $U \arr X$ is not continuous, since all subsets of $X$ formed by closed points are closed, while only the finite sets are closed in $U$.

\end{remark}

\begin{definition}\label{def:subcanonical}
A topology $\cT$ on a category $\cC$ is called \emph{subcanonical}\index{topology!subcanonical}\index{subcanonical!topology} if every representable functor on $\cC$ is a sheaf with respect to $\cT$.

A \emph{subcanonical site}\index{site!subcanonical}\index{subcanonical!site} is a category endowed with a subcanonical topology.
\end{definition}

There are examples of sites that are not subcanonical (we have just seen one in Remark~\ref{rmk:really-need-finiteness}), but I have never had dealings with any of them.

The name ``subcanonical'' comes from the fact that on a category $\cC$ there is a topology, known as the \emph{canonical topology}\index{topology!canonical}\index{canonical!topology}, which is the finest topology in which every representable functor is a sheaf. We will not need this fact.

\begin{definition}\label{def:comma-topology}
Let $\cC$ be a site, $S$ an object of $\cC$. We define the \emph{comma topology} \index{topology!comma}\index{comma!topology} on the comma category $(\cC/S)$ as the topology in which a covering of  an object $U \arr S$ of $(\cC/S)$ is a collection of arrows
   \[
   \xymatrix@-15pt{
   U_i\ar[rr]^{f_i}\ar[rd] && U\ar[ld]\\
   &S
   }
   \]
such that the collection $\{f_i \colon U_i \arr U\}$ is a covering in $\cC$. In other words, the coverings of $U \arr S$ are simply the coverings of $U$.
\end{definition}

It is very easy to check that the comma topology is in fact a topology.

For example, if $\cC$ is the category of all schemes (or, equivalently, the category of schemes over $\ZZ$), then $(\cC/S)$ is the category of schemes over $S$, and the comma topology induced by the fpqc topology on $(\cC/S)$ is the fpqc topology. Analogous statements hold for the Zariski, \'etale and fppf topology.

\begin{proposition}\label{prop:comma-subcanonical}
If $\cC$ is a subcanonical site and $S$ is an object of $\cC$, then $(\cC/S)$ is also subcanonical.
\end{proposition}

\begin{proof}

We need to show that for any covering $\{U_i \to U\}$ in $(\cC/S)$ the sequence
   \[
   \hom_S(U,X) \arr \prod_i\hom_S(U_i, X)
   \doublelong{}{} \prod _{i,j}\hom_S(U_i \times_U U_j, X)
   \]
is an equalizer. The injectivity of the function
   \[
   \hom_S(U,X) \arr \prod_i\hom_S(U_i, X)
   \]
is clear, since $\hom_S(U,X)$ injects into $\hom(U,X)$, $\prod_i\hom_S(U_i, X)$ injects into $\prod_i\hom(U_i, X)$, and  $\hom(U,X)$ injects into $\prod_i\hom(U_i, X)$, because $\h_{X}$ is a sheaf. On the other hand, let us suppose that we are given an element $(a_i)$ of the product $\prod_i\hom_S(U_i, X)$, with the property that for all pairs $i$, $j$ of indices the equality $\pr_1^*a_i = \pr_2^*a_j$ holds in $\hom_S(U_i \times_U U_j, X)$. Then there exists a morphism $a \in \hom(U, X)$ such that the composite $U_i \arr U \xrightarrow{a} X$ coincides with $a_i$ for all, and we only have to check that $a$ is a morphism of $S$-objects. But the composite $U_i \arr U \xrightarrow{a} X \arr S$ coincides with the structure morphism $U_i \arr S$ for all $i$; since $\hom(-,S)$ is a sheaf on the category $\cC$, so that $\hom(U,S)$ injects into $\prod_i\hom(U_i,S)$, this implies that the composite $U \xrightarrow{a} X \arr S$ is the structure morphism of $U$, and this completes the proof.
\end{proof}

\begin{proof}[Proof of Theorem~\ref{thm:rep-fppf}]

We will use the following useful criterion.

\begin{lemma}\label{lem:criterion-sheaf}
\index{sheaf!characterization of fpqc}
\index{fpqc!characterization of fpqc sheaves}
Let $S$ be a scheme, $F \colon \catsch{S}\op \arr \catset$ a functor. Suppose that $F$ satisfies the following two conditions.
\begin{enumeratei}

\item $F$ is a sheaf in the global Zariski topology.

\item Whenever $V \arr U$ is a faithfully flat morphism of affine $S$-schemes, the diagram
   \[
   FU\arr FV\doublelong{\pr_1^*}{\pr_2^*}
   F(V \times_U V)
   \]
is an equalizer.
\end{enumeratei}

Then $F$ is a sheaf in the fpqc topology.
\end{lemma}

\begin{proof}The proof will be divided into several steps.

\steps

\step[: reduction to the case of a single morphism]

Take a covering $\{U_i \to U\}$ of schemes over $S$ in the fpqc topology, and set $V = \coprod_i U_i$. The induced morphism $V \arr U$ is fpqc. Since $F$ is a Zariski sheaf, the function $FV \arr \prod_i FU_i$ induced by restrictions is an isomorphism. We have a commutative diagram of sets
   \[
   \xymatrix{
   FU\ar[r]\ar@{=}[d] & FV\ar[d]\ar@<3pt>[r]^-{\pr_1^*}
   \ar@<-3pt>[r]_-{\pr_2^*} & F(V \times_U V)\ar[d]\\
   FU\ar[r] & {}\prod_i FU_i\ar@<3pt>[r]^-{\pr_1^*}
   \ar@<-3pt>[r]_-{\pr_2^*} & {}\prod_{i,j}F(U_i \times_U U_j)
   }
   \]
where the columns are bijections; hence to show that the bottom row is an equalizer it is enough to show that the top row is an equalizer. In other words, we have shown that it is enough to consider coverings $\{V \arr U\}$ consisting of a single morphism. Similarly, to check that $F$ is separated we may limit ourselves to considering coverings consisting of a single morphism.

This argument also shows that if $\{U_i \to U\}$ is a finite covering, such that $U$ and the $U_i$ are affine, then the diagram
   \[
   FU\arr \prod_i
   FU_i\doublelong{\pr_1^*}{\pr_2^*}
   \prod_{i,j} F(U_i \times_U U_j)
   \]
is an equalizer. In fact, in this case the finite disjoint union $\coprod_i U_i$ is also affine.

\step[: proof that $F$ is separated]

Now we are given an fpqc morphism $f \colon V \arr U$; take an open covering $\{V_i\}$ of $V$ by open quasi-compact subsets, whose image $U_i = fV_i$ is open and affine. Write each $V_i$ as a union of finitely many affine open subschemes $V_{ia}$. Consider the commutative diagram of restriction functions
   \[
   \xymatrix@C-5pt{
   FU\ar[r]\ar[d]
   & FV \ar[d]\\
   {}\prod_i FU_i\ar[r]
   & {}\prod_i \prod_a FV_{ia}\hsmash{.}
   }
   \]

Its columns are injective, because $F$ is a sheaf in the Zariski topology. On the other hand, the second row is also injective, because each of the restriction morphisms $FU_i \arr \prod_{a,b}F(V_{ia} \times_U V_{ib})$
is injective. Hence the restriction function $FU \arr FV$ is injective, so $F$ is separated.

\step[: the case of a morphism from a quasi-compact scheme onto an affine scheme] Let $f\colon V \arr U$ a faithfully flat morphism, with $V$ quasi-compact and $U$ affine. Let $b\in FV$ be an element such that
   \[
   \pr_{1}^{*}b = \pr_{2}^{*}b \in F(V \times_{U} V).
   \]
We need to show that there exists an element $a \in FU$ such that $f^{*}a = b \in FV$.

Let ${V_{i}}$ be a finite covering of $V$ by open affine subschemes; then $\{V_{i} \arr U\}$ is a finite fpqc covering of $U$ by affine subschemes, hence the sequence
   \[
   FU \arr \prod_{i}FV_{i} \doublelong{\pr_{1}^{*}}{\pr_{2}^{*}} 
   \prod_{i,j}F(V_{i}\times_{U}V_{j})
   \]
is an equalizer.

For each $i$ denote by $b_{i}$ the restriction of $b$ to $FV_{i}$; then $\pr_{1}^{*} b_{i} \in F(V_{i}\times_{U}V_{j})$ is the restriction of $\pr_{1}^{*}b \in FV$ to $F(V_{i}\times_{U}V_{j})$, while $\pr_{2}^{*} b_{j} \in F(V_{i}\times_{U}V_{j})$ is the restriction of $\pr_{2}^{*}b \in FV$ to $F(V_{i}\times_{U}V_{j})$. Hence $\pr_{1}^{*} b_{i} = \pr_{2}^{*}b$ for all $i$ and $j$, so there exists some $a \in FU$ whose pullback to $FV_{i}$ is $b_{i}$ for all $i$. Then the restrictions of $f^{*}a$ and $b$ to $FV_{i}$ coincide for all $i$, so $f^{*}a = b$, because $F$ is a sheaf in the Zariski topology.

\step[: the case of a morphism to an affine scheme]
Let $V \arr U$ be an fpqc scheme, where $U$ is affine. Let $f\colon V \arr U$ a faithfully flat morphism, with $V$ quasi-compact and $U$ affine. Let $b\in FV$ be an element such that
   \[
   \pr_{1}^{*}b = \pr_{2}^{*}b \in F(V \times_{U} V).
   \]
We need to show that there exists an element $a \in FU$ such that $f^{*}a = b \in FV$.

Let $\{V_{i}\}$ be an open covering of $V$ by quasi-compact open subschemes, such that the projection $V_{i} \arr U$ is surjective for all $i$. For each $i$, denote by $b_{i}$ the restriction of $b$ to $FV_{i}$. The restriction morphism $f\mid_{V_{i}}\colon V_{i} \arr V$ is fpqc, hence by the previous step there exists $a_{i} \in FU$ such that $(f\mid_{V_{i}})^{*}a_{i} = b_{i}$. However, I claim that $a_{i} = a_{j}$ for all $i$ and $j$. In fact, the morphism $V_{i} \cup V_{j} \arr U$ is also fqpc, and $V_{i} \cup V_{j}$ is quasi-compact: hence there exists $a_{ij}$ in $FU$ whose pullback to $F(V_{i} \cup V_{j})$ is the restriction of $b$. Since the pullbacks of $a_{ij}$ to $FV_{i}$ and $FV_{j}$ coincide with $b_{i}$ and $b_{j}$ respectively, we have that $a_{i} = a_{ij} = a_{j}$.

Hence the pullback of $a \eqdef a_{i}$ to $FV_{i}$ is $b_{i}$ for all $i$; it follows $f^{*}a$ is $b$, as desired.

\step[: the general case]
Now $f\colon V \arr U$ is an arbitrary fpqc morphism. Let $\{U_{i}\}$ be an covering of $U$ by open affine subschemes, and denote by $V_{i}$ the inverse image of $U_{i}$ in $V$. We have a diagram of restriction functions
   \[
   \xymatrix@C-5pt{
   FU\ar[r]\ar[d]
   & FV \ar@<3pt>[r]\ar@<-3pt>[r]\ar[d]
   & F(V \times_U V)\ar[d]\\
   {}\prod_i FU_i\ar[r]\ar@<3pt>[d]\ar@<-3pt>[d]
   & {}\prod_i FV_{i}\ar@<3pt>[r]\ar@<-3pt>[r]
   \ar@<3pt>[d]\ar@<-3pt>[d]
   & {}\prod_i
   F(V_{i} \times_{U_{i}} V_{i})\\
   {}\prod_{i,j} F(U_i \cap U_j)\ar[r]
   & {}\prod_{i,j}
   F(V_{i}\cap V_{j})\hsmash{.}
   }
   \]
The columns are equalizers, because $F$ is a sheaf in the Zariski topology; furthermore the second row is also an equalizers, because each of the diagrams
   \[
   FU_i \arr FV_{i} \doublelong{}{} F(V_{i} \times_{U_{i}} V_{i})
   \]
is an equalizer, by the previous step, and the product of equalizers is an equalizer. Finally, the bottom row is injective, because $F$ is separated, and the result follows from a simple diagram chasing.
\end{proof}

To prove Theorem~\ref{thm:rep-fppf} we need to check that if $F = \h_X$, where $X$ is an $S$-scheme, then the second condition of Lemma~\ref{lem:criterion-sheaf} is satisfied. First of all, by Proposition~\ref{prop:comma-subcanonical} it is enough to prove the result in case $S = \spec \ZZ$, that is, when $\catsch{S}$ is simply the category of all schemes. So for the rest of the proof we only need to work with morphism of schemes, without worrying about base schemes.

We will assume at first that $X$ is affine. Set $U = \spec A$, $V = \spec B$, $X = \spec R$. In this case the result is an easy consequence of the following lemma. Consider the ring homomorphism $f \colon A \arr B$ corresponding to the morphism $V \arr U$, and the two homomorphisms of $A$-algebras $e_1, e_2\colon B \arr B \otimes_A B$ defined by $e_1(b) = b \otimes 1$ and $e_2(b) = 1 \otimes b$; these correspond to the two projections $V \times_U V \arr V$.

\begin{lemma}\label{lem:exact-affine}
The sequence
   \[
   0 \arr A \overset{f}\larr B
   \xrightarrow{e_1 - e_2} B \otimes_A B
   \]
is exact.
\end{lemma}
\begin{proof}

The injectivity of $f$ is clear, because $B$ is faithfully flat over $A$. Also, it is clear that the image of $f$ is contained in the kernel of $e_1 - e_2$, so we have only to show that the kernel of $e_1 - e_2$ is contained in the image of $f$.

Assume that there exists a homomorphism of $A$-algebras $g \colon B \arr A$ (in other words, assume that the morphism $V \arr U$ has a section). Then the composite $g \circ f \colon A\arr A$ is the identity. Take an element $b \in \ker (e_1 - e_2)$; by definition, this means that $b \otimes 1 = 1 \otimes b$ in $B \otimes_A B$. By applying the homomorphism $g \otimes \id_B \colon B \otimes_A B \arr A \otimes_A B = B$ to both members of the equality we obtain that $f(g b) = b$, hence $b \in \im f$.

In general, there will be no section $U \arr V$; however, suppose that there exists a faithfully flat $A$ algebra $A \arr A'$, such that the homomorphism $f \otimes \id_{A'} \colon A' \arr B \otimes A'$ obtained by base change has a section $B \otimes A' \arr A'$ as before. Set $B' = B \otimes A'$. Then there is a natural isomorphism of $A'$-algebras $B' \otimes_{A'} B' \simeq (B \otimes_A B) \otimes_A A'$, making the diagram
   \[
   \xymatrix@C+5pt{
   0\ar[r] & A'\ar[r]^-{f \otimes \id_{A'}}\ar@{=}[d]
   &B'\ar@{=}[d]\ar[rr]^-{e'_1 - e'_2}
   && B' \otimes_{A'} B'\ar[d]\\
   0\ar[r] & A'\ar[r]^-{f \otimes \id_{A'}}
   &B'\ar[rr]^-{(e_1 - e_2) \otimes \id_{A'}}
   && (B \otimes_A B) \otimes_A A'
   }
   \]
commutative. The top row is exact, because of the existence of a section, and so the bottom row is exact. The thesis follows, because $A'$ is faithfully flat over $A$.

But to find such homomorphism $A \arr A'$ it is enough to set $A' = B$; the product $B \otimes_A A' \arr A'$ defined by $b \otimes b' \mapsto bb'$ gives the desired section. In geometric terms, the diagonal $V \arr V \times_U V$ gives a section of the first projection $V \times_U V \arr V$.
\end{proof}

To finish the proof of Theorem~\ref{thm:rep-fppf} in the case that $X$ is affine, recall that morphisms of schemes $U \arr X$, $V \arr X$ and $V \times_U V \arr X$ correspond to ring homomorphisms $R \arr A$, $R \arr B$ and $R \arr B \otimes_A B$; then the result is immediate from the lemma above. This proves that $\h_X$ is a sheaf when $X$ is affine.

If $X$ is not necessarily affine, write $X = \cup_i X_i$ as a union of affine open subschemes.

Let us show that $\h_X$ is separated.  Given a covering $V \arr U$, take two morphisms $f, g \colon U \arr X$ such that the two composites $V \arr U \arr X$ are equal. Since $V \arr U$ is surjective, $f$ and $g$ coincide set-theoretically, so we can set $U_i = f^{-1}X_i = g^{-1}X_i$, and call $V_i$ the inverse images of $U_i$ in $V$. The two composites
   \[
   V_i \arr U_i
   \doublelong{f\mid_{U_i}}{g\mid_{U_i}} X_i
   \]
coincide, and $X_i$ is affine; hence $f\mid_{U_i} = g\mid_{U_i}$ for all $i$, so $f = g$, as desired.

To complete the proof, suppose that $g \colon V \arr X$ is a morphism with the property that the two composites
   \[
   V \times_U V \doublelong{\pr_1}{\pr_2} V 
   \overset{g} \larr X
   \]
are equal; we need to show that $g$ factors through $U$. The morphism $V \arr U$ is surjective, so, from Lemma~\ref{lem:product-set} below, $g$ factors through $U$ set-theoretically. Since $U$ has the quotient topology induced by the morphism $V \arr U$ (Proposition~\ref{prop:flat->quotient-topology}), we get that the resulting function $f \colon U \arr X$ is continuous.

Set $U_i = f ^{-1} X_i$ and $V_i = g ^{-1}V_i$ for all $i$. The composites
   \[
   V_i \times_U V_i \doublelong{\pr_1}{\pr_2} V_i 
   \overset{g\mid_{V_i}}\larr
   V_i \arr X_i
   \]
coincide, and $X_i$ is affine, so $g\mid_{V_i} \colon V_i \arr X$ factors uniquely through a morphism $f_i \colon U_i \arr X_i$. We have 
   \[
   f_i\mid_{U_i \cap U_j} = f_j\mid_{U_i \cap U_j} \colon 
   U_i \cap U_j \arr X,
   \]
because $\h_X$ is separated; hence the $f_i$ glue together to give the desired factorization $V \arr U \arr X$.
\end{proof}

\begin{lemma}\label{lem:product-set}
Let $f_1 \colon X_1 \arr Y$ and $f_2 \colon X_2 \arr Y$ be morphisms of schemes. If $x_1$ and $x_2$ are points  of $X_1$ and $X_2$ respectively, and $f_1(x_1) = f_2(x_2)$, then there exists a point $z$ in the fibered product $X_1 \times_{Y} X_2$ such that $\pr_1(z) = x_1$ and $\pr_2(z) = x_2$.
\end{lemma}

\begin{proof}

Set $y = f(x_1) = f_2(x_2) \in Y$. Consider the extensions $k(y) \subseteq k(x_1)$ and $k(y) \subseteq k(x_2)$; the tensor product $k(x_1) \otimes_{k(y)} k(x_2)$ is not 0, because the tensor product of two vector spaces over a field is never $0$, unless one of the vector spaces is $0$. Hence $k(x_1) \otimes_{k(y)} k(x_2)$ has a maximal ideal; the quotient field $K$ is an extension of $k(y)$ containing both $k(x_1)$ and $k(x_2)$. The two composites $\spec K \arr \spec k(x_1) \arr X_1 \xrightarrow{f_1} U$ and $\spec K \arr \spec k(x_2)\arr X_2 \xrightarrow{f_2} U$ coincide, so we get a morphism $\spec K \arr X_{1} \times_Y X_{2}$. We take $z$ to be the image of $\spec K$ in $X_1 \times_{Y} X_2$.
\end{proof}

The proof of Theorem~\ref{thm:rep-fppf} is now complete.

\subsection{The sheafification of a functor}

The usual construction of the sheafification of a presheaf of sets on a topological space carries over to this more general context.

\begin{definition}\label{def:sheafification}
Let $\cC$ be a site, $F \colon \cC\op \arr \catset$ a functor. A \emph{sheafification}\index{sheafification} of $F$ is a sheaf $F\ass \colon \cC\op \arr \catset$, together with a natural transformation $F \arr F\ass$, such that:

\begin{enumeratei}

\item given an object $U$ of $\cC$ and two elements $\xi$ and $\eta$ of $FU$ whose images $\xi\ass$ and $\eta\ass$ in  $F\ass U$ are the same, there exists a covering $\{\sigma_i\colon U_i \arr U\}$ such that $\sigma_i^* \xi = \sigma_i^* \eta$, and

\item for each object $U$ of $\cC$ and each $\overline \xi \in F\ass (U)$, there exists a covering $\{\sigma_i\colon U_i \arr U\}$ and elements $\xi_i \in F(U_i)$ such that $\xi_i\ass = \sigma_i^* \overline{\xi}$.

\end{enumeratei}
\end{definition}

\begin{theorem}\call{thm:sheafification}
\index{sheafification!existence of}
\index{theorem!existence of sheafifications}
Let $\cC$ be a site, $F \colon \cC\op \arr \catset$ a functor.

\begin{enumeratei}

\itemref{1} If $F\ass \colon \cC\op \arr \catset$ is a sheafification of $F$, any morphism from $F$ to a sheaf factors uniquely through $F\ass$.

\itemref{2} There exists a sheafification $F \arr F\ass$, which is unique up to a canonical isomorphism.

\itemref{3} The natural transformation $F \arr F\ass$ is injective (that is, each function $FU \arr F\ass U$ is injective) if and only if $F$ is separated.

\end{enumeratei}
\end{theorem}

\begin{proof}[Sketch of proof]

For part~\refpart{thm:sheafification}{1}, we leave to the reader to check uniqueness of the factorization.

For existence, let $\phi \colon F \arr G$ be a natural transformation from $F$ to a sheaf $G \colon  \cC\op \arr \catset$. Given an element $\overline{\xi}$ of $F\ass U$, we want to define the image of $\xi$ in $GU$. There exists a covering $\{\sigma_i \colon U_i \arr U\}$ and elements $\xi_i \in FU_i$, such that the image of $\xi_i$ in $F\ass U_i$ is $\sigma_i^* \overline{\xi}$.  Set $\eta_i= \phi\xi_i \in GU_{i}$. The pullbacks $\pr_1^* \eta_i$ and $\pr_2^* \eta_j$ in $GU_{ij}$ both have as their image in $F\ass U_{ij}$ the pullback of $\overline{\xi}$; hence there is a covering $\{U_{ij \alpha} \arr U _{ij}\}$ such that the pullbacks of $\pr_1^* \xi_i$ and $\pr_2^* \xi_j$ in $F U_{ij\alpha}$ coincide for each $\alpha$. By applying $\phi$, and keeping in mind that it is a natural transformation, and that $G$ is a sheaf, we see that the pullbacks of $\eta_i$ and $\eta_j$ to $G U_{ij}$ are the same, for any pair of indices $i$ and $j$. Hence there an element $\eta$ of $GU$ whose pullback to each $GU_i$ is $\eta_i$. 

We leave to the reader to verify that this $\eta \in GU$ only depends on $\xi$, and that by sending each $\xi$ to the corresponding $\eta$ we define a natural transformation $F\ass \arr G$, whose composition with the given morphism $F \arr F\ass$ is $\phi$.

Let us prove part~\refpart{thm:sheafification}{2}. For each object $U$ of $\cC$, we define an equivalence relation $\sim$ on $F U$ as follows. Given two elements $a$ and $b$ of $F U$, we write $a \sim b$ if there is a covering $U_i \arr U$ such that the pullbacks of $a$ and $b$ to each $ U_i$ coincide. We check easily that this is an equivalence relation, and we define $F\sep U  = F U/\sim$. We also verify that if $V \arr U$ is an arrow in $\cC$, the pullback $F U \arr F V$ is compatible with the equivalence relations, yielding a pullback $F\sep U \arr F\sep V$. This defines the functor $F\sep$ with the surjective morphism $F \arr F\sep$. It is straightforward to verify that $F\sep$ is separated, and that every natural transformation from $F$ to a separated functor factors uniquely through $F\sep$.

To construct $F\ass$, we  take for each object $U$ of $\cC$ the set of pairs $(\{U_i \to U\}, \{a_i\})$, where $\{U_i \to U\}$ is a covering, and $\{a_i\}$ is a set of elements with $a_i \in F\sep U_i$ such that the pullback of $a_i$ and $a_j$ to $F\sep(U_i \times_U U_j)$, along the first and second projection respectively, coincide. On this set we impose an equivalence relation, by declaring $(\{U_i \to U\}, \{a_i\})$ to be equivalent to $(\{V_j \to U\}, \{b_j\})$ when the restrictions of $a_i$ and $b_j$ to $F\sep(U_i \times_U V_j)$, along the first and second projection respectively, coincide. To verify the transitivity of this relation we need to use the fact that the functor $F\sep$ is separated.

For each $U$, we denote by $F\ass U$ the set of equivalence classes. If $V \arr U$ is an arrow, we define a function $F\ass U \arr F\ass V$ by associating with the class of a pair $(\{U_i \to U\}, \{a_i\})$ in $F\ass U$ the class of the pair $(\{U_i \times_U V\}, p_i^* a_i)$, where $p_i \colon U_i\times_U V \arr U_i$ is the projection. Once we have checked that this is well defined, we obtain a functor $F\ass \colon \cC\op \arr \catset$. There is also a natural transformation $F\sep \arr F\ass$, obtained by sending an element $a \in F\sep U$ into $(\{U = U\}, a)$. Then one verifies that $F\ass$ is a sheaf, and that the composite of the natural transformations $F \arr F\sep$ and $F\sep \arr F\ass$ has the desired universal property.

The uniqueness up to a canonical isomorphism follows immediately from part~\refpart{thm:sheafification}{1}. Part~\refpart{thm:sheafification}{3} follows easily from the definition.
\end{proof}

A slicker, but equivalent, definition is as follows. Consider the set $\{S_{i}\}$ of sieves belonging to $\cT$ on an object $U$ of $\cC$. These form a ordered set: we set $i \leq j$ if $S_{j} \subseteq S_{i}$. According to Proposition~\ref{prop:directed-sieves}, this is a direct system, that is, given two indices $i$ and $j$ there is some $k$ such that $k \geq i$ and $k \geq j$. Then $F\ass U$ is in a canonical bijective correspondence with the direct limit $\indlim_{i}\hom_{\cC}(S_{i}, F\sep)$.

\chapter{Fibered categories}\label{ch:fibered}

\section{Fibered categories}\label{sec:fibered}

\subsection{Definition and first properties}

In this Section we will fix a category $\cC$; the topology will play no role. We will study  categories over $\cC$, that is, categories $\cF$ equipped with a functor $\p_{\cF}\colon \cF \arr \cC$.

We will draw several commutative diagrams involving objects of $\cC$ and $\cF$; an arrow going from an object $\xi$ of $\cF$ to an object  $U$ of $\cC$ will be of type ``$\xi\mapsto U$'', and will mean that $\p_{\cF} \xi = U$. Furthermore the commutativity of the diagram
   \[
   \xymatrix{{}\xi\ar@{|->}[d]\ar[r]^\phi &{}\eta\ar@{|->}[d] \\
   U\ar[r]^f & V}
   \]
will mean that $\p_{\cF} \phi = f$.

\begin{definition} Let $\cF$ be a category over $\cC$. An arrow  $\phi\colon \xi \arr \eta$ of $\cF$ is \emph{cartesian}\index{arrow!cartesian}\index{cartesian arrow} if for any arrow $\psi\colon \zeta \arr \eta$ in $\cF$ and any arrow $h \colon \p_{\cF} \zeta \arr \p_{\cF} \xi$ in $\cC$ with $\p_{\cF} \phi \circ  h = \p_{\cF} \psi$, there exists a unique arrow $\theta \colon \zeta \arr \xi$ with $\p_{\cF} \theta  = h$ and $\phi\circ \theta = \psi$, as in the commutative diagram
   \[
   \xymatrix@R=6pt{
   {\zeta}\ar@{|->}[dd] \ar@{-->}[rd]_\theta
   \ar@/^/[rrd]^{\psi} \\
   & {}\xi\ar@{|->}[dd]\ar[r]_{\phi}
   & {}\eta\ar@{|->}[dd] \\
   {}\p_{\cF} \zeta\ar[rd]_h\ar@/^/[rrd]|(.487)\hole\\
   & {}\p_{\cF} \xi \ar[r]
   & {}\p_{\cF} \eta \hsmash{.}
   }
   \]

If $\xi \arr \eta$ is a cartesian arrow of $\cF$ mapping to an arrow $U \arr V$ of $\cC$, we also say that $\xi$ is \emph{a pullback of $\eta$ to $U$}\index{pullback}.
\end{definition}

\begin{remark}\label{rmk:unique-pullback}
The definition that we give of cartesian arrow is more restrictive than the definition in \cite{sga1};  our cartesian arrows are called \emph{strongly cartesian} in \cite{gray-fibered}. However, the resulting notions of fibered category coincide.
\end{remark}

\begin{remark}
Given two pullbacks $\phi\colon \xi \arr \eta$ and $\widetilde\phi \colon \widetilde\xi \arr \eta $ of $\eta$ to $U$, the unique arrow $\theta\colon \widetilde \xi \arr \xi$ that fits into the diagram
   
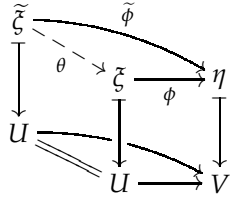
\[
   \xymatrix@R=6pt{
   {\widetilde\xi}\ar@{|->}[dd]\ar@{-->}[rd]_\theta
   \ar@/^/[rrd]^{\widetilde\phi} \\
   & {}\xi\ar@{|->}[dd]\ar[r]_{\phi}
   &{}\eta\ar@{|->}[dd] \\
   U\ar@{=}[rd]\ar@/^/[rrd]|(.482)\hole\\
   & U \ar[r] & V
   }
   \]
is an isomorphism; the inverse is the arrow $\xi \arr \widetilde{\xi}$ obtained by exchanging $\xi$ and $\widetilde{\xi}$ in the diagram above.

In other words, a pullback is unique, up to a unique isomorphism.
\end{remark}

The following facts are easy to prove, and are left to the reader.

\begin{proposition}\call{prop:comp-cartesian}
\index{arrow!cartesian!properties}
\index{cartesian arrow!properties}
\hfil
\begin{enumeratei}

\itemref{1} If $\cF$ is a category over $\cC$, the composite of cartesian arrows in $\cF$ is cartesian.

\itemref{4} If $\xi \arr \eta$ and $\eta \arr \zeta$ are arrows in $\cF$ and $\eta \arr\zeta$ is cartesian, then $\xi \arr\eta$ is cartesian if and only if the composite $\xi \arr \zeta$ is cartesian.

\itemref{3} An arrow in $\cF$ whose image in $\cC$ is an isomorphism is cartesian if and only if it is an isomorphism.

\itemref{2} Let $\p_{\cG}\colon \cG \arr \cC$ and $F \colon \cF \arr \cG$ be functors, $\phi \colon \xi \arr \eta$ an arrow in $\cF$. If $\phi$ is cartesian over its image $F \phi \colon F \xi \arr F \eta$ in $\cG$ and $F \phi$ is cartesian over its image $\p_{\cG}F\phi \colon \p_{\cG}F\xi \arr \p_{\cG}F \eta$ in $\cC$, then $\phi$ is cartesian over its image $\p_{\cG}F\phi$ in $\cC$.
\end{enumeratei}
\end{proposition}

\begin{definition} A \emph{fibered category over $\cC$}\index{category!fibered}\index{fibered!category} is a category $\cF$ over $\cC$, such that given an arrow $f \colon U \arr V$ in $\cC$ and an object $\eta$ of $\cF$ mapping to $V$, there is a cartesian arrow $\phi\colon \xi \arr \eta$ with $\p_{\cF} \phi = f$.
\end{definition}

In other words, in a fibered category $\cF \arr \cC$  we can pull back objects of $\cF$ along any arrow of $\cC$.

\begin{definition}
If $\cF$ and $\cG$ are fibered categories over $\cC$, then a \emph{morphism of fibered categories}\index{category!fibered!morphism of}\index{morphism!of fibered categories} $F \colon \cF \arr \cG$ is a functor such that:

\begin{enumeratei}
\item $F$ is \emph{base-preserving}\index{functor!base-preserving}\index{base-preserving!functor}, that is, $\p_{\cG} \circ F = \p_{\cF}$; \item $F$ sends cartesian arrows to cartesian arrows.
\end{enumeratei}
\end{definition}

Notice that in the definition above the equality $\p_{\cG} \circ F = \p_{\cF}$ must be interpreted as an actual equality. In other words, the existence of an isomorphism of functors between $\p_{\cG} \circ F$ and $\p_{\cF}$ is not enough.
 
\begin{proposition}\call{prop:comp-fibered}
Let there be given two functors $\cF \arr \cG$ and $\cG \arr \cC$. If $\cF$ is fibered over $\cG$ and $\cG$ is fibered over $\cC$, then $\cF$ is fibered over $\cC$.
\end{proposition}

\begin{proof}
This follows from Proposition~\refall{prop:comp-cartesian}{2}.
\end{proof}

\subsection{Fibered categories as pseudo-functors}

\begin{definition} Let $\cF$ be a fibered category over $\cC$. Given an object $U$ of $\cC$, the \emph{fiber}\index{category!fibered!fiber of}\index{fiber of a fibered category} $\cF(U)$ of $\cF$ over $U$ is the subcategory of $\cF$ whose objects are the objects $\xi$ of $\cF$ with $\p_{\cF} \xi = U$, and whose arrows are arrows $\phi$ in $\cF$ with $\p_{\cF} \phi = \id_U$.
\end{definition}

By definition, if $F \colon \cF \arr \cG$ is a morphism of fibered categories over $\cC$ and $U$ is an object of $\cC$, the functor $F$ sends $\cF(U)$ to $\cG(U)$, so we have a restriction functor $F_U \colon \cF(U) \arr \cG(U)$.

Notice that formally we could give the same definition of a fiber for any functor $\p_{\cF} \colon \cF \arr \cC$, without assuming that $\cF$ is fibered over $\cC$. However, we would end up with a useless notion. For example, it may very well happen that we have two objects $U$ and $V$ of $\cC$ which are isomorphic, but such that $\cF(U)$ is empty while $\cF(V)$ is not. This kind of pathology does not arise for fibered categories, and here is why.

Let $\cF$ be a category fibered over $\cC$, and $f \colon U \arr V$ an arrow in $\cC$. For each object $\eta$ over $V$, we choose a pullback $\phi_\eta \colon f^*\eta \arr \eta$ of $\eta$ to $U$. We define a functor $f^* \colon \cF(V) \arr \cF(U)$ by sending each object $\eta$ of $\cF(V)$ to $f^*\eta $, and each arrow $\beta\colon \eta \arr \eta'$ of $\cF(V)$ to the unique arrow $f^*\beta \colon f^*\eta \arr f^*\eta'$ in $\cF(U)$ making the diagram
   \[
   \xymatrix{
   f^* \eta\ar[r]\ar@{-->}[d]^{f^* \beta}
   &{}\eta\ar[d]^\beta\\
   f^* \eta'\ar[r]&{}\eta'\\}
   \] commute.

\begin{definition}
A \emph{cleavage}\index{cleavage} of a fibered category $\cF \arr \cC$ consists of a class $K$ of cartesian arrows in $\cF$ such that for each arrow $f \colon U \arr V$ in $\cC$ and each object $\eta$ in $\cF(V)$ there exists a unique arrow in $K$ with target $\eta$ mapping to $f$ in $\cC$.
\end{definition}

By the axiom of choice, every fibered category has a cleavage. Given a fibered category $\cF \arr \cC$ with a cleavage, we associate with each object $U$ of $\cC$ a category $\cF(U)$, and to each arrow $f \colon U \arr V$ a functor $f^* \colon \cF(V) \arr \cF(U)$, constructed as above. It is very tempting to believe that in this way we have defined a functor from $\cC$ to the category of categories; however, this is not quite correct. First of all, pullbacks $\id_U^* \colon \cF(U) \arr \cF(U)$ are not necessarily identities. Of course we could just choose all pullbacks along identities to be identities on the fiber categories: this would certainly work, but it is not very natural, as there are often natural defined pullbacks where this does not happen (in Example~\ref{ex:arrows} and many others). What happens in general is that, when $U$ is an object of $\cC$ and $\xi$ an object of $\cF(U)$, we have the pullback $\isoiden_U(\xi)\colon\id_U^* \xi \arr \xi$ is an isomorphism, because of Proposition~\refall{prop:comp-cartesian}{3}, and this defines an isomorphism  of functors $\isoiden_U\colon \id_U^* \simeq \id_{\cF(U)}$.

A more serious problem is the following. Suppose that we have two arrows $f \colon U \arr V$ and $g \colon V \arr W$ in $\cC$, and an object $\zeta$ of $\cF$ over $W$. Then $f^*g^*\zeta$ is a pullback of $\zeta$ to $U$; however, pullbacks are not unique, so there is no reason why $f^*g^*\zeta$ should coincide with $(gf)^*\zeta$. However, there is a canonical isomorphism $\isoass_{f,g}(\zeta) \colon f^*g^*\zeta \simeq (gf)^*\zeta$ in $\cF(U)$, because both are pullbacks, and this gives an isomorphism $\isoass_{f,g} \colon f^*g^* \simeq (gf)^*$ of functors $\cF(W) \arr \cF(U)$.

So, after choosing a cleavage a fibered category almost gives a functor from $\cC$ to the category of categories, but not quite. The point is that the category of categories is not just a category, but what is known as a 2-category; that is, its arrows are functors, but two functors between the  same two categories in turn form a category, the arrows being natural transformations of functors. Thus there are 1-arrows (functors) between objects (categories), but there are also 2-arrows (natural transformations) between 1-arrows.

What we get instead of a functor is what is called a \emph{pseudo-functor}, or, in a more modern terminology, a \emph{lax 2-functor}.

\begin{definition}\label{def:pseudo-functor}
A \emph{pseudo-functor}\index{pseudo-functor} $\Phi$ on $\cC$ consists of the following data.

\begin{enumeratei}

\item For each object $U$ of $\cC$ a category $\Phi U$.

\item For each arrow $f \colon U \arr V$ a functor $f^* \colon \Phi V \arr \Phi U$.

\item For each object $U$ of $\cC$ an isomorphism $\isoiden_U\colon \id_U^* \simeq \id_{\Phi U}$ of functors $\Phi U \arr \Phi U$.

\item For each pair of arrows $U \xrightarrow{f} V \xrightarrow{g} W$ an isomorphism
   \[
   \isoass_{f,g} \colon f^*g^* \simeq (gf)^* \colon \Phi W \arr \Phi U
   \]
of functors $\Phi W \arr \Phi U$.

\end{enumeratei}

These data are required to satisfy the following conditions.

\begin{enumeratea}

\item If $f \colon U \arr V$ is an arrow in $\cC$ and $\eta$ is an object of $\Phi V$, we have
   \[
   \isoass_{\id_U, f}(\eta) = \isoiden_U (f^* \eta) \colon 
   \id_U^* f^* \eta \arr f^* \eta
   \]
and
   \[
    \isoass_{f, \id_V}(\eta) = f^*\isoiden_V (\eta) \colon 
   f^*\id_V^*  \eta \arr f^* \eta.
   \]

\item Whenever we have arrows $U \xrightarrow{f} V \xrightarrow{g} W \xrightarrow{h} T$ and an object $\theta$ of $\Phi(T)$, the diagram
   \[
   \xymatrix@C+30pt{
   f^*g^*h^* \theta \ar[r]^{\isoass_{f,g}(h^* \theta)}
   \ar[d]^{f^* \isoass_{g,h}(\theta)} &
   (gf)^*h^* \theta \ar[d]^{\isoass_{gf,h}(\theta)}\\
   f^*(hg)^* \theta \ar[r]^{\isoass_{f, hg}(\theta)}&
   (hgf)^* \theta
   }
   \]
commutes.
\end{enumeratea}
\end{definition}

In this definition we only consider (contravariant) pseudo-functors into the category of categories. Of course, there is a much more general notion of pseudo-functor with values in a 2-category, which we will not use at all.

A functor $\Phi$ from $\cS$ into the category of categories can be considered as a pseudo-functor, in which every $\isoiden_U$ is the identity on $\Phi U$, and every $\isoass_{f,g}$ is the identity on $f^{*}g^{*} = (gf)^{*}$.

We have seen how to associate with a fibered category over $\cC$, equipped with a cleavage, the data for a pseudo-functor; we still have to check that the two conditions of the definition are satisfied.

\begin{proposition}
A fibered category over $\cC$ with a cleavage defines a pseudo-functor on $\cC$.
\end{proposition}

\begin{proof}
We have to check that the two conditions are satisfied. Let us do this for condition~(b) (the argument for condition~(a) is very similar). The point is that $f^*g^*h^* \zeta$ and $(hgf)^* \zeta$ are both pullbacks of $\zeta$, and so, by the definition of cartesian arrow, there is a unique arrow $f^*g^*h^* \zeta \arr (hgf)^* \zeta$ lying over the identity on $U$, and making the diagram
   \[
   \xymatrix@C-15pt{
   f^*g^*h^* \zeta \ar[rr]\ar[rd] && (hgf)^* \zeta\ar[ld]\\
   &{}\zeta
   }
   \]
commutative. But one sees immediately that both $\isoass_{gf,h}(\zeta) \circ \isoass_{f,g}(h^* \zeta)$ and $\isoass_{f, hg}(\zeta) \circ f^* \isoass_{g,h}(\zeta)$ satisfy this condition.
\end{proof}


It is easy to see when a cleavage defines a functor from $\cC$ into the category of categories.

\begin{definition}
A cleavage on a fibered category is a \emph{splitting}\index{splitting} if it contains all the identities, and it is closed under composition.

A fibered category endowed with a splitting is called \emph{split}\index{category!fibered!split}\index{split fibered category}.
\end{definition}

\begin{proposition}
The pseudo-functor associated with a cleavage is a functor if and only if the cleavage is a splitting.
\end{proposition}

The proof is immediate.

In general a fibered category does not admit a splitting.

\begin{example}\label{ex:fibered-nosplitting}
Every group $G$ can be considered as a category with one object, where the set of arrows is exactly $G$, and the composition is given by the operation in $G$. A group homomorphism $G \arr H$ can be considered as a functor. An arrow in $G$ (that is, an element of $G$) is always cartesian; hence $G$ is fibered over $H$ if and only if $G \arr H$ is surjective.

Given a surjective homomorphism $G \arr H$, a cleavage $K$ is a subset of $G$ that maps bijectively onto $H$; and a cleavage is a splitting if and only if $K$ is a subgroup of $G$. So, a splitting is a splitting $H \arr  G$ of the homomorphism $G \arr H$, in the usual sense of a group homomorphism such that the composite $H \arr G \arr H$ is the identity on $H$. But of course such a splitting does not always exist.
\end{example}

Despite this, every fibered category is equivalent to a split fibered category (Theorem~\ref{thm:equivalent-split}).

\subsection{The fibered category associated with a pseudo-functor} \label{subsec:presheaf-cat}

So a fibered category with a cleavage defines a pseudo-functor. Conversely, from a pseudo-functor on $\cC$ one gets a fibered category over $\cC$ with a cleavage. First of all, let us analyze the case that the pseudo-functor is simply a functor $\Phi \colon \cC\op \arr \cat{Cat}$  into the category of categories, considered as a 1-category. This means that with each object $U$ of $\cC$ we associate a category $\Phi U$, and with each arrow $f \colon U \arr V$ a functor $\Phi f \colon \Phi V \arr \Phi U$, in such a way that $\Phi\id_U \colon \Phi U \arr \Phi U$ is the identity, and $\Phi(g\circ f) = \Phi f \circ \Phi g$ every time we have two composable arrows $f$ and $g$ in $\cC$.

With this $\Phi$, we can associate a fibered category $\cF \arr \cC$, such that for any object $U$ in $\cC$ the fiber $\cF(U)$ is canonically equivalent to the category $\Phi U$. An object of $\cF$ is a pair $(\xi, U)$ where $U$ is an object of $\cC$ and $\xi$ is an object of $\Phi U$. An arrow $(a, f) \colon (\xi, U) \arr (\eta, V)$ in $\cF$ consists of an arrow $f \colon U \arr V$ in $\cC$, together with an arrow $a \colon \xi \arr \Phi f(\eta)$ in $\Phi U$.

The composition is defined as follows: if
   \[
   (a, f) \colon (\xi, U) \arr (\eta, V)\quad
   \text{and}
   \quad (b, g) \colon (\eta, V) \arr (\zeta, W)
   \]
are two arrows, then 
   \[
   (b, g) \circ (a, f) = (\Phi f(b) \circ a, g \circ f) \colon 
   (\xi, U) \arr (\zeta, W).
   \]

There is an obvious functor $\cF \arr \cC$ that sends an object $(\xi, U)$ into $U$ and an arrow $(a, f)$ into $f$; I claim that this functor makes $\cF$ into a fibered category over $\cC$. In fact, given an arrow $f \colon U \arr V$ in $\cC$ and an object $(\eta, V)$ in $\cF(V)$, then $\bigl(\Phi f (\eta), U\bigr)$ is an object of $\cF(U)$, and it is easy to check that the pair $\bigl(\id_{\Phi f (\eta)}, f\bigr)$ gives a cartesian arrow $\bigl(\Phi f (\eta), U\bigr) \arr (\eta, V)$.

The fiber of $\cF$ is canonically isomorphic to the category $\Phi U$: the isomorphism $\cF(U) \arr \Phi U$ is obtained at the level of objects by sending $(\xi, U)$ to $\xi$, and at the level of arrows by sending $(a, \id_U)$ to $a$. The collection of all the arrows of type $\bigl(\id_{\Phi f (\eta)}, f\bigr)$ gives a splitting.

The general case is similar, only much more confusing. Consider a pseudo-functor $\Phi$ on $\cC$. As before, we define the objects of $\cF$ to be pairs $(\xi, U)$ where $U$ is an object of $\cC$ and $\xi$ is an object of $\cF(U)$. Again, an arrow $(a, f) \colon (\xi, U) \arr (\eta, V)$ in $\cF$ consists of an arrow $f \colon U \arr V$ in $\cC$, together with an arrow $a \colon \xi \arr f^*(\eta)$ in $\Phi U$.

Given two arrows $(a, f) \colon (\xi, U) \arr (\eta, V)$ and $(b, g) \colon (\eta, V) \arr (\zeta, W)$, we define  the composite $(b,g) \circ (a,f)$ as the pair $(b \cdot a, gf)$, where $b \cdot a = \isoass_{f,g}(\zeta) \circ f^* b \circ a$ is the composite
   \[
   \xi \overset a \larr f^* \eta
   \overset{f^* b} \larr f^*g^* \zeta
   \xrightarrow{\isoass_{f,g}(\zeta)} (gf)^* \zeta
   \]
in $\Phi U$.

Let us check that composition is associative. Given three arrows
   \[
   \xymatrix{
   (\xi, U)\ar[r]^{(a, f)}
   & (\eta, V)\ar[r]^{(b,g)}
   &(\zeta, W)\ar[r]^{(c, h)}
   &(\theta, T)
   }
   \]
we have to show that
   \[
   (c,h) \circ \bigl((b,g) \circ (a,f)\bigr)
   \eqdef \bigl(c \cdot (b \cdot a), hgf\bigr)
   \]
equals
   \[
   \bigl((c,h) \circ (b,g)\bigr) \circ (a,f) \eqdef
   \bigl((c \cdot b) \cdot a, hgf\bigr).
   \]
By the definition of the composition,  we have
   \begin{align*}
   c \cdot (b \cdot a) &=
   \isoass_{gf,h} (\theta) \circ (gf)^* c \circ (b \cdot a)\\
   & =  \isoass_{gf,h} (\theta) \circ (gf)^* c \circ 
   \isoass_{f,g}(\zeta) \circ f^* b \circ a
   \end{align*}
while
   \begin{align*}
   (c \cdot b) \cdot a &=
   \isoass_{f,hg}(\theta) \circ f^*(c \cdot b) \circ a\\
   & =  \isoass_{f,hg} (\theta) \circ f^*\isoass_{g,h}(\theta) \circ
   f^*g^*c \circ f^*b \circ a;
   \end{align*}
hence it is enough to show that the diagram
   \[
   \xymatrix@C+30pt{
   f^*g^* \zeta \ar[r]^{f^*g^*c} \ar[d]^{\isoass_{f,g}(\zeta)}&
   f^*g^*h^*\theta \ar[r]^{f^*\isoass_{g,h}(\theta)}
   \ar[d]^{\isoass_{f,g}(h^{*}\theta)} &
   f^*(hg)^* \theta\ar[d]^{\isoass_{f,hg}(\theta)} \\
   (gf)^*\zeta \ar[r]^{(gf)^*c} &
   (gf)^*h^*\theta \ar[r]^{\isoass_{gf,h}(\theta)} &
   (hgf)^*\theta
   }
   \]
commutes. But the commutativity of the first square follows from the fact that $\isoass_{f,g}$ is a natural transformation of functors, while that of the second is condition~(b) in Definition~\ref{def:pseudo-functor}.

Given an object $(\xi, U)$ of $\cF$, we have the isomorphism $\isoiden_U(\xi) \colon \id_U^* \xi \arr \xi$; we define the identity $\id_{(\xi,U)} \colon (\xi,U) \arr (\xi, U)$ as $\id_{(\xi,U)} = \bigl(\isoiden_U(\xi) ^{-1}, \id_U\bigr)$. To check that this is neutral with respect to composition, take an arrow $(a, f) \colon (\xi,U) \arr (\eta, V)$; we have 
   \[
   (a,f) \circ \bigl(\isoiden_U(\xi)^{-1}, \id_U\bigr) =
   \bigl(a \cdot \isoiden_U(\xi)^{-1}, f\bigr)
   \]
and
   \[
   a \cdot \isoiden_U(\xi)^{-1} = \isoass_{\id_U,f}(\eta)
   \circ \id_U^* a \circ \isoiden_U(\xi)^{-1}.
   \]
But condition~(a) of Definition~\ref{def:pseudo-functor} says that $\isoass_{\id_U,f}(\eta)$ equals $ \isoiden_U (f^* \eta)$, while the diagram
   \[
   \xymatrix@C+=40pt{
   {}\id_U^*\xi \ar[r]^-{\isoiden_U(\xi)} \ar[d]^{\id_U^*a} &
   {}\xi \ar[d]^{a}\\
   {}\id_U^*f^*\xi \ar[r]^-{\isoiden_U(f^*\eta)} &
   {}f^*\eta
   }
   \]
commutes, because $\isoiden_U$ is a natural transformation. This implies that $a \cdot \isoiden_U(\xi)^{-1} = a$, and therefore $(a,f) \circ \bigl(\isoiden_U(\xi)^{-1}, \id_U\bigr) = (a,f)$.

A similar argument shows that $\bigl(\isoiden_U(\xi)^{-1}, \id_U\bigr)$ is also a left identity.

Hence $\cF$ is a category. There is an obvious functor $\p_{\cF} \colon \cF \arr \cC$ sending an object $(\xi, U)$ to $U$ and an arrow $(a,f)$ to $f$. I claim that this makes $\cF$ into a category fibered over $\cC$.

Take an arrow $f \colon U \arr V$ of $\cC$, and an object $(\eta, V)$ of $\cF$ over $V$. I claim that the arrow
   \[
   (\id_{f^*\eta}, f) \colon (f^*\eta, U) \arr
   (\eta, V)
   \]
is cartesian. To prove this, suppose that we are given a diagram
   \[
   \xymatrix@R-10pt@C+10pt{
   (\zeta,W) \ar@{|->}[dd] \ar@{-->}[rd]_{(c,h)}
   \ar@/^/[rrd]^{(b,g)} \\
   & (f^*\eta, U) \ar@{|->}[dd]\ar[r]_{(\id_{f^*\eta},f)}
   & (\eta,V) \ar@{|->}[dd] \\
   W\ar[rd]_h\ar@/^/[rrd]^<<<<<<<<g|(.4965)\hole\\
   & U \ar[r]_f
   & V
   }
   \]
(without the dotted arrow); we need to show that there is a unique arrow $(c,h)$ that can be inserted in the diagram. But it is easy to show that
   \[
   (\id_{f^*\eta},f) \circ (c,h) =
   (\isoass_{h,f}(\eta) \circ c, fh),
   \]
and this tells us that the one and only arrow that fits into the diagram is $(\isoass_{h,f}(\eta)^{-1} \circ b, h)$.

This shows that $\cF$ is fibered over $\cC$, and also gives us a cleavage.

Finally, let us notice that for all objects $U$ of $\cC$ there is functor $\cF(U) \arr \Phi U$, sending an object $(\xi, U)$ to $\xi$ and an arrow $(a, f)$ into $a$. This is an isomorphism of categories.

The cleavage constructed above gives, for each arrow $f \colon U \arr V$, functors $f^* \colon \cF(V) \arr \cF(U)$. If we identify each $\cF(U)$ with $\Phi U$ via the isomorphism above, then these functors correspond to the $f^* \colon \Phi V \arr \Phi U$. Hence if we start with a pseudo-functor, we construct the associated fibered category with a cleavage, and then we take the associated pseudo-functor, this is isomorphic to the original pseudo-functor (in the obvious sense).

Conversely, it is easy to see that if we start from a fibered category with a cleavage, construct the associated pseudo-functor, and then take the associated fibered category with a cleavage, we get something isomorphic to the original fibered category with a cleavage (again in the obvious sense). So really giving a pseudo-functor is the same as giving a fibered category with a cleavage.

On the other hand, since cartesian pullbacks are unique up to a unique isomorphism (Remark~\ref{rmk:unique-pullback}), also cleavages are unique up to a unique isomorphism. This means that, in a sense that one could make precise, the theory of fibered categories is equivalent to the theory of pseudo-functors. On the other hand, as was already remarked in \cite[Remarque, pp.~193--194]{sga1}, often the choice of a cleavage hinders more than it helps.

\section{Examples of fibered categories}\label{sec:examples-fibered}

\begin{example}\label{ex:arrows} Assume that $\cC$ has fibered products. Let $\Arr \cC$ be the category of arrows in $\cC$; its objects are the arrows in $\cC$, while an arrow from $f\colon X \arr U$ to $g \colon Y \arr V$ is a commutative diagram
   \[
   \xymatrix{
   X\ar[r] \ar[d]^f & Y \ar[d]^g\\
   U\ar[r] & V
   }
   \]
The functor $\p_{\Arr \cC}\colon \Arr \cC\arr \cC$ sends each arrow $S \arr U$ to its codomain $U$, and each commutative diagram to its bottom row.

I claim that $\Arr \cC$ is a fibered category over $\cC$. In fact, it easy to check that the cartesian diagrams are precisely the cartesian squares, so the statement follows from the fact that $\cC$ has fibered products.

\end{example}

\begin{definition}\label{def:stable-class}
A class $\cP$ of arrows in a category $\cC$ is \emph{stable}\index{class of arrows!stable}\index{stable class of arrows} if the following two conditions hold.

\begin{enumeratea}

\item If $f \colon X \arr U$ is in $\cP$, and $\phi \colon X' \simeq X$, $\psi \colon U \simeq U'$ are isomorphisms, the composite
   \[
  \psi \circ f \circ \phi \colon  X' \arr U'
   \]
is in $\cP$.

\item Given an arrow $Y \arr V$ in $\cP$ and any other arrow $U \arr V$, then a fibered product $U \times_V Y$ exists, and the projection $U \times_V Y \arr U$ is in $\cP$.

\end{enumeratea} 
\end{definition}

\begin{example}\label{ex:restrictedarrows} As a variant of the
example above, let $\cP$ be a stable class of arrows. The arrows in $\cP$ are the objects in a category, again denoted by $\cP$, in which an arrow from $X \arr U$ to $Y \arr V$ is a commutative diagram
   \[
   \xymatrix{
   X \ar[r]\ar[d] & Y \ar[d] \\
   U \ar[r]       & V        \\
   }
   \]

It is easy to see that this is a fibered category over $\cC$; the cartesian arrows are precisely the cartesian diagrams.

\end{example}

\begin{example}\label{ex:classifying} Let $G$ a topological group. The \emph{classifying stack}\index{classifying stack of a topological group} of $G$ is the fibered category $\cB G \arr \cattop$ over the category of topological spaces, whose objects are principal $G$ bundles $P \arr U$, and whose arrows $(\phi, f)$ from $P \arr U$ to $Q \arr V$  are commutative diagrams
   \[
   \xymatrix{
   P\ar[r]^\phi\ar[d] & Q\ar[d]\\
   U\ar[r]^f & V
   }
   \]
where the function $\phi$ is $G$-equivariant. The functor $\cB G \arr \cattop$ sends a principal bundle $P \arr U$ into the topological space $U$,  and an arrow $(\phi, f)$ into $f$.

This fibered category $\cB G \arr \cattop$ has the property that each of its arrows is cartesian.
\end{example}

\begin{example}
Here is an interesting example, suggested by one of the participants in the school. Consider the forgetful functor $F \colon \cattop \arr \catset$ that associates with each topological space $X$ its underlying set $FX$, and to each continuous function the function itself.

I claim that this makes $\cattop$ fibered over $\catset$. Suppose that you have a topological space $Y$, a set $U$ and a function $f \colon U \arr FY$. Denote by $X$ the set $U$ with the \emph{initial topology}, in which the open sets are the inverse images of the open subsets of $Y$; this is the coarsest topology that makes $f$ continuous. If $T$ is a topological space, a function $T \arr X$ is continuous if and only if the composite $T \arr X \arr Y$ is continuous; this means that $f \colon X \arr Y$ is a cartesian arrow over the given arrow $f \colon U \arr FY$.

The fiber of $\cattop$ over a set $U$ is the partially ordered set of topologies on $U$, made into a category in the usual way.

Notice that in this example the category $\cattop$ has a canonical splitting over $\catset$.
\end{example}

We are interested in categories of sheaves. The simplest example is the fibered category of sheaves on objects of a site, defined as follows.

\begin{example}\label{ex:fibered-sheaves}
\index{fibered!category!of sheaves}
\index{category!fibered!of sheaves}
Let $\cC$ be a site, $\cT$ its topology. We will refer to a sheaf in the category $(\cC/X)$, endowed with the comma topology (Definition~\ref{def:comma-topology}) as a \emph{sheaf on $X$}, and denote the category of sheaves on $X$ by $\sh{X}$.

If $f \colon X \arr Y$ is an arrow in $\cC$, there is a corresponding restriction functor $f^* \colon \sh{Y} \arr \sh{X}$, defined as follows.

If $F$ is a sheaf on $Y$ and $U \arr X$ is an object of $(\cC/X)$, we define $f^*F(U \arr X) \eqdef F(U \arr Y)$, where $U \arr Y$ is the composite of $U \arr X$ with $f$.

If $U \arr X$ and $V \arr X$ are objects of $(\cC/X)$ and $\phi \colon U \arr V$ is an arrow in $(\cC/X)$, then $\phi$ is also an arrow from $U \arr Y$ to $V \arr Y$, hence it induces a function $\phi^* \colon f^*F(U \arr X) = F(U \arr Y) \arr F(V \arr Y) = f^*F(V \arr X)$. This gives $f^*F$ the structure of a functor $(\cC/X)\op \arr \catset$. One sees easily that $f^*F$ is a sheaf on $X$.

If $\phi \colon F \arr G$ is a natural transformation of sheaves on $(\cC/Y)$, there is an induced natural transformation $f^*\phi \colon f^*F \arr f^*G$ of sheaves on $(\cC/X)$, defined in the obvious way. This defines a functor $f^* \colon \sh Y \arr \sh X$.

It is immediate to check that, if $f \colon X \arr Y$ and $g \colon Y \arr Z$ are arrows in $\cC$, we have an \emph{equality} of functors $(gf)^* = f^*g^* \colon (\cC/Z) \arr (\cC/X)$. Furthermore $\id_X^* \colon \sh X \arr \sh X$ is the identity. This means that we have defined a functor from $\cC$ to the category of categories, sending an object $X$ into the category of sheaves on $(\cC/X)$. According to the result of \S\ref{subsec:presheaf-cat}, this yields a category $\catsh{\cC} \arr \cC$, whose fiber over $X$ is $\sh{X}$.
\end{example}

There are many variants on this example, by considering sheaves in abelian groups, rings, and so on.

This example is particularly simple, because it is defined by a functor. In most of the cases that we are interested in, the sheaves on a given object will be defined in a site that is not the one inherited from the base category $\cC$; this creates some difficulties, and forces one to use the unpleasant machinery of pseudo-functors. On the other hand, this discrepancy between the topology on the base and the topology in which the sheaves are defined is what makes descent theory for \qc sheaves so much more than an exercise in formalism.

Let us consider directly the example we are interested in, that is, fibered categories of \qc sheaves.

\subsection{The fibered category of \qc sheaves}
\label{subsec:fibered-quasi-coherent}
\index{fibered!category!of quasi-coherent sheaves}
\index{category!fibered!of quasi-coherent sheaves}
Here $\cC$ will be the category $\catsch{S}$ of schemes over a fixed base scheme $S$. For each scheme $U$ we define $\qcoh U$ to be the category of \qc sheaves on $U$. Given a morphism $f \colon U \arr V$, we have a functor $f^* \colon \qcoh V \arr \qcoh U$. Unfortunately, given two morphisms $U \xrightarrow{f} V \xrightarrow{g} W$, the pullback $(gf)^* \colon \qcoh W \arr \qcoh U$ does not coincide with the composite $f^*g^* \colon \qcoh W \arr \qcoh U$, but it is only  canonically isomorphic to it. This may induce one to suspect that we are in the presence of a pseudo-functor; and this is indeed the case.

The neatest way to prove this is probably by exploiting the fact that the pushforward $f_* \colon \qcoh U \arr \qcoh V$ is functorial, that is, $(gf)_*$ equals $g_* f_*$ on the nose, and $f^*$ is a left adjoint to $f_*$. This means that, given \qc sheaves $\cM$ on $U$ and $\cN$ on $V$, there is a canonical isomorphism of groups
   \[
   \Theta_f(\cN,\cM)\colon
   \hom_{\cO_V}(\cN, f_*\cM) \simeq
   \hom_{\cO_U}(f^* \cN, \cM)
   \]
that is natural in $\cM$ and $\cN$. More explicitly, there are two functors
   \[
   \qcoh U\op \times \qcoh V \arr \catgrp
   \]
defined by
   \[
   (\cM,\cN) \mapsto
   \hom_{\cO_V}(\cN, f_*\cM)
   \]
and
   \[
   (\cM,\cN) \mapsto
   \hom_{\cO_U}(f^* \cN, \cM);
   \]
then $\Theta_f$ defines a natural isomorphism from the first to
the second.

Equivalently, if $\alpha \colon \cM \arr \cM'$ and $\beta \colon \cN \arr \cN'$ are homomorphisms of \qc sheaves on $U$ and $V$ respectively, the diagrams
   \[
   \xymatrix@C+20pt{
   \hom_{\cO_V}(\cN, f_*\cM)
   \ar[d]^{f_* \alpha \circ -}
   \ar[r]^{\Theta_f(\cN,\cM)}  &
   \hom_{\cO_U}(f^* \cN, \cM)
   \ar[d]^{\alpha \circ -} \\
   \hom_{\cO_V}(\cN, f_*\cM')
   \ar[r]^{\Theta_f(\cN,\cM')} &
   \hom_{\cO_U}(f^* \cN, \cM')
   }
   \]
and
   \[
   \xymatrix@C+20pt{
   \hom_{\cO_V}(\cN', f_*\cM)
   \ar[r]^{\Theta_f(\cN',\cM)} \ar[d]^{- \circ \beta} &
   \hom_{\cO_U}(f^* \cN', \cM)
   \ar[d]^{- \circ f^*\beta}  \\
   \hom_{\cO_V}(\cN, f_*\cM)
   \ar[r]^{\Theta_f(\cN,\cM)} &
   \hom_{\cO_U}(f^{*}\cN, \cM)
   }
   \]
commute.

If $U$ is a scheme over $S$ and $\cN$  a \qc sheaf on $U$, then the pushforward functor $(\id_U)_* \colon \qcoh U \arr \qcoh U$ is the identity (this has to be interpreted literally, I am not simply asserting the existence of a canonical isomorphism between $(\id_U)_*$ and the identity on $\qcoh U$). Now, if $\cM$ is a \qc sheaf on $U$, there is a canonical adjunction isomorphism
   \[
   \Theta_{\id_U}(\cM, -) \colon
   \hom_{\cO_U}(\cM, (\id_U)_*-) = 
   \hom_{\cO_U}(\cM, -) \simeq
   \hom_{\cO_U}(\id_U^*\cM, -)
   \]
of functors from $\qcoh U$ to $\catset$. By the dual version of Yoneda's lemma (Remark~\ref{rmk:dual-Yoneda}) this corresponds to an isomorphism $\isoiden_U (\cM) \colon \id_U^*\cM \simeq \cM$. This is easily seen to be functorial, and therefore defines an isomorphism
   \[
   \isoiden_U \colon \id_U^* \simeq \id_{\qcoh U}   
   \]
of functors from $\qcoh U$ to itself. This isomorphism is the usual one: a section $s \in \cM(A)$, for some open subset $A \subseteq U$, yields a section $\id_U^* s \in \cM(\id_U ^{-1} A) = \cM(A)$, and $\isoiden_U(\cM)$ sends $\id_U^* s$ to $s$. This is the first piece of data that we need.

For the second, consider two morphisms $U \xrightarrow{f} V \xrightarrow{g} W$ and a \qc sheaf $\cP$ on $W$. We have the chain of isomorphisms of functors $\qcoh U \arr \catgrp$
   \begin{alignat*}{3}
   \hom_{\cO_U}\bigl((gf)^*\cP, -\bigr)
   &{}\simeq \hom_{\cO_W}\bigl(\cP, (gf)_*-\bigr)
   && \quad\bigl(\text{this is $\Theta_{gf}(\cP,-)^{-1}$}\bigr)\\
   &{}=\hom_{\cO_W}\bigl(\cP, g_*f_*-\bigr)\\
   &{} \simeq \hom_{\cO_V}(g^*\cP, f_*-)
   && \quad\bigl(\text{this is $\Theta_{g}(\cP,f_*-)$}\bigr)\\
   &{} \simeq \hom_{\cO_U}(g^*f^*\cP, -)
   && \quad\bigl(\text{this is $\Theta_{f}(g^*\cP,-)$}\bigr);
   \end{alignat*}
the composite
   \begin{align*}
   \Theta_{f}(g^*\cP,-) \circ \Theta_{g}(\cP, f_*-)
   \circ {}&\Theta_{gf}(\cP,-)^{-1} \colon\\
   &\hom_{\cO_U}\bigl((gf)^*\cP, -\bigr)
   \simeq \hom_{\cO_U}(g^*f^*\cP, -)
   \end{align*}
corresponds, again because of the covariant Yoneda lemma, to  an isomorphism $\isoass_{f,g}(\cP) \colon f^*g^*\cP\simeq (gf)^*\cP$. These give an isomorphism $\isoass_{f,g} \colon f^*g^* \simeq (gf)^*$ of functors $\qcoh W \arr \qcoh U$. Once again, this is the usual isomorphism: given a section $s \in \cP(A)$ for some open subset $A \subseteq W$, there are two sections
   \[
   (gf)^*s \in (gf)^*\cP\bigl((gf)^{-1} A\bigr) =
   (gf)^*\cP(f ^{-1} g ^{-1} A)
   \]
and
   \[
   f^* g^* s \in f^* g^* \cP(f ^{-1} g ^{-1} A);
   \]
the isomorphism $\isoass_{f,g}(\cP)$ sends $f^* g^* s$ into $(gf)^*s$. Since the sections of type $(gf)^*s$ generate $(gf)^*\cP$ as a sheaf of $\cO_U$-modules, this characterizes $\isoass_{f,g}(\cP)$ uniquely.

We have to check that the $\isoiden_U$ and $\isoass_{f,g}$ satisfy the conditions of Definition~\ref{def:pseudo-functor}. This can be done directly at the level of sections, or using the definition of the two isomorphisms via the covariant Yoneda lemma; we will follow the second route. Take a morphism of schemes $f \colon U \arr V$. We need to prove that for any \qc sheaf $\cN$ on $V$ we have the equality
   \[
   \isoass_{\id_U,f}(\cN) = \isoiden_U(f^*\cN) \colon 
   \id_U^*f^* \cN \arr f^*\cN.
   \]
This is straightforward: by the covariant Yoneda lemma, it is enough to show that $\isoass_{id_U,f}(\cN)$ and $\isoiden_U(f^*\cN)$ induce the same natural transformation
   \[
   \hom_{\cO_U}(f^* \cN, - ) \arr
   \hom_{\cO_U}(\id_U^*f^* \cN, - ).
   \]
But by definition the natural transformation induced by $\isoiden_U(f^*\cN)$ is
   \[
   \Theta_{\id_U}(f^* \cN, -),
   \]
while that induced by $ \isoass_{\id_U,f}(\cN) $ is 
   \[
   \Theta_{\id_U}(f^* \cN, -) \circ
   \Theta_f(\cN,(\id_U)_*-) \circ
   \Theta_f (\cN, -)^{-1} = 
   \Theta_{\id_U}(f^* \cN, -).
   \]
Similar arguments works for the second part of the first condition and for the second condition.

The fibered category on $\catsch{S}$ associated with this pseudo-functor is the fibered category of \qc sheaves, and will be denoted by $\catqc{S} \arr \catsch{S}$.

There are many variants on  this example. For example, one can define the fibered category of sheaves of $\cO$-modules over the category of ringed topological spaces in exactly the same way.

\section{Categories fibered in groupoids}

\begin{definition} A \emph{category fibered in groupoids over $\cC$}\index{category!fibered!in groupoids}\index{fibered!category!in groupoids} is a category $\cF$ fibered over $\cC$, such that the category $\cF(U)$ is a groupoid for any object $U$ of $\cC$.
\end{definition}

In the literature one often finds a different definition of a category fibered in groupoids.

\begin{proposition}\call{prop:char-fibered-in-groupoids} Let $\cF$ be a category over $\cC$. Then $\cF$ is fibered in groupoids over $\cC$ if and only if the following two conditions hold.

\begin{enumeratei}
\itemref{1} Every arrow in $\cF$ is cartesian.

\itemref{2} Given an object $\eta$ of $\cF$ and an arrow $f \colon U \arr \p_\cF \eta$ of $\cC$, there exists an arrow $\phi \colon \xi \arr \eta$ of $\cF$ with $\p_{\cF} \phi = f$.
\end{enumeratei}
\end{proposition}

\begin{proof} Suppose that these two conditions hold: then clearly $\cF$ is fibered over $\cC$. Also, if $\phi\colon \xi \arr \eta $ is an arrow of $\cF(U)$ for some object $U$ of $\cC$, then we see from condition~\refall{prop:char-fibered-in-groupoids}{1} that there exists an arrow $\psi\colon \eta \arr \xi$ with $\p_{\cF} \psi = \id_U$ and $\phi\psi = \id_\eta$; that is, every arrow in $\cF(U)$ has a right inverse. But this right inverse $\psi$ also must also have a right inverse, and then the right inverse of $\psi$ must be $\phi$. This proves that every arrow in $\cF(U)$ is invertible.

Conversely, assume that $\cF$ is fibered over $\cC$, and each $\cF(U)$ is a groupoid. Condition~\refpart{prop:char-fibered-in-groupoids}{2} is trivially verified. To check condition~\refpart{prop:char-fibered-in-groupoids}{1}, let $\phi\colon \xi \arr \eta$ be an arrow in $\cC$ mapping to $f \colon U \arr V$ in $\cC$. Choose a pullback $\phi'\colon \xi' \arr \eta$ of $\eta$ to $U$; by definition there will be an arrow $\alpha \colon \xi \arr \xi'$ in $\cF(U)$ such that $\phi'\alpha = \phi$. Since $\cF(U)$ is a a groupoid, $\alpha$ will be an isomorphism, and this implies that $\phi$ is cartesian.
\end{proof}

\begin{corollary} Any base-preserving functor from a fibered category to a category fibered in groupoids is a morphism.
\end{corollary}

\begin{proof} This is clear, since every arrow in a category fibered in groupoids is cartesian.
\end{proof}

Of the examples of Section~\ref{sec:examples-fibered}, \ref{ex:arrows} and \ref{ex:restrictedarrows} are not in general fibered in groupoids, while the classifying stack of a topological group introduced in \ref{ex:classifying} is always fibered in groupoids.

\section{Functors and categories fibered in sets}\label{sec:fiberedsets}

The notion of category generalizes the notion of set: a set can be thought of as a category in which every arrow is an identity. Furthermore functors between sets are simply functions.

Similarly, fibered categories are generalizations of functors.

\begin{definition} A \emph{category fibered in sets over $\cC$}\index{category!fibered!in sets}\index{fibered!category!in sets} is a category $\cF$ fibered over $\cC$, such that for any object $U$ of $\cC$ the category $\cF(U)$ is a set.
\end{definition}

Here is a useful characterization of categories fibered in sets.

\begin{proposition}\label{prop:characterization-fibered-sets}
\index{category!fibered!in sets, characterization}
\index{fibered!category!in sets, characterization}
Let $\cF$ be a category over $\cC$. Then $\cF$ is fibered in sets if and only if for any object $\eta$ of $\cF$ and any arrow $f \colon U \arr  \p_{\cF} \eta$ of $\cC$, there is a unique arrow $\phi\colon \xi \arr \eta$ of $\cF$ with $\p_{\cF} \phi = f$.
\end{proposition}

\begin{proof} Suppose that $\cF$ is fibered in sets. Given $\eta$ and $f \colon U \arr \p_{\cF} \eta$ as above, pick a cartesian arrow $\xi \arr \eta$ over $f$. If $\xi' \arr \eta$ is any other arrow over $f$, by definition there exists an arrow $\xi' \arr \xi$ in $\cF(U)$ making the diagram 
   \[
   \xymatrix@-15pt{{}\xi' \ar@{-->}[rr]\ar[rd] && {}\xi\ar[ld]\\
             & {}\eta \\}
   \]
commutative. Since $\cF(U)$ is a set, it follows that this arrow $\xi' \arr \xi$ is the identity, so the two arrows $\xi\arr \eta$ and $\xi' \arr \eta$ coincide.

Conversely, assume that the condition holds. Given a diagram
   \[
   \xymatrix@R=6pt{{\zeta}\ar@{|->}[dd] \ar@{-->}[rd]_\theta
   \ar@/^/[rrd]^{\psi} \\
   & {}\xi\ar@{|->}[dd]\ar[r]_{\phi}
   &{}\eta\ar@{|->}[dd] \\
   {}\p_{\cF} \zeta\ar[rd]_h\ar@/^/[rrd]|(.487)\hole\\
   & {}\p_{\cF} \xi \ar[r]
   & {}\p_{\cF} \eta
   }
   \]
the condition implies that the only arrow $\theta\colon \zeta \arr \xi$ over $h$ makes the diagram commutative; so the category $\cF$ is fibered.

It is obvious that the condition implies that $\cF(U)$ is a set for all $U$.
\end{proof}

So, for categories fibered in sets the pullback of an object of $\cF$ along an arrow of $\cC$ is strictly unique. It follows from this that when $\cF$ is fibered in sets over $\cC$ and $f \colon U \arr V$ is an arrow in $\cC$, the pullback map $f^* \colon \cF(V) \arr \cF(U)$ is uniquely defined, and the composition rule $f^*g^* = (gf)^*$ holds.  Also for any object $U$ of $\cC$ we have that $\id_U^* \colon \cF(U) \arr \cF(U)$ is the identity. This means that we have defined a functor $\Phi_\cF \colon \cC\op \arr \catset$ by sending each object $U$ of $\cC$ to $\cF(U)$, and each arrow $f \colon U \arr V$ of $\cC$ to the function $f^* \colon \cF(V) \arr \cF(U)$.

Furthermore, if $F \colon \cF \arr \cG$ is a morphism of categories fibered in sets, because of the condition that $\p_{\cG} \circ F = \p_{\cF}$, then every arrow in $\cF(U)$, for some object $U$ of $\cC$, will be sent to $\cG(U)$. So we get a function $F_U \colon \cF(U) \arr \cG(U)$. It is immediate to check that this gives a natural transformation $\phi_F \colon \Phi_\cF \arr \Phi_\cG$.

There is a category of categories fibered in sets over $\cC$, where the arrows are morphisms of fibered categories; the construction above gives a functor from this category to the category of functors $\cC\op \arr \catset$.

\begin{proposition}
\index{category!fibered!in sets, equivalence with functors}
\index{fibered!category!in sets, equivalence with functors}

This is an equivalence of the category of categories fibered in sets over $\cC$ and the category of functors $\cC\op \arr \catset$.
\end{proposition}

\begin{proof} An inverse functor is obtained by the construction of \S\ref{subsec:presheaf-cat}. Consider a functor $\Phi \colon \cC\op \arr \catset $: we construct a category fibered in sets $\cF_\Phi$ as follows. The objects of $\cF_\Phi$ will be pairs $(U, \xi)$, where $U$ is an object of $\cC$, and $\xi\in \Phi U$. An arrow from $(U, \xi)$ to $(V, \eta)$ is an an arrow $f \colon U \arr V$ of $\cC$ with the property that $\Phi f \eta = \xi$. It follows from Proposition~\ref{prop:characterization-fibered-sets} that $\cF_\Phi$ is fibered in sets over $\cC$.

With each natural transformation of functors $\phi\colon \Phi \arr \Phi'$ we associate a morphism $F_\phi \colon \cF_\Phi \arr \cF_{\Phi'}$. An object $(U, \xi)$ of $\cF_\Phi$ will be sent to $(U, \phi_U \xi)$. An arrow $f \colon (U, \xi) \arr (V, \eta)$ in $\cF_\Phi$ is is simply an arrow $f \colon U \arr V$ in $\cC$, with the property that $\Phi f(\eta) = \xi $. This implies that $\Phi'(f) \phi_V(\eta) = \phi_U \Phi(f) (\eta) = \phi_V \xi$, so the same $f$ will yield an arrow $f \colon (U, \phi_U \xi) \arr (V, \phi_V \eta)$. 

We leave it the reader to check that this defines a functor from the category of functors to the category of categories fibered in sets.
\end{proof}

So, any functor $\cC\op  \arr \catset $ will give an example of a fibered category over $\cC$.

\begin{remark}
It is interesting to notice that if $F \colon \cC\op \arr \catset$ is a functor and $\cF \arr \cC$ the associated category fibered in sets, then an object $(X, \xi)$ of $\cF$ is universal pair for the functor $F$ if and only if it is a terminal object for $\cF$. Hence $F$ is representable if and only if $\cF$ has a terminal object.
\end{remark}

In particular, given an object $X$ of $\cC$, we have the representable functor 
   \[
   \h_X\colon \cC\op \arr \catset,
   \]
defined on objects by the rule $\h_X U = \hom_{\cC}(U,X)$. The category fibered in sets over $\cC$ associated with this functor is the comma category $(\cC/X)$, and the functor $(\cC/X) \arr \cC$ is the functor that forgets the arrow into $X$.
 
So the situation is the following. From Yoneda's lemma we see that the category $\cC$ is embedded into the category of functors $\cC\op \arr \catset$, while the category of functors is embedded into the category of fibered categories.

From now we will identify a functor $F \colon \cC\op \arr \catset$ with the corresponding category fibered in sets over $\cC$, and will (inconsistently) call a category fibered in sets simply ``a functor''.

\subsection{Categories fibered over an object}
\label{subsec:fibered-over}

\begin{proposition}\label{prop:fibered-over}

Let $\cG$ be a category fibered in sets over $\cC$, $\cF$ another category, $F \colon \cF \arr \cG$ a functor. Then $\cF$ is fibered over $\cG$ if and only if it is fibered over $\cC$ via the composite $\p_{\cG} \circ F \colon \cF \arr \cC$.

Furthermore, $\cF$ is fibered in groupoids over $\cG$ if and only if it fibered in groupoids over $\cC$, and is fibered in sets over $\cG$ if and only if it fibered in sets over $\cC$.
\end{proposition}

\begin{proof}
One sees immediately that an arrow of $\cF$ is cartesian over its image in $\cG$ if and only if it is cartesian over its image in $\cC$, and the first statement follows from this.

Furthermore, one sees that the fiber of $\cF$ over an object $U$ of $\cC$ is the disjoint union, as a category, of the fibers of $\cF$ over all the objects of $\cG$ over $U$; these fiber are groupoids, or sets, if and only if their disjoint union is.
\end{proof}

This can be used as follows. Suppose that $S$ is an object of $\cC$, and consider the category fibered in sets $(\cC/S) \arr \cC$, corresponding to the representable functor $\h_S \colon \cC\op \arr \catset$. By Proposition~\ref{prop:fibered-over}, a fibered category $\cF \arr (\cC/S)$ is the same as a fibered category $\cF \arr \cC$, together with a morphism $\cF\arr (\cC/S)$ of categories fibered over $\cC$.

In particular, categories fibered in sets correspond to functors; hence we get that giving a functor $(\cC/S)\op \arr \catset$ is equivalent to assigning a functor $\cC\op \arr \catset$ together with a natural transformation $F \arr \h_{S}$. Describing this process for functors seems less natural than for fibered categories in general.

Given a functor $F \colon (\cC/S)\op \arr \catset$, this corresponds to a category fibered in sets $F \arr (\cC/S)$, which can be composed with the forgetful functor $(\cC/S) \arr \cC$ to get a category fibered in sets $F \arr \cC$, which in turn corresponds to a functor $F' \colon \cC\op \arr \catset$. What is this functor? One minute's thought will convince you that it can be described as follows: $F'(U)$ is the disjoint union of the $F(U \to S)$ for all the arrows $U \to S$ in $\cC$. The action of $F'$ on arrows is the obvious one.

\subsection{Fibered subcategories}

\begin{definition}
Let $\cF\arr \cC$ be a fibered category. A \emph{fibered subcategory}\index{subcategory!fibered}\index{fibered!subcategory} $\cG$ of $\cF$ is a subcategory of $\cF$, such that the composite $\cG \into \cF \arr \cC$ makes $\cG$ into a fibered category over $\cC$, and such that any cartesian arrow in $\cG$ is also cartesian in $\cF$.
\end{definition}

The last condition is equivalent to requiring that the inclusion
$\cG \into \cF$ is a morphism of fibered categories.

\begin{example}
Let $\cF \arr \cC$ be a fibered category, $\cG$ a full subcategory of $\cF$, with the property that if $\eta$ is an object of $\cG$ and $\xi \arr \eta$ is a cartesian arrow in $\cF$, then $\xi$ is also is $\cG$. Then $\cG$ is a fibered subcategory of $\cF$; the cartesian arrows in $\cG$ are the cartesian arrows in $\cF$ whose target is in $\cG$.
\end{example}

So, for example, the category of locally free sheaves is a fibered subcategory of the fibered category $\catqc{S}$ over $\catsch{S}$. 

Here is an interesting example.

\begin{definition}\label{def:ass-groupoids}
Let $\cF \arr \cC$ be a fibered category. The \emph{category fibered in groupoids associated with $\cF$}\index{category!fibered!$\cF\cart$}\index{fibered!category!$\cF\cart$}\index{$\cF\cart$} is the subcategory $\cF\cart$ of $\cF$, whose objects are all the objects of $\cF$, and whose arrows are the cartesian arrows of $\cF$.
\end{definition}

\begin{proposition}
If $\cF \arr \cC$ is a fibered category, then $\cF\cart \arr \cC$ is fibered in groupoids.

Furthermore, if $F\colon \cG \arr \cF$ is a morphism of fibered categories and $\cG$ is fibered in groupoids, then the image of $F$ is in $\cF\cart$.
\end{proposition}

The proof is immediate from Proposition~\ref{prop:char-fibered-in-groupoids} and from Proposition \refall{prop:comp-cartesian}{4}.

\section{Equivalences of fibered categories}

\subsection{Natural transformations of functors}

The fact that fibered categories are categories, and not functors, has strong implications, and does cause difficulties. As usual, the main problem is that functors between categories can be isomorphic without being equal; in other words, functors between two fixed categories form a category, the arrows being given by natural transformations.

\begin{definition}\label{def:morfunctors} Let $\cF$ and $\cG$ be two categories fibered over $\cC$, $F$, $G \colon \cF \arr \cG$ two morphisms. A \emph{base-preserving natural transformation}\index{natural transformation!base-preserving} \index{base-preserving!natural transformations} $\alpha\colon F \arr G$ is a natural transformation such that for any object $\xi$ of $\cF$, the arrow $\alpha_\xi \colon F \xi \arr G \xi$ is in $\cG(U)$, where $U \eqdef \p_{\cF} \xi = \p_{\cG}(F \xi) = \p_{\cG}(G \xi)$.

An \emph{isomorphism of $F$ with $G$} is a base-preserving natural transformation $F \arr G$ which is an isomorphism of functors. 
\end{definition}

It is immediate to check that the inverse of a base-preserving isomorphism is also base-preserving.

There is a category whose objects are the morphism from $\cF$ to $\cG$, and the arrows are base-preserving natural transformations; we denote it by $\hom_\cC(\cF, \cG)$.

\subsection{Equivalences}

\begin{definition} Let $\cF$ and $\cG$ be two fibered categories over $\cC$. An \emph{equivalence}\index{equivalence!of fibered categories}, of $\cF$ with $\cG$ is a morphism $F \colon \cF \arr \cG$, such that there exists another morphism $G \colon \cG \arr \cF$, together with isomorphisms of $G \circ F$ with $\id_\cF$ and of $F \circ G$ with $\id_\cG$.

We call $G$ simply an \emph{inverse} to $F$.
\end{definition}

\begin{proposition}\label{prop:equiv->equivalenthom}
\index{category!fibered!equivalence of, characterization}
\index{fibered!category!equivalence of, characterization}
\index{equivalence!of fibered categories!characterization}
Suppose that $\cF$, $\cF'$, $\cG$ and $\cG'$ are categories fibered over $\cC$. Suppose  that $F \colon \cF' \arr \cF$ and $G \colon \cG \arr \cG'$ are equivalences. Then there an equivalence of categories
   \[
   \hom_{\cC}(\cF, \cG) \arr
   \hom_{\cC}(\cF', \cG')
   \]
that sends each $\Phi \colon \cF \arr \cG$ into the composite
   \[
   G \circ \Phi \circ F \colon \cF' \arr \cG'.
   \]
\end{proposition}

The proof is left as an exercise to the reader.

The following is the basic criterion for checking whether a morphism of fibered categories is an equivalence.

\begin{proposition}\label{prop:char-equivalence} Let $F \colon \cF \arr \cG$ be a morphism of fibered categories. Then $F$ is an equivalence if and only if the restriction $F_U \colon \cF(U) \arr \cG(U)$ is an equivalence of categories for any object $U$ of $\cC$.
\end{proposition}

\begin{proof} Suppose that $G \colon \cG \arr \cF$ is an inverse to $F$; the two isomorphisms $F \circ G\simeq \id_{\cG}$ and $G \circ F\simeq \id_{\cF}$  restrict to isomorphisms $F_U \circ G_U\simeq \id_{\cG(U)}$ and $G_U \circ F_U\simeq \id_{\cF(U)}$, so $G_U$ is an inverse to $F_U$.

Conversely, we assume that $F_U \colon \cF(U) \arr \cG(U)$ is an equivalence of categories for any object $U$ of $\cC$, and construct an inverse $G \colon \cG \arr \cF$. Here is the main fact that we are going to need.

\begin{lemma}\label{lem:fully-fibered} Let $F \colon \cF \arr \cG$ be a morphism of fibered categories such that every restriction $F_U \colon \cF(U) \arr \cG(U)$ is fully faithful. Then the functor $F$ is fully faithful.
\end{lemma}

\begin{proof} We need to show that, given two objects $\xi'$ and $\eta'$ of $\cF$ and an arrow $\phi \colon F \xi' \arr F \eta'$ in $\cG$, there is a unique arrow $\phi' \colon \xi' \arr \eta'$ in $\cF$ with $F \phi' = \phi$. Set $\xi = F \xi'$ and $\eta = F \eta'$. Let $\eta'_1 \arr \eta'$ be a pullback of $\eta'$ to $U$, $\eta_1 = F \eta'_1$. Then the image $\eta_1 \arr \eta$ of $\eta_1' \arr \eta'$ is cartesian, so every morphism $\xi \arr \eta$ factors uniquely as $\xi \arr \eta_1 \arr \eta$, where the arrow $\xi \arr \eta_1$ is in $\cG(U)$. Analogously all arrows $\xi' \arr \eta'$ factor uniquely through $\eta'_1$; since every arrow $\xi \arr \eta_1$ in $\cG(U)$ lifts uniquely to an arrow $\xi' \arr \eta'_1$ in $\cF(U)$, we have proved the Lemma.
\end{proof}

For any object $\xi$ of $\cG$ pick an object $G \xi$ of $\cF(U)$, where $U = \p_{\cG} \xi$, together with an isomorphism $\alpha_\xi \colon \xi \simeq F(G \xi)$ in $\cG(U)$; these $G \xi$ and $\alpha_\xi$ exist because $F_U \colon \cF(U) \arr \cG(U)$ is an equivalence of categories. 

Now, if $\phi\colon \xi \arr \eta$ is an arrow in $\cG$, by the Lemma there is a unique arrow $G \phi \colon G \xi \arr G \eta$ such that $F(G \phi) = \alpha_\eta\circ \phi \circ \alpha_\xi^{-1}$, that is, such that the diagram
   \[
   \xymatrix@C+8pt{{}\xi\ar[r]^{\phi}\ar[d]^{\alpha_\xi}& 
   {}\eta\ar[d]^{\alpha_\eta}\\
   F(G \xi)\ar[r]^{F(G \phi)}&F(G \eta)}
   \] commutes.

These operations define a functor $G \colon  \cG \arr \cF$. It is immediate to check that by sending each object $\xi$ to the isomorphism $\alpha_\xi \colon \xi \simeq F(G \xi)$ we define an isomorphism of functors $\id_{\mathcal F} \simeq F \circ G \colon \cG \arr \cG$. 

We only have left to check that $G \circ F \colon \cF \arr \cF$ is isomorphic to the identity $\id_{\mathcal F}$. 

Fix an object $\xi'$ of $\cF$ over an object $U$ of $\cC$; we have a canonical isomorphism $\alpha_{F \xi'} \colon F \xi' \simeq F(G(F \xi'))$ in $\cG(U)$. Since $F_U$ is fully faithful there is a unique isomorphism $\beta_{\xi'} \colon \xi' \simeq G(F \xi')$ in $\cF(U)$ such that $F \beta_{\xi'} = \alpha_{F \xi'}$; one checks easily that this defines an isomorphism of functors $\beta\colon G \circ F \simeq \id_{\mathcal G}$.
\end{proof}

\subsection{Categories fibered in equivalence relations}\label{subsec:fibered-eq-relations}

As we remarked in \S\ref{sec:fiberedsets}, the notion of category generalizes the notion of set.

It is also possible to characterize the categories that are equivalent to a set: these are the equivalence relations.

Suppose that $R \subseteq X \times X$ is an equivalence relation on a set $X$. We can produce a category $(X,R)$ in which $X$ is the set of objects, $R$ is the set of arrows, and the source and target maps $R \arr X$ are given by the first and second projection. Then given $x$ and $y$ in $X$, there is precisely one arrow $(x,y)$ if $x$ and $y$ are in the same equivalence class, while there is none if they are not. Then transitivity assures us that we can compose arrows, while reflexivity tell us that over each object $x \in X$ there is a unique arrow $(x,x)$, which is the identity. Finally symmetry tells us that any arrow $(x,y)$ has an inverse $(y,x)$. So, $(X, R)$ is groupoid such that from a given object to another there is at most one arrow.

Conversely, given a groupoid such that from a given object to another there is at most one arrow, if denote by $X$ the set of objects and by $R$ the set of arrows, the source and target maps induce an injective map $R \arr X \times X$, which gives an equivalence relation on $X$.

So an equivalence relation can be thought of as a groupoid such that from a given object to another there is at most one arrow. Equivalently, an equivalence relation is a groupoid in which the only arrow from an object to itself is the identity.

\begin{proposition}\label{prop:equivalentset}
A category is equivalent to a set if and only if it is an equivalence relation.
\end{proposition}

\begin{proof}
If a category is equivalent to a set, it is immediate to see that it is an equivalence relation. If $(X, R)$  is an equivalence relation and $X/R$ is the set of isomorphism classes of objects, that is, the set of equivalence classes, one checks easily that the function $X \arr X/R$ gives a functor that is fully faithful and essentially surjective, so it is an equivalence.
\end{proof}

There is an analogous result for fibered categories.

\begin{definition} A category $\cF$ over $\cC$ is a \emph{quasi-functor}\index{quasi-functor}, or is \emph{fibered in equivalence relations}\index{category!fibered!in equivalence relations}, if it is fibered, and each fiber $\cF(U)$ is an equivalence relation.
\end{definition}

We have the following characterization of quasi-functors.

\begin{proposition}\label{prop:characterization-quasifunctors} A category $\cF$ over $\cC$ is a quasi-functor if and only if the following two conditions hold.

\begin{enumeratei}
\item Given an object $\eta$ of $\cF$ and an arrow $f \colon U \arr \p_ \cF \eta$ of $\cC$, there exists an arrow $\phi \colon \xi \arr \eta$ of $\cF$ with $\p_{\cF} \phi = f$. Furthermore, given any other arrow $\phi'\colon \xi' \arr \eta$ with $\p_{\cF}\phi' = f$, there exists $\alpha\colon \xi \arr \xi'$ in $\cF(U)$ such that $\phi'\alpha = \phi$.

\item  Given two objects $\xi$ and $\eta$ of $\cF$ and an arrow $f \colon \p_{\cF} \xi \arr \p_{\cF} \eta$ of $\cC$, there exists at most one arrow $\xi \arr \eta$ over $f$.
\end{enumeratei}
\end{proposition}

The easy proof is left to the reader.

\begin{proposition}\label{prop:equivalent-functor}
\index{category!fibered!in equivalence relations, characterization}
\index{fibered!category!in equivalence relations, characterization}
\index{quasi-functor!characterization}
A fibered category over $\cC$ is a quasi-functor if and only if it is equivalent to a functor.
\end{proposition}

\begin{proof} This is an application of Proposition~\ref{prop:char-equivalence}.

Suppose that a fibered category $\cF$ is equivalent to a functor $\Phi$; then every category $\cF(U)$ is equivalent to the set $\Phi U$, so $\cF$ is fibered in equivalence relations over $\cC$ by Proposition~\ref{prop:equivalentset}.

Conversely, assume that $\cF$ is fibered in equivalence relations. In particular it is fibered in groupoids, so every arrow in $\cF$ is cartesian, by Proposition~\ref{prop:char-fibered-in-groupoids}. For each object $U$ of $\cC$, denote by $\Phi U$ the set of isomorphism classes of elements in $\cF(U)$. Given an arrow $f \colon U \arr V$ in $\cC$, two isomorphic objects $\eta$ and $\eta'$ of $\cF(V)$, and two pullbacks $\xi$ and $\xi'$ of $\eta$ and $\eta'$ to $\cF(U)$, we have that $\xi$ and $\xi'$ are isomorphic in $\cF(U)$; this gives a well defined function $f^* \colon \Phi V \arr \Phi U$ that sends an isomorphism class $[\eta]$ in $\cF(V)$ into the isomorphism class of pullbacks of $\eta$. It is easy to see that this gives $\Phi$  the structure of a functor $\cC\op \arr \catset$. If we think of $\Phi$ as a category fibered in sets, we get by construction a morphism $\cF \arr \Phi$. Its restriction $\cF(U) \arr \Phi U$ is an equivalence for each object $U$ of $\cC$, so by Proposition~\ref{prop:char-equivalence} the morphism $\cF \arr \Phi$ is an equivalence.
\end{proof}

Here are a few useful facts.

\begin{proposition}\call{prop:->hom}\hfil
\begin{enumeratei}

\itemref{1} If $\cG$ is fibered in groupoids, then $\hom_\cC(\cF, \cG)$ is a groupoid.

\itemref{2} If $\cG$ is a quasi-functor, then $\hom_\cC(\cF, \cG)$ is an equivalence relation.

\itemref{3} If $\cG$ is a functor, then $\hom_\cC(\cF, \cG)$ is a set.

\end{enumeratei}
\end{proposition}

We leave the easy proofs to the reader.

In 2-categorical terms, part \refpart{prop:->hom}{3} says that the 2-category of categories fibered in sets is in fact just a 1-category, while part \refpart{prop:->hom}{2} says that the 2-category of quasi-functors is equivalent to a 1-category.

\section{Objects as fibered categories and the 2-Yoneda Lemma}

\subsection{Representable fibered categories}

In \S\ref{sec:repfunctors} we have seen how we can embed a category $\cC$ into the functor category $\func{\cC}$, while in \S\ref{sec:fiberedsets} we have seen how to embed the category $\func{\cC}$ into the 2-category of fibered categories over $\cC$. By composing these embeddings we have embedded $\cC$ into the 2-category of fibered categories: an object $X$ of $\cC$ is sent to the fibered category $(\cC/X) \arr \cC$. Furthermore, an arrow $f \colon X \arr Y$ goes to the morphism of fibered categories $(\cC/f) \colon (\cC/X) \arr (\cC/Y)$ that sends an object $U \arr X$ of $(\cC/X)$ to the composite $U \to X \overset f \to Y$. The functor $(\cC/f)$ sends an arrow
   \[
   \xymatrix@-15pt{
   U\ar[rd]\ar[rr]&&V\ar[ld]\\
   &X&
   }
   \]
of $(\cC/X)$ to the commutative diagram obtained by composing both sides with $f \colon X \arr Y$. 

This is the 2-categorical version of the weak Yoneda lemma.

\begin{named}{The weak 2-Yoneda Lemma}
\index{2-Yoneda Lemma!weak version}\index{Lemma!2-Yoneda!weak version}
The function that sends each arrow $f \colon X \arr Y$ to the morphism $(\cC/f) \colon (\cC/X) \arr (\cC/Y)$ is a bijection.
\end{named}

\begin{definition}
A fibered category over $\cC$ is \emph{representable}\index{category!fibered!representable}\index{representable!fibered category} if it is equivalent to a category of the form $(\cC/X)$.
\end{definition}

So a representable category is necessarily a quasi-functor, by Proposition~\ref{prop:equivalent-functor}. However, we should be careful: if $\cF$ and $\cG$ are fibered categories, equivalent to $(\cC/X)$ and $(\cC/Y)$ for two objects $X$ and $Y$ of $\cC$, then
   \[
   \hom(X, Y) = \hom\bigl( (\cC/X),
   (\cC/Y)\bigr),
   \]
and according to Proposition~\ref{prop:equiv->equivalenthom} we have an equivalence of categories
   \[
   \hom\bigl( (\cC/X), (\cC/Y)\bigr)
   \simeq \hom_{\cC}(\cF, \cG);
   \]
but $\hom_{\cC}(\cF, \cG)$ need not be a set, it could very well be an equivalence relation.

\subsection{The 2-categorical Yoneda lemma}

As in the case of functors, we have a stronger version of the 2-categorical Yoneda lemma. Suppose that $\cF$ is a category fibered over $\cC$, and that $X$ is an object of $\cC$. Let there be given a morphism $F \colon (\cC/X) \arr \cF$; with this we can associate an object $F(\id_X) \in \cF(X)$. Also, to each base-preserving natural transformation $\alpha\colon F \arr G$ of functors $F, G \colon (\cC/X) \arr \cF$ we associate the arrow $\alpha_{\id_X} \colon F(\id_X) \arr G(\id_X)$. This defines a functor
   \[
   \hom_{\cC}\bigl((\cC/X), \cF\bigr)
   \arr \cF(X).
   \]

Conversely, given an object $\xi \in \cF(X)$ we get a functor $F_\xi \colon (\cC/X) \arr \cF$ as follows. Given an object $\phi\colon U \arr X$ of $(\cC/X)$, we define $F_\xi (\phi) = \phi^* \xi \in \cF(U)$; with an arrow
   \[
    \xymatrix@-6pt{
   U\ar[rd]^\phi\ar[rr]^f&&V\ar[ld]_\psi\\
   &X&
   }
   \]
in $(\cC/X)$ we associate the only arrow $\theta \colon \phi^* \xi \arr \psi^* \xi$ in $\cF(U)$ making the diagram 
   \[
   \xymatrix@R=6pt{{\phi^* \xi}{}\ar@{|->}[dd]
   \ar@{-->}[rd]_\theta
   \ar@/^/[rrd] \\
   & \psi^*\xi\ar@{|->}[dd]\ar[r]s
   &\xi\ar@{|->}[dd] \\
   U\ar[rd]_f\ar@/^/[rrd]|(.505)\hole^>>>>>\phi\\
   & V \ar[r]_\psi
   & X
   }
   \]
commutative. We leave it to the reader to check that $F_\xi$ is indeed a functor.

\begin{named}{2-Yoneda Lemma}\label{2-yoneda}
\index{2-Yoneda Lemma!strong version}\index{Lemma!2-Yoneda!strong version}
The two functors above define an equivalence of categories
   \[
   \hom_{\cC}\bigl((\cC/X), \cF\bigr)
   \simeq \cF(X).
   \]
\end{named}

\begin{proof}
To check that the composite
   \[
   \cF(X) \arr 
   \hom_{\cC}\bigl((\cC/X), \cF\bigr)
   \arr \cF(X)
   \]
is isomorphic to the identity, notice that for any object $\xi \in \cF(X)$, the composite applied to $\xi$ yields $F_\xi (\id_{X}) = \id_X^* \xi$, which is canonically isomorphic to $\xi$. It is easy to see that this defines an isomorphism of functors.

For the composite
   \[
   \hom_{\cC}\bigl((\cC/X), \cF\bigr)
   \arr \cF(X)
   \arr
   \hom_{\cC}\bigl((\cC/X), \cF\bigr)
   \]
take a morphism $F \colon (\cC/X) \arr \cF$ and set $\xi = F(\id_X)$. We need to produce a base-preserving isomorphism of functors of $F$ with $F_\xi$. The identity $\id_X$ is a terminal object in the category $(\cC/X)$, hence for any object $\phi \colon U \arr X$ there is a unique arrow $\phi \arr \id_X$, which is clearly cartesian. Hence it will remain cartesian after applying $F$, because $F$ is a morphism: this means that $F(\phi)$ is a pullback of $\xi = F(\id_X)$ along $\phi \colon U \arr X$, so there is a canonical isomorphism $F_\xi(\phi) = \phi^*\xi \simeq F(\phi)$ in $\cF(U)$. It is easy to check that this defines a base-preserving isomorphism of functors, and this ends the proof.
\end{proof}

We have identified an object $X$ with the functor $\h_X \colon \cC\op \arr \catset$ it represents, and we have identified the functor $\h_X$ with the corresponding category $(\cC/X)$: so, to be consistent, we have to identify $X$ and $(\cC/X)$. So, we will write $X$ for $(\cC/X)$.

As for functors, the strong form of the 2-Yoneda Lemma can be used to reformulate the condition of representability. A morphism $(\cC/X) \arr \cF$ corresponds to an object $\xi\in \cF(X)$, which in turn defines the functor $F'\colon (\cC/X) \arr \cF$ described above; this is isomorphic to the original functor $F$. Then $F'$ is an equivalence if and only if for each object $U$ of $\cC$ the restriction
   \[
   F'_U \colon \hom_{\cC}(U, X) = (\cC/X)(U)
   \arr \cF(U)
   \]
that sends each $f \colon U \arr X$ to the pullback $f^* \xi \in \cF(U)$, is an equivalence of categories. Since $\hom_{\cC}(U, X)$ is a set, this is equivalent to saying that $\cF(U)$ is a groupoid, and each object of $\cF(U)$ is isomorphic to the image of a unique element of $\hom_{\cC}(U,X)$ via a unique isomorphism. Since the isomorphisms $\rho \simeq f^* \xi$ in $U$ correspond to cartesian arrows $\rho \arr \xi$, and in a groupoid all arrows are cartesian, this means that $\cF$ is fibered in groupoids, and for each $\rho\in \cF(U)$ there exists a unique arrow $\rho \arr \xi$. We have proved the following.

\begin{proposition}
\index{representable!fibered category!characterization}
\index{fibered!category!representable, characterization}
\index{category!fibered!representable, characterization}
A fibered category $\cF$ over $\cC$ is representable if and only if $\cF$ is fibered in groupoids, and there is an object $U$ of $\cC$ and an object $\xi$ of $\cF(U)$, such that for any object $\rho$ of $\cF$ there exists a unique arrow $\rho \arr \xi$ in $\cF$.
\end{proposition}

\subsection{Splitting a fibered category}

As we have seen in Example~\ref{ex:fibered-nosplitting}, a fibered category does not necessarily admit a splitting. However, a fibered category is always equivalent to a split fibered category.

\begin{theorem}\label{thm:equivalent-split}
\index{theorem!existence of a split fibered category equivalent to a given fibered category}
Let $\cF \arr \cC$ be a fibered category. Then there exists a canonically defined split fibered category $\widetilde{\cF} \arr \cC$ and an equivalence of fibered categories of $\widetilde{\cF}$ with $\cF$.
\end{theorem}

\begin{proof}
In this proof, if $U$ is an object of $\cC$, we will identify the functor $\h_{U}$ with the comma category $(\cC/U)$. We have a functor $\hom(-, \cF) \colon \cC\op \arr \cat{Cat}$ from the category $\cC\op$ into the category of all categories. If $U$ is an object of $\cC$ this functor will send $U$ into the category $\hom_{\cC}(\h_{U}, \cF)$ of base-preserving natural transformations. An arrow $U \arr V$ corresponds to a natural transformation $\h_{U} \arr \h_{V}$, and this induces a functor $\hom_{\cC}(\h_{V}, \cF) \arr \hom_{\cC}(\h_{U}, \cF)$.

Let us denote by $\widetilde{\cF}$ the fibered category associated with this functor: by definition, $\widetilde{\cF}$ comes with a splitting. There is an obvious morphism $\widetilde{\cF} \arr \cF$, sending an object $\phi \colon \h_{U} \arr \cF$ into $\phi(\id_{U}) \in \cF(U)$. According to the 2-Yoneda Lemma, and the criterion of Proposition~\ref{prop:char-equivalence}, this is an equivalence.
\end{proof}

It is an interesting exercise to figure out what this construction yields in the case of a surjective group homomorphism $G \arr H$, as in Example~\ref{ex:fibered-nosplitting}.

\section{The functors of arrows of a fibered category}
\label{sec:func-arrows}

Suppose that $\cF \arr \cC$ is a fibered category; if $U$ is an object in $\cC$ and $\xi$, $\eta$ are objects of $\cF(U)$, we denote by $\hom_U(\xi, \eta)$ the set of arrows from $\xi$ to $\eta$ in $\cF(U)$.

Let $\xi$ and $\eta$ be two objects of $\cF$ over the same object $S$ of $\cC$. Let $u_1 \colon U_1 \arr S$  and $u_2 \colon U_2 \arr S$ be arrows in $\cC$; these are objects of the comma category $(\cC/S)$. Suppose that $\xi_i \arr \xi$ and $\eta_i \arr \eta$ are pullbacks along $u_i \colon U_i \arr S$ for $i=1$, 2. For each arrow $f \colon U_1 \arr U_2$ in $(\cC/S)$, by definition of pullback there are two arrows, each unique,  $\alpha_f \colon \xi_1 \arr \xi_2$ and $\beta_f \colon \eta_1 \arr \eta_2$, such that and the two diagrams
   \[
   \xymatrix@-15pt{\xi_1\ar[rr]^{\alpha_f}\ar[rd]
   && \xi_2 \ar[ld]\\ &\xi &}\qquad
   \text{and}\qquad
   \xymatrix@-15pt{\eta_1 \ar[rr]^{\beta_f}\ar[rd]
   && \eta_2\ar[ld]\\ &\eta &}
   \]
commute. By Proposition~\refall{prop:comp-cartesian}{4} the arrows $\alpha_f$ and $\beta_f$ are cartesian; we define a pullback function
   \[
   f^* \colon \hom_{U_2}(\xi_2, \eta_2) \arr
   \hom_{U_1}(\xi_1, \eta_1)
   \]
in which $f^* \phi$ is defined as the only arrow $f^* \phi \colon \xi_1 \arr \eta_1$ in $\cF(U_1)$ making the diagram
   \[
   \xymatrix{
   \xi_1\ar[r]^{f^*\phi}\ar[d]^{\alpha_f}
   & \eta_1\ar[d]^{\beta_f}\\
   \xi_2\ar[r]^{\phi} & \eta_2
   }
   \]
commute. If we are given a third arrow $g \colon U_2 \arr U_3$ in $(\cC/S)$ with pullbacks $\xi_3\arr \xi$ and $\eta_3 \arr \eta$, we have arrows $\alpha_g \colon \xi_2 \arr \xi_3$ and $\beta_g \colon \eta_2 \arr \eta_3$; it is immediate to check that
   \[
   \alpha_{gf} = \alpha_g \circ \alpha_f
   \colon \xi_1 \arr \xi_3\quad
   \text{and}\quad
   \beta_{gf} = \beta_g \circ \beta_f
   \colon \eta_1 \arr \eta_3
   \]
and this implies that
   \[
   (gf)^* = f^* g^* \colon \hom_{U_3}(\xi_3, \eta_3)
   \arr \hom_{U_1}(\xi_1, \eta_1).
   \]

After choosing a cleavage for $\cF$, we can define a functor
   \[
   \underhom_S(\xi, \eta) \colon (\cC/S)\op \arr
   \catset
   \]
by sending each object $u \colon U \arr S$ into the set $\hom_U(u^* \xi, u^* \eta)$ of arrows in the category $\cF(U)$. An arrow $f \colon U_1 \arr U_2$ from $u_1 \colon U_1 \arr S$ to $u_2 \colon U_2 \arr S$ yields a function
   \[
   f^* \colon \hom_{U_2}(u^*_2\xi, u^*_2\eta)
   \arr \hom_{U_1}(u^*_1\xi, u^*_1\eta);
   \]
and this defines the effect of $\underhom_S(\xi, \eta)$ on arrows.

It is easy to check that the functor $\underhom_S(\xi, \eta)$ is independent of the choice of a cleavage, in the sense that cleavages give canonically isomorphic functors. Suppose that we have chosen for each $f \colon U \arr V$ and each object $\zeta$ in $\cF(V)$ another pullback $f^\vee \zeta \arr \zeta$: then there is a canonical isomorphism $u^* \eta \simeq u^\vee \eta$ in $\cF(U)$ for each arrow $u \colon U \arr S$, and this gives a bijective correspondence
   \[
   \hom_U(u^* \xi, u^* \eta) \simeq \hom_U(u^\vee
   \xi, u^\vee \eta),
   \]
yielding an isomorphism of the functors of arrows defined by the two pullbacks.

In fact, $\underhom_S(\xi, \eta)$ can be more naturally defined as a quasi-functor
   \[
   \curshom_S(\xi, \eta) \arr (\cC/S);
   \]
this does not require any choice of cleavages.

From this point of view, the objects of $\curshom_S(\xi, \eta)$ over some object $u \colon U \arr S$ of $(\cC/S)$ are triples
   \[
   (\xi_1 \arr \xi, \eta_1 \arr \eta, \phi),
   \]
where $\xi_1 \arr \xi$ and $\eta_1 \arr \eta$ are cartesian arrows of $\cF$ over $u$, and $\phi \colon \xi_1 \arr \eta_1$ is an arrow in $\cF(U)$. An arrow from $(\xi_1 \arr \xi, \eta_1 \arr \eta, \phi_1)$ over $u_1 \colon U_1 \arr S$ and $(\xi_2 \arr \xi, \eta_2 \arr \eta, \phi_2)$ over $u_2 \colon U_2 \arr S$ is an arrow $f \colon U_1 \arr U_2$ in $(\cC/S)$ such that $f^*\phi_2 = \phi_1$.

From Proposition~\ref{prop:characterization-quasifunctors} we see that $\curshom_S(\xi, \eta)$ is a quasi-functor over $\cC$, and therefore, by Proposition~\ref{prop:equivalent-functor}, it is equivalent to a functor: of course this is the functor $\underhom_S(\xi, \eta)$ obtained by the previous construction.

This can be proved as follows: the objects of $\underhom_S(\xi, \eta)$, thought of as a category fibered in sets over $(\cC/S)$ are pairs $(\phi, u \colon U \arr S)$, where $u \colon U \arr S$ is an object of $(\cC/S)$ and $\phi \colon u^* \xi \arr u^* \eta$ is an arrow in $\cF(U)$; this also gives an object $(u^* \xi \arr \xi, u^*\eta \arr \eta, \phi)$ of $\curshom_S(\xi, \eta)$ over $U$. The arrows between objects of $\underhom_S(\xi, \eta)$ are precisely the arrows between the corresponding objects of $\curshom_S(\xi, \eta)$, so we have an embedding of $\underhom_S(\xi, \eta)$ into $\curshom_S(\xi, \eta)$. But every object of $\curshom_S(\xi, \eta)$ is isomorphic to an object of $\underhom_S(\xi, \eta)$, hence the two fibered categories are equivalent.

\section{Equivariant objects in fibered categories}
\label{sec:fibered-actions}

The notion of an equivariant sheaf of modules on a scheme with the action of a group scheme, as defined in \cite[Chapter~1, \S~3]{git}, or in \cite{thomason87}), is somewhat involved and counterintuitive. The intuition is that if we are given the action of a group scheme $G$ on a scheme $X$, an equivariant sheaf should be a sheaf $\cF$, together with an action of $G$ on the \emph{pair} $(X, \cF)$, which is compatible with the action of $G$ on $X$. Since the pair $(X, \cF)$ is an object of the fibered category of sheaves of modules, the language of fibered categories is very well suited for expressing this concept.

Let $G \colon \cC\op \arr \catgrp$ be a functor, $\cF \arr \cC$ a fibered category, $X$ an object of $\cC$ with an action of $G$ (see \S\ref{subsec:actions}).

\begin{definition}\call{def:equivariant-object}
A \emph{$G$-equivariant object}\index{equivariant!object} of $\cF(X)$ is an object $\rho$ of $\cF(X)$, together with an action of $G(U)$ on the set $\hom_{\cF}(\xi, \rho)$ for any $\xi \in \cF(U)$, such that the following two conditions are satisfied.
\begin{enumeratei}

\itemref{1} For any arrow  $\phi \colon \xi \arr \eta$ of $\cF$ mapping to an arrow $f \colon U \arr V$, the induced function $\phi^* \colon \hom_{\cF}(\eta, \rho) \arr \hom_{\cF}(\xi, \rho)$ is equivariant with respect to the group homomorphism $f^* \colon G(V) \arr G(U)$.

\itemref{2} The function $ \hom_{\cF}(\xi, \rho) \arr  \hom_{\cC}(U, X)$ induced by $\p_{\cF}$ is $G(U)$-equivariant.

\end{enumeratei}
An arrow $u \colon \rho \arr \sigma$ in $\cF(X)$ is \emph{$G$-equivariant} if it has the property that the induced function $\hom_{\cF}(\xi, \rho) \arr \hom_{\cF}(\xi, \sigma)$ is $G(U)$-equivariant for all $U$  and all $\xi\in \cF(U)$.
\end{definition}

The first condition can be expressed by saying that the data define an action of $G \circ \p_{\cF} \colon \cF\op \arr \catgrp$ on the object $\rho$ of $\cF$. In other words, for any object $\xi$ of $\cF$, the action defines a set theoretic action
   \[
   (G \circ \p_{\cF})(\xi) \times \h_{\rho}(\xi)
   \arr \h_{\rho}(\xi)
   \]
and this action is required to give a natural transformation of functors $\cF\op \arr \catset$
   \[
   (G \circ \p_{\cF})\times \h_{\rho}
   \arr \h_{\rho}.
   \]

The second condition can be thought of as saying that the action of $G$ on $\rho$ is compatible with the action of $G$ on $X$.

The $G$-equivariant objects over $X$ are the objects of a category $\cF^G(X)$, in which the arrows are the equivariant arrows in $\cF(X)$. 

It is not hard to define the fibered category of $G$-equivariant objects of $\cF$ over the category of $G$-equivariant objects of $\cC$, but we will not do this.

Now assume that $G$ is a group object in $\cC$ acting on an object $X$ of $\cC$, corresponding to  an arrow $\alpha \colon G \times X \arr X$, as in Proposition~\ref{prop:diagram-action}. Take a fibered category $\cF \arr \cC$: the category $\cF^{G}(X)$ of equivariant objects over $X$ has a different description. Choose a cleavage for $\cF$.

Let $\rho$ be an object of $\cF^{G}(X)$. Consider the pullback $\pr_2^*\rho \in \cF(G \times X)$, and the functor $\h_{\pr_2^*\rho} \colon \cF \arr \catset$ it represents. If $\phi \colon \xi \arr \pr_2^* \rho$ is an arrow in $\cF$, we obtain an arrow $\xi \arr \rho$ by composing $\phi$ with the given cartesian arrow $\pr_2^* \rho \arr \rho$, and an arrow $\p_{\cF}\xi \arr G$ by composing $\p_{\cF}\phi \colon \p_{\cF} \xi \arr G \times X$ with the projection $G \times X \arr G$. This defines a natural transformation $\h_{\pr_2^*\rho} \arr (G \circ \p_{\cF}) \times \h_\rho$.

The fact that the canonical arrow $p \colon \pr_2^* \rho \arr \rho$ is cartesian implies that each pair consisting of an arrow $\xi \arr \rho$ in $\cF$ and an arrow $\p_{\cF}\xi \arr G$ in $\cC$ comes from a unique arrow $\xi \arr \pr_{2}^{*}\rho$. This means that the natural transformation above is in fact an isomorphism $\h_{\pr_2^*\rho} \simeq (G \circ \p_{\cF}) \times \h_\rho$ of functors $\cF\op \arr \catset$. Hence, by Yoneda's lemma, a morphism $(G \circ \p_{\cF}) \times \h_\rho \arr \h_\rho$ corresponds to an arrow $\beta\colon \pr_2^*\rho \arr \rho$. Condition~\refpart{def:equivariant-object}{2} of Definition~\ref{def:equivariant-object} can be expressed as saying that $\p_{\cF}\beta = \alpha$.

There are two conditions that define an action. First consider the natural transformation $\h_{\rho} \arr (G \circ \p_{\cF}) \times \h_{\rho}$ that sends an object $u \in  \h_{\rho}(\xi)$ to the pair $(1, u) \in G(\p_{\cF}\xi) \times \h_{\rho}(\xi)$; this corresponds to an arrow $\epsilon_{\rho} \colon \rho \arr \pr_{2}^{*}\rho$, whose composite with $p \colon  \pr_{2}^{*}\rho \arr \rho$ is the identity $\id_{\rho}$, and whose image in $\cC$ is the arrow $\epsilon_{X}\colon X = \pt \times X \arr G \times X$ induced by $\iden_{G} \colon \pt \arr G$. Since $p$ is cartesian, these two conditions characterize $\epsilon_{\rho}$ uniquely.

The first condition that defines an action of $(G \circ \p_{\cF})$ on $\rho$ (see Proposition~\ref{prop:diagram-action}) is that the composite $\h_{\rho} \arr (G \circ \p_{\cF}) \times \h_{\rho} \arr \h_{\rho}$ be the identity; and this is equivalent to saying that the composite $\beta \circ \epsilon_{\rho} \colon \rho \arr \rho$ is the identity $\id_{\rho}$.

The second condition can be expressed similarly. The functor
   \[
   (G \circ \p_{\cF}) \times (G \circ \p_{\cF}) \times \h_{\rho} \colon \cF 
   \arr \catgrp
   \]
is represented by the pullback $\pr_{3}^{*}\rho$ of $\rho$ along the third projection $\pr_{3} \colon G \times G \times X \arr X$. Now, given any arrow $f \colon G \times G \times X \arr G \times X$ whose composite with $\pr_{2}\colon F \times X \arr X$ equals $\pr_{2}\colon G \times G \times X \arr X$, there is a unique arrow $\widetilde{f} \colon \pr_{3}^{*}\rho \arr \p_{2}^{*}\rho$ mapping to $f$, such that the composite $p \circ \widetilde{f} \colon  \pr_{3}^{*}\rho \arr \rho$ equals the canonical arrow $q \colon \pr_{3}^{*}\rho \arr \rho$. Then it is an easy matter to convince oneself that the second condition that defines an action is equivalent to the commutativity of the diagram
   \begin{equation}\label{eq:lifts}
   \xymatrix@C+20pt{
   {}\pr_{3}^{*}\rho \ar[r]^{\widetilde{\mul_{G}\times \id_{X}}}
      \ar[d]^{\widetilde{\id_{G}\times \alpha}} &
   {}\pr_{2}^{*}\rho \ar[d]^{\beta}\\
   {}\pr_{2}^{*}\rho \ar[r]^{\beta}&
   {}\rho
   }
   \end{equation}

This essentially proves the following fact (we leave the easy details to the reader).

\begin{proposition}
Let $\rho$ be an object of $\cF(X)$. To give $\rho$ the structure of a  $G$-equivariant object is the same as assigning an arrow $\beta \colon  \pr_{2}^{*}\rho \arr \rho$ with $\p_{\cF}\beta = \alpha$, satisfying the following two conditions.
\begin{enumeratei}

\item $\beta \circ \epsilon_{\rho} = \id_{\rho}$.

\item The diagram (\ref{eq:lifts}) commutes.

\end{enumeratei}

Furthermore, if $\rho$ and $\rho'$ are $G$-equivariant objects, and we denote by $\beta \colon \pr_{2}^{*}\rho \arr \rho$ and $\beta' \colon \pr_{2}^{*}\rho' \arr \rho'$ the corresponding arrows, then $u \colon \rho \arr \rho'$ is $G$-equivariant if and only if the diagram
   \[
   \xymatrix{
   {}\pr_{2}^{*}\rho \ar[r]^{\beta} \ar[d]^{\pr_{2}^{*}u}&
   \rho              \ar[d]^{u}  \\
   {}\pr_{2}^{*}\rho' \ar[r]^{\beta'}&
   \rho'
   }
   \]
commutes.

\end{proposition}

This can be restated further, to make it look more like the classical definition of an equivariant sheaf. First of all, let us notice that if an arrow $\beta \colon \pr_2^*\rho \arr \rho$ corresponds to a $G$-equivariant structure on $\rho$, then it is cartesian. This can be shown as follows. 

There is an automorphism $i$ of $G \times X$, defined in functorial terms by the equation $i(g,x) = (g^{-1}, gx)$ whenever $U$ is an object of $\cC$, $g \in G(U)$ and $x \in X(U)$; this has the property that 
   \[
   \pr_{2} \circ i = \alpha \colon G \times X \arr X.
   \]
Analogously one can use the action of $(G \circ \p_{\cF})$ on $\rho$ to define an automorphism $I$ of $(G \circ \p_{\cF}) \times \h_\rho$, hence an automorphism of $\pr_{2}^{*}\rho$, whose composite with the canonical arrow $p \colon \pr_{2}^{*}\rho \arr \rho$ equals $\beta$. Since $I$ is an isomorphism, hence is cartesian, the canonical arrow $p$ is cartesian, and the composite of cartesian arrow is cartesian, it follows that $\beta$ is cartesian.

Now, start from a \emph{cartesian} arrow $\beta \colon \pr_{2}^{*} \rho \arr \rho$ with $\p_{\cF}\beta = \alpha$. Assume that the diagram~(\ref{eq:lifts}) is commutative. I claim that in this case we also have $\beta \circ \epsilon_{\rho} = \id_{\rho}$. This can be checked in several ways: here is one.

The arrow $\beta$ corresponds to a natural transformation $(G \circ \p_{\cF}) \times \h_{\rho} \arr \h_{\rho}$. The commutativity of the diagram~(\ref{eq:lifts}) expresses the fact that $(g_{1}g_{2})x = g_{1}(g_{2}x)$ for any object $\xi$ of $\cF$ and any $g_{1}, g_{2} \in G(\p_{\cF}\xi)$ and $x \in \h_{\rho}\xi$.  The arrow $\beta \circ \epsilon_{\rho} \colon \rho \arr \rho$ corresponds to the natural transformation $\h_{\rho} \arr \h_{\rho}$ given by multiplication by the identity: and, because of the previous identity, this is an idempotent endomorphism of $\h_{\rho}$. Hence $\beta \circ \epsilon_{\rho}$ is an idempotent arrow in $\hom_{\cF}(\rho, \rho)$.

On the other hand, Proposition~\refall{prop:comp-cartesian}{4} implies that $\epsilon_{\rho}$ is a cartesian arrow, so $\beta \circ \epsilon_{\rho}$ is also cartesian. But $\beta \circ \epsilon_{\rho}$ maps to $\id_{X}$ in $\cC$, hence is an isomorphism: and the only idempotent isomorphism is the identity.

This allows us to rewrite the conditions as follows.

\begin{proposition}
Let $\rho$ be an object of $\cF(X)$. To give $\rho$ the structure of a  $G$-equivariant object is the same as assigning a cartesian arrow $\beta \colon  \pr_{2}^{*}\rho \arr \rho$ with $\p_{\cF}\beta = \alpha$, such that the diagram (\ref{eq:lifts}) commutes.

Furthermore, let $\rho$ and $\rho'$ be $G$-equivariant objects, and denote by $\beta \colon \pr_{2}^{*}\rho \arr \rho$ and $\beta' \colon \pr_{2}^{*}\rho' \arr \rho'$ the corresponding arrows, Then $u \colon \rho \arr \rho'$ is $G$-equivariant if and only if the diagram
   \[
   \xymatrix{
   {}\pr_{2}^{*}\rho \ar[r]^{\beta} \ar[d]^{\pr_{2}^{*}u}&
   \rho              \ar[d]^{u}  \\
   {}\pr_{2}^{*}\rho' \ar[r]^{\beta'}&
   \rho'
   }
   \]
commutes.

\end{proposition}

A final restatement is obtained via a cleavage, in the language of pseudo-functors. Recall that an arrow $\pr_{2}^{*}\rho \arr \rho$ mapping to $\alpha$ in $\cC$ corresponds to an arrow $\beta\colon \phi \colon \pr_{2}^{*}\rho \arr \alpha^{*} \rho$ in $\cF(G \times X)$, and that $\beta$ is cartesian if and only if $\phi$ is an isomorphism.

We also have the equalities
   \[
   \xymatrix{
   \pr_{3} =
   \pr_{2} \circ (\mul_{G} \times \id_{X}) =  \pr_{2} \circ \pr_{23}
   \colon G \times G \times X \arr X,
   }
   \]
   \[
   \xymatrix{
   A \eqdef
   \alpha \circ (\mul_{G} \times \id_{X})
    = \alpha \circ (\id_{G} \times \alpha)
   \colon G \times G \times X \arr X,
   }
   \]
and
   \[
   B \eqdef 
   \pr_{2} \circ (\id_{G} \times \alpha)
    = \alpha \circ \pr_{23}
   \colon G \times G \times X \arr X.
   \]

 We leave it to the reader to unwind the various definitions and check that the following is equivalent to the previous statement.

\begin{proposition}\label{prop:definition-equivariant}
Let $\rho$ be an object of $\cF(X)$. To give $\rho$ the structure of a  $G$-equivariant object is the same as assigning an isomorphism $\phi \colon \pr_{2}^{*}\rho \simeq \alpha^{*}\rho$ in $\cF(G \times X)$, such that the diagram
   \[
   \xymatrix@R+10pt{
   {}\pr_{3}^{*}\rho
      \ar[rr]^-{(\mul_{G} \times \id_{X})^{*}\phi}
      \ar[rd]_{\pr_{23}^{*}\phi}
   && A^{*}\rho\\
   & B^{*}\rho\ar[ru]_{(\id_{G}\times \alpha)^{*} \phi}
   }
   \]
commutes.

\end{proposition}

When applied to the fibered category of sheaves of some kind (for example, \qc sheaves) one gets precisely the usual definition of an equivariant sheaf.

\chapter{Stacks}\label{ch:stacks}

\section{Descent of objects of fibered categories}

\subsection{Gluing continuous maps and topological spaces}

The  following is the archetypal example of descent.
Take $\cat{Cont}$ to be the category of continuous maps (that is, the category of arrows in $\cattop$, as in Example~\ref{ex:arrows}); this category is fibered over $\cattop$ via the functor $\p_{\cat{Cont}} \colon \cat{Cont} \arr \cattop$ sending each continuous map to its codomain. Now, suppose that $f \colon X \arr U$ and $g \colon Y \arr U$ are two objects of $\cat{Cont}$ mapping to the same object $U$ in $\cattop$; we want to construct a continuous map $\phi\colon X \arr Y$ over $U$, that is, an arrow in $\cat{Cont}(U) = \catover{Top}{U}$. Suppose that we are given an open covering $\{U_i\}$ of $U$, and continuous maps $\phi_i \colon f^{-1} U_i \arr g^{-1} U_i$ over $U_i$; assume furthermore that the restriction of $\phi_i$ and $\phi_j$ to $f ^{-1}(U_i \cap U_j) \arr g^{-1}(U_i \cap U_j)$ coincide. Then there is a unique continuous map $\phi \colon X \arr Y$ over $U$ whose restriction to each $f^{-1}U_i$ coincides with $\phi_i$.

This can be written as follows. The category $\cat{Cont}$ is fibered over $\cattop$, and if $f\colon V \arr U$ is a continuous map, $X \arr U$ an object of $\cat{Cont} (U) = \catover{Top}{U}$, then a pullback of $X \arr U$ to $V$ is given by the projection $V \times_U X \arr V$. The functor $f^* \colon \cat{Cont}(U) \arr \cat{Cont}(V)$ sends each object $X \arr U$ to $V \times_U X \arr V$, and each arrow in $\catover{Top}{U}$, given by  continuous function $\phi\colon X \arr Y$ over $U$, to the continuous function $f^* \phi = \id_V \times_U f \colon V \times_U X \arr V \times_U Y$. 

Suppose that we are given two topological spaces $X$ and $Y$ with continuous maps $X \arr S$ and $Y \arr S$. Consider the functor
   \[
   \underhom_S(X, Y) \colon \catover{Top}{S} \arr \catset
   \]
from the category of topological spaces over $S$, defined in Section~\ref{sec:func-arrows}. This sends each arrow $U \arr S$ to the set of continuous maps $\hom_U(U \times_S X, U \times_S Y)$ over $U$. The action on arrows is obtained as follows: given a continuous function $f \colon V \arr U$ over $S$, we send each continuous function $\phi \colon U \times_S X \arr U \times_S Y$ to the function
   \[
   f^* \phi = \id_V \times \phi \colon V
   \times_S X = V \times_U(U \times_S X) \arr V \times_U(U
   \times_S Y) = V \times_S Y.
   \]

Then the fact that continuous functions can be constructed locally and then glued together can be expressed by saying that the functor
   \[
   \underhom_S(X, Y) \colon \catover{Top}{S}
   \op \arr \catset
   \]
is a sheaf in the classical topology of $\cattop$.

But there is more: not only can we construct continuous functions locally: we can also do this for spaces, although this is more complicated.

\begin{proposition}\label{prop:gluing-top}
Suppose that we are given a topological space $U$ with an open covering $\{U_i\}$; for each triple of indices $i$, $j$ and $k$ choose fibered products $U_{ij} = U_i \cap U_j$ and $U_{ijk} = U_i \cap U_j \cap U_k$. Assume that for each $i$ we have a continuous map $u_i \colon X_i \arr U_i$, and that for each pair of indices $i$ and $j$ we have a homeomorphism $\phi_{ij} \colon u_j^{-1}U_{ij} \simeq u_i^{-1}U_{ij}$ over $U_{ij}$, satisfying the cocycle condition
   \[
   \phi_{ik} = \phi_{ij} \circ \phi_{jk} \colon 
   u_k^{-1}U_{ijk} \arr 
   u_j^{-1}U_{ijk} \arr
   u_i^{-1}U_{ijk}.
   \]
Then there exists a continuous map $u \colon X \arr U$, together with isomorphisms $\phi_i \colon u^{-1} U_i \simeq X_i$, such that $\phi_{ij} = \phi_i \circ \phi_j^{-1} \colon u_j^{-1}U_{ij} \arr u^{-1}U_{ij} \arr u_iU_{ij}$ for all $i$ and $j$.
\end{proposition}

\begin{proof}
Consider the disjoint union $U'$ of the $U_i$; the fibered product $U' \times_U U'$ is the disjoint union of the $U_{ij}$. The disjoint union $X'$ of the $X_i$, maps to $U'$; consider the subset $R \subseteq X' \times X'$ consisting of pairs $(x_i, x_j) \in X_i \times X_j \subseteq X' \times X'$ such that $x_i = \phi_{ij} x_j$. I claim that $R$ is an equivalence relation in $X'$. Notice that the cocycle condition $\phi_{ii} = \phi_{ii} \circ \phi_{ii}$ implies that $\phi_{ii}$ is the identity on $X_i$, and this shows that the equivalence relation is reflexive. The fact that $\phi_{ii} = \phi_{ij} \circ \phi_{ji}$, and therefore $\phi_{ji} = \phi_{ij}^{-1}$, prove that it is symmetric; and transitivity follows directly from the general cocycle condition. We define $X$ to be the quotient $X'/R$.

If two points of $X'$ are equivalent, then their images in $U$ coincide; so there is an induced continuous map $u \colon X \arr U$. The restriction to $X_i \subseteq X'$ of the projection $X'\arr X$ gives a continuous map $\phi_i \colon X_i \arr u^{-1}U_i$, which is easily checked to be a homeomorphism. One also sees that $\phi_{ij} = \phi_i \circ \phi_j^{-1}$, and this completes the proof.
\end{proof}

The fact that we can glue continuous maps and topological spaces says that $\cat{Cont}$ is a \emph{stack} over $\cattop$.

\subsection{The category of descent data}\label{subsec:descent-data}

Let $\cC$ be a site. We have seen that a fibered category over $\cC$ should be thought of as a functor from $\cC$ to the category of categories, that is, as a presheaf of categories over $\cC$. A stack is, morally, a sheaf of categories over $\cC$.

Let $\cF$ be a category fibered over $\cC$. We fix a cleavage; but we will also indicate how the definitions can be given without resorting to the choice of a cleavage.

Given a covering $\{\sigma \colon U_i \arr U\}$, set $U_{ij} = U_i \times_U U_j$ and $U_{ijk} = U_i \times_U U_j \times_U U_k$ for each triple of indices $i$, $j$ and $k$.

\begin{definition}

Let $\cU = \{\sigma_i \colon U_i \arr U\}$ be a covering in $\cC$. An \emph{object with descent data}\index{descent data!object with} $(\{\xi_i\}, \{\phi_{ij}\})$ on $\cU$, is a collection of objects $\xi_i \in \cF(U_i)$, together with isomorphisms $\phi_{ij} \colon \pr_2^* \xi_j \simeq \pr_1^* \xi_i$ in $\cF(U_i \times_U U_j)$, such that the following cocycle condition is satisfied.

For any triple of indices $i$, $j$ and $k$, we have the equality
   \[
   \mathrm{pr}_{13}^* \phi_{ik} = \mathrm{pr}_{12}^* \phi_{ij}
   \circ \mathrm{pr}_{23}^* \phi_{jk} \colon \pr_3^*\xi_k
   \arr \pr_1^*\xi_i
   \]
where the $\mathrm{pr}_{ab}$ and $\pr_a$ are projections on the $a\th$ and $b\th$ factor, or the $a\th$ factor respectively.

The isomorphisms $\phi_{ij}$ are called \emph{transition isomorphisms}\index{transition isomorphisms} of the object with descent data.

An arrow between objects with descent data
   \[
   \{\alpha_i\} \colon (\{\xi_i\}, \{\phi_{ij}\}) \arr
   (\{\eta_i\}, \{\psi_{ij}\})
   \]
is a collection of arrows $\alpha_i\colon \xi_i \arr \eta_i$ in $\cF(U_i)$, with the property that for each pair of indices $i$, $j$, the diagram
   \[
   \xymatrix@C+15pt{
   {}\pr_2^* \xi_j \ar[r]^{\pr_2^* \alpha_j} \ar[d]^{\phi_{ij}}
   & {}\pr_2^* \eta_j\ar[d]^{\psi_{ij}}\\
   {}\pr_1^* \xi_i \ar[r]^{\pr_1^* \alpha_i}&
   {}\pr_1^* \eta_i
   }
   \]
commutes.

\end{definition}

In understanding the definition above it may be useful to contemplate the cube 
   \begin{equation}\label{eq:basic-cube}
   \xymatrix{
   &U_{ijk}\ar[rr]^{\pr_{23}}\ar'[d]^{\pr_{13}}[dd]
   \ar[ld]_{\pr_{12}}
   &&U_{jk}\ar[ld]\ar[dd]\\
   U_{ij}\ar[rr]\ar[dd]
   && U_j \ar[dd]\\
   &U_{ik} \ar[ld]\ar'[r][rr]
   &&U_k\ar[ld]\\
   U_i \ar[rr]&& U
   }
   \end{equation}
in which all arrows are given by projections, and every face is cartesian.

There is an obvious way of composing morphisms, which makes objects with descent data the objects of a category, denoted by $\cF(\cU) = \cF(\{U_i \to U\})$.

\begin{remark}
This category does not depend on the choice of fibered products $U_{ij}$ and $U_{ijk}$, in the sense that with different choices we get isomorphic categories.
\end{remark}

For each object $\xi$ of $\cF(U)$ we can construct an object with descent data on a covering $\{\sigma_i \colon U_i \arr U\}$ as follows. The objects are the pullbacks $\sigma_i^* \xi$; the isomorphisms  $\phi_{ij} \colon \pr_2^* \sigma_j^* \xi \simeq \pr_1^* \sigma_i^* \xi$ are the isomorphisms that come from the fact that both $\pr_2^* \sigma_j^* \xi$ and $\pr_1^* \sigma_i^* \xi$ are pullbacks of $\xi$ to $U_{ij}$. If we identify $\pr_2^* \sigma_j^* \xi$ with $\pr_1^* \sigma_i^* \xi$, as is commonly done, then the $\phi_{ij}$ are identities.

Given an arrow $\alpha \colon \xi \arr \eta$ in $\cF(U)$, we get arrows $\sigma_i^*\alpha \colon \sigma_i^* \xi \arr \sigma_i^* \eta$, yielding an arrow from the object with descent associated with $\xi$ to the one associated with $\eta$. This defines a functor $\cF(U) \arr \cF(\{U_i \to U\})$.

It is important to notice that these constructions do not depend on the choice of a cleavage, in the following sense. Given a different cleavage,  for each covering $\{U_i \to U\}$ there is a canonical isomorphism of the resulting categories $\cF(\{U_i \to U\})$; and the functors $\cF(U) \arr \cF(\{U_i \to U\})$ commute with these equivalences.

Here is a definition of the category of descent data that does not depend on choosing of a cleavage. Let $\{U_i \to U\}_{i \in I}$ be a covering. We define an object with descent data to be a triple of sets
   \[
   (\{\xi_i\}_{i \in I}, \{\xi_{ij}\}_{i, j\in I},
   \{\xi_{ijk}\}_{i, j, k\in I}),
   \]
where each $\xi_\alpha$ is an object of $\cF(U_{\alpha})$, plus, for each triple of indices $i$, $j$ and $k$, a commutative diagram
   \[
   \xymatrix@-15pt{
   &\xi_{ijk}\ar[rr]\ar'[d][dd] \ar[ld] &&\xi_{jk}\ar[ld]\ar[dd]\\
   \xi_{ij}\ar[dd]\ar[rr] && \xi_j\\
   &\xi_{ik} \ar[ld]\ar[rr] &&\xi_k\\
   \xi_i&&
   }
   \]
in which every arrow is cartesian, and such that when applying $\p_{\cF}$ every arrow maps to the appropriate projection in the diagram (\ref{eq:basic-cube}). These form the objects of a category $\cF\desc(\{U_i \to U\})$.

An arrow
   \[
   \{\phi_i\}_{i \in I} \colon (\{\xi_i\}, \{\xi_{ij}\},\{\xi_{ijk}\})
   \arr (\{\eta_i\}, \{\eta_{ij}\}, \{\eta_{ijk}\})
   \]
consists of set of arrows with $\phi_i\ \colon \xi_i \arr \eta_i $ in $\cF(U_i)$, such that for every pair of indices $i$ and $j$ we have
   \[
   \pr_1^* \phi_i = \pr_2 ^* \phi_j \colon 
   \xi_{ij} \arr \eta_{ij}.
   \]

Alternatively, and perhaps more naturally, we could define an arrow as a triple $((\{\phi_i\}_{i \in I}, \{\phi_{ij}\}_{i, j\in I}, \{\phi_{ijk}\}_{i, j, k\in I}))$, where $\phi_\alpha \colon \xi_\alpha \arr \eta_\alpha$ is an arrow in $\cF(U_{\alpha})$ for each $\alpha$ in $I$, $I\times I$ or $I \times I \times I$, with the obvious compatibility conditions with the various arrows involved in the definition of an object. We leave it to the reader to check that these two definitions of an arrow are equivalent.

Once we have chosen a cleavage, there is a functor from $\cF\desc (\{U_i \to U\})$ to $\cF(\{U_i \to U\})$. Given an object $(\{\xi_i\}, \{\xi_{ij}\}, \{\xi_{ijk}\})$ of $\cF\desc(\{U_i \to U\})$, the arrows $\xi_{ij} \arr \xi_i$ and $\xi_{ij} \arr\xi_j$ induce isomorphisms $\xi_{ij} \simeq \pr_1^* \xi_i$ and $\xi_{ij} \simeq \pr_2^* \xi_j$; the resulting isomorphism $\pr_2^* \xi_j \simeq \pr_1^* \xi_i$ is easily seen to satisfy the cocycle condition, thus defining an object of $\cF(\{U_i \to U\})$. An arrow $\{\phi_i\}$ in $\cF\desc (\{U_i \to U\})$ is already an arrow in $\cF(\{U_i \to U\})$.

It is not hard to check that this functor is an equivalence of categories.

We can not define a functor $\cF(U) \arr \cF\desc (\{U_i \to U\})$ directly, without the choice of a cleavage. However, let us define another category
   \[
   \cF\comp (\{U_i \to U\}),
   \]
in which the objects are quadruples $(\xi, \{\xi_i\}, \{\xi_{ij}\}, \{\xi_{ijk}\}))$, where $\xi$ is an object of $\cF(U)$ and each $\xi_\alpha$ is an object of $\cF(U_{\alpha})$, plus a commutative cube
   \[
   \xymatrix@-15pt{
   &\xi_{ijk}\ar[rr]\ar'[d][dd]
   \ar[ld]
   &&\xi_{jk}\ar[ld]\ar[dd]\\
   \xi_{ij}\ar[dd]\ar[rr]
   && \xi_j\ar[dd]\\
   &\xi_{ik} \ar[ld]\ar'[r][rr]
   &&\xi_k\ar[ld]\\
   \xi_i\ar[rr]&&\xi
   }
   \]
in $\cF$ for all the triples of indices, in which all the arrows are cartesian, and whose image in $\cC$ is the cube (\ref{eq:basic-cube}) above. An arrow from $(\xi, \{\xi_i\}, \{\xi_{ij}\}, \{\xi_{ijk}\}))$ to $(\eta, \{\eta_i\}, \{\eta_{ij}\}, \{\eta_{ijk}\}))$ can be indifferently defined as an arrow $\phi \colon \xi \arr \eta$ in $\cF(U)$, or as collections of arrows $\xi \arr \eta$, $\xi_i \arr \eta_i$, $\xi_{ij} \arr \eta _{ij}$ and $\xi_{ijk} \arr \eta_{ijk}$ satisfying the obvious commutativity conditions.

There is a functor from $\cF\comp (\{U_i \to U\})$ to $\cF(U)$ that sends a whole object $(\xi, \{\xi_i\}, \{\xi_{ij}\}, \{\xi_{ijk}\}))$ to $\xi$, and is easily seen to be an equivalence. There is also a functor from $\cF\comp (\{U_i \to U\})$ to $\cF\desc (\{U_i \to U\})$ that forgets the object of $\cF(U)$. This takes the place of the functor from $\cF(U)$ to $\cF(\{U_i \to U\})$ defined using cleavages.

\begin{remark}
Of course, if one really wants to be consistent, one should not assume that the category $\cC$ has a canonical choice of fibered products, and not suppose that the $U_{ij}$ and the $U_{ijk}$ are given a priori, but allow them to be arbitrary fibered products.
\end{remark}

The most elegant definition of objects with descent data is one that uses sieves; it does not require choosing anything. Let $\cU = \{U_i \to U\}$ be a covering in $\cC$. The sieve  $\h_{\cU} \colon \cC\op \arr \catset$ is a functor, whose associated category fibered in sets is the full subcategory of $(\cC/U)$, whose objects are arrows $T \arr U$ that factor through some $U_i \arr U$. According to our principle that functors and categories fibered in sets should be identified, we denote by $\h_{\cU}$ this category. By the same principle, we also denote by $\h_{U}$ the category $(\cC/U)$.

There is a functor $\hom_\cC\bigl(\h_{\cU}, \cF\bigr) \arr \cF\desc(\cU)$, defined as follows. Suppose that we are given a morphism $F \colon \h_{\cU} \arr \catset$. For any triple of indices $i$, $j$ and $k$ we have objects $U_i \arr U$, $U_{ij} \arr U$ and $U_{ijk} \arr U$ of $\h_{\cU}$, and each of the projections of (\ref{eq:basic-cube}) not landing in $U$ is an arrow in $\h_{\cU}$. Hence we can apply $F$ and get a diagram
   \[
   \xymatrix@R-15pt@C-25pt{
   &F(U_{ijk})\ar[rr]\ar'[d][dd]
   \ar[ld]
   &&F(U_{jk})\ar[ld]\ar[dd]\\
   F(U_{ij})\ar[dd]\ar[rr]
   && F(U_j)\\
   &F(U_{ik}) \ar[ld]\ar[rr]
   &&F(U_k)\\
   F(U_i)&&
   }
   \]
giving an object of $\cF\desc (\cU)$. This extends to a functor
   \[
   \hom(\h_{\cU}, \cF) \arr \cF\desc(\cU)
   \]
in the obvious way.

Also, consider the functor 
   \[
   \hom_\cC(\h_{U},\cF)
   \arr
  \hom_\cC(\h_{\cU},\cF)
   \]
 induced by the embedding $\h_{\cU} \subseteq \h_{U}$; after choosing a cleavage, it is easy to verify that the composite of functors
   \[
   \hom_\cC(\h_{U},\cF)
   \arr
   \hom_\cC (\h_{\cU},\cF)
   \simeq
   \cF\desc(\cU)
   \simeq \cF(\cU)
   \]
is isomorphic to the composite
   \[
   \hom_\cC(\h_{U},\cF)
   \simeq
   \cF(U)
   \arr
   \cF(\cU)
   \]
where the first functor is the equivalence of the 2-Yoneda Lemma.

The following generalizes Proposition~\ref{prop:funnychar-gluing}.

\begin{proposition}\label{prop:funnychar-descent-data}
\index{object with descent data!via sieves}
\index{descent data!object with!via sieves}
The functor $\hom(\h_{\cU}, \cF) \arr \cF\desc(\cU)$ is an equivalence.
\end{proposition}

\begin{proof}
Let us construct a functor $\cF\desc(\cU) \arr \hom(\h_{\cU}, \cF)$. Set $\cU = \{U_{i} \to U\}$, and, for each $T \arr U$ in $\h_{\cU}T$, choose a factorization $T \arr U_{i} \arr U$. Assume that we given an object $(\{\xi_i\}, \{\xi_{ij}\},\{\xi_{ijk}\})$ of $\cF\desc(U)$; for each arrow $T \arr U$ in the category $\h_{\cU}$ we get an object $\xi_{T}$ of $T$ by pulling back $\xi_{i}$ along the chosen arrow $T \arr U_{i}$. This defines a function from the set of objects of $\h_{\cU}$ to $\cF$. 

Given an arrow $T'\arr T \arr U$ in $\h_{\cU}$, chose a factorization $T' \arr U_{j} \arr U$ of the composite $T' \arr U$. This, together with the composite $T' \arr T \arr U_{i}$ yields an arrow $T' \arr U_{ij}$ fitting into a diagram
   \[
   \xymatrix{
   T' \ar[r] \ar[d] & U_{ij} \ar[r]\ar[d] & U_{j} \ar[d] \\
   T  \ar[r]        & U_{i}  \ar[r]       & U\hsmash{.}
   }
   \]
Since the given arrow $\xi_{ij} \arr \xi_{j}$ is cartesian, the canonical arrow $\xi_{T'} \arr \xi_{j}$, that is given by definition, because $\xi_{T'}$ is a pullback of $\xi_{j}$, will factor uniquely as $\xi_{T'} \arr \xi_{ij} \arr \xi_{j}$, in such a way that $\xi_{T'} \arr \xi_{ij}$ maps to $T' \arr U_{ij}$. Now the composite $\xi_{T'} \arr \xi_{ij} \arr \xi_{i}$ will factor as $\xi_{T'} \arr \xi_{T} \arr \xi_{i}$ for a unique arrow $\xi_{T'} \arr \xi_{T}$ mapping to the given arrow $T' \arr T$ in $\cC$. According to Proposition~\refall{prop:comp-cartesian}{4}, the arrow $\xi_{T'} \arr \xi_{T}$ is cartesian. 

These two functions, on objects and on arrows, define a morphism $\h_{\cU} \arr \cF$. We need to check that the composites
   \[
   \cF\desc(\cU) \arr \hom(\h_{\cU}, \cF) \arr \cF\desc(\cU)
   \]
and
   \[
   \hom(\h_{\cU}, \cF) \arr \cF\desc(\cU) \arr \hom(\h_{\cU}, \cF)
   \]
are isomorphic to the identities. This is straightforward, and left to the reader.
\end{proof}

If we choose a cleavage, the composite of functors
   \[
   \hom_\cC(\h_{U},\cF)
   \arr
   \hom_\cC (\h_{\cU},\cF)
   \simeq
   \cF\desc(\cU)
   \simeq \cF(\cU)
   \]
is isomorphic to the composite
   \[
   \hom_\cC(\h_{U},\cF)
   \simeq
   \cF(U)
   \arr
   \cF(\cU)
   \]
where the first functor is the equivalence of the 2-Yoneda Lemma.

\subsection{Fibered categories with descent}
\label{sec:fibered-descent}

\begin{definition} Let $\cF \arr \cC$ be a fibered category on a site $\cC$.

\begin{enumeratei}

\item $\cF$ is a \emph{prestack}\index{prestack} over $\cC$ if for each covering $\{U_i \to U\}$ in $\cC$, the functor $\cF(U) \arr \cF(\{U_i \to U\})$ is fully faithful.

\item $\cF$ is a \emph{stack}\index{stack} over $\cC$ if for each covering $\{U_i \to U\}$ in $\cC$, the functor $\cF(U) \arr \cF(\{U_i \to U\})$ is an equivalence of categories.

\end{enumeratei}
\end{definition}

Concretely, for $\cF$ to be a prestack means the following. Let $U$ be an object of $\cC$, $\xi$ and $\eta$ objects of $\cF(U)$, $\{U_i \to U\}$ a covering, $\xi_i$ and $\eta_i$ pullbacks of $\xi$ and $\eta$ to $U_i$, $\xi_{ij}$ and $\eta_{ij}$ pullbacks of $\xi$ and $\eta$ to $U_{ij}$. Suppose that there are arrows $\alpha_i \colon \xi_i \arr \eta_i$ in $\cF(U_i)$, such that $\pr_1^* \alpha_i = \pr_2^* \alpha_j \colon \xi_{ij}  \arr \eta_{ij}$ for all $i$ and $j$. Then there is a unique arrow $\alpha \colon \xi \arr \eta$ in $\cF(U)$, whose pullback to $\xi_i \arr \eta_i$ is $\alpha_i$ for all $i$.

This condition can be restated using the functor of arrows of Section~\ref{sec:func-arrows}, and the comma topology on the category $(\cC/S)$(Definition~\ref{def:comma-topology}).

\begin{proposition}
\index{prestack!characterization of}
Let $\cF$ be a fibered category over a site $\cC$. Then $\cF$ is a prestack if and only if for any object $S$ of $\cC$ and any two objects $\xi$ and $\eta$ in $\cF(S)$, the functor $\underhom_S(\xi, \eta) \colon (\cC/S)\op \arr \catset$ is a sheaf in the comma topology.
\end{proposition}

\begin{proof}
Let us prove the first part. Assume that  for any object $S$ of $\cC$ and any two objects $\xi$ and $\eta$ in $\cF(S)$, the functor $\underhom_S(\xi, \eta) \colon (\cC/S)\op \arr \catset$ is a sheaf. Take an object $U$ of $\cC$, a covering $\{U_i\arr U\}$, and two objects $\xi$ and $\eta$ of $\cF(U)$. If we denote by $(\{\xi_i\}, (\alpha_{ij}))$ and $(\{\eta_i\}, (\beta_{ij}))$ the descent data associated with $\xi$ and $\eta$ respectively, we see easily that the arrows in $\cF(\{U_i \to U\})$ are the collections of arrows $\{\phi_i \colon \xi_i \arr \eta_i\}$ such that the restrictions of $\phi_i$ and $\phi_j$ to the pullbacks of $\xi$ and $\eta$ to $U_{ij}$ coincide. The fact that $\underhom_U(\xi, \eta )$ is a sheaf ensures that this comes from a unique arrow $\xi \arr \eta$ in $\cF(U)$; but this means precisely that the functor $\cF(U) \arr \cF(\{U_i \to U\})$ is fully faithful.

The proof of the opposite implication is similar, and left to the reader.
\end{proof}

\begin{definition}
An object with descent data $(\{\xi_i\}, \{\phi_{ij}\})$ in $\cF(\{U_i \to U\})$ is \emph{effective}\index{object with descent data!effective} if it is isomorphic to the image of an object of $\cF(U)$.
\end{definition}

Here is another way of saying this: an object with descent data $(\{\xi_i\}, \{\phi_{ij}\})$ in $\cF(\{U_i \to U\})$ is effective if there exists an object $\xi$ of $\cF(U)$, together with cartesian arrows $\xi_i \arr \xi$ over $\sigma_i \colon U_i \arr U$, such that the diagram
   \[
   \xymatrix@C-10pt@R-10pt{
   {}\pr_2^* \xi_j \ar[rr]^{\phi_{ij}} \ar[d] &&
   {}\pr_1^* \xi_i \ar[d] \\
   \xi_j \ar[rd] &&
   \xi_i \ar[ld]\\
   &\xi
   }
   \]
commutes for all $i$ and $j$. In fact, the cartesian arrows $\xi_i \arr \xi$ correspond to isomorphisms $\xi_i \simeq \sigma_i^* \xi$ in $\cF(U_i)$; and the commutativity of the diagram above  is easily seen to be equivalent to the cocycle condition.

Clearly, $\cF$ is a stack if and only if it is a prestack, and all objects with descent data in $\cF$ are effective.
\index{stack!characterization of}

Stacks are the correct generalization of sheaves, and give the right notion of ``sheaf of categories''. We should of course prove the following statement.

\begin{proposition}
Let $\cC$ be a site, $F \colon \cC\op \arr \catset$ a functor; we can also consider it as a category fibered in sets $F \arr \cC$.

\begin{enumeratei}

\item $F$ is a prestack if and only if it is a separated functor.

\item $F$ is stack if and only if it is a sheaf.

\end{enumeratei}
\end{proposition}

\begin{proof}
Consider a covering $\{U_i \to U\}$. The fiber of the category $F \arr \cC $ over $U$ is precisely the set $F(U)$, while the category $F(\{U_i \to U\})$ is the set of elements $(\xi_i) \in \prod_iF(U_i)$ such that the pullbacks of $\xi_i$ and  $\xi_j$ to $F(U_i \times_U U_j)$, via the first and second projections $U_i \times_U U_j \arr U_i$ and $U_i \times_U U_j \arr U_i$, coincide. The functor $F(U) \arr F(\{U_i \to U\})$ is the function that sends each element $\xi\in F(U)$ to the collection of restrictions $(\xi\mid_{U_i})$.

Now, to say that a function, thought of as a functor between discrete categories, is fully faithful is equivalent to saying that it is injective; while to say that it is an equivalence means that it is a bijection. From this both statements follow.
\end{proof}

\begin{remark}
The terminology here, due to Grothendieck, is a little unfortunate. Fibered categories are a generalization of functors: however, a presheaf is simply a functor, and thus, by analogy, a prestack should be simply a fibered category. What we call a prestack should be called a separated prestack.

I have decided to stick with Grothendieck's terminology, mostly because there is a notion of ``separated stack'' in the theory of algebraic stacks, and using the more rational term ``separated prestack'' would make ``separated stack'' pleonastic.
\end{remark}

\begin{example}
\index{fibered!category!of sheaves, as a stack}
\index{category!fibered!of sheaves, as a stack}
Let $\cC$ be a site. Then I claim that the fibered category $\catsh{\cC} \arr \cC$, defined in Example~\ref{ex:fibered-sheaves}, is a stack.

Here is a sketch of proof. Let $F$ and $G$ be two sheaves on an object $X$ of $\cC$: to show that $\catsh{\cC}$ is a prestack we want to show that $\underhom_{X}(F,G)\colon (\cC/X)\op \arr \catset$ is a sheaf.

For each arrow $U \arr X$, let us denote by $F_{U}$ and $G_{U}$ the restrictions of $F$ and $G$ to $U$. Let $\{U_{i} \arr U\}$ be a covering, $\phi_{i}\colon F_{U_{i}} \arr G_{U_{i}}$ a morphism of sheaves on $(\cC/U_{i})$, such that the restrictions of $\phi_{i}$ and $\phi_{j}$ to $(\cC/U_{ij})$ coincide. Denote by $\phi_{ij}\colon F_{U_{ij}} \arr G_{U_{ij}}$ this restriction. If $T \to U$ is an arrow, set $T_{i} = U_{i}\times_{U} T$, and consider the covering $\{T_{i} \to T\}$. Each $T_{i} \arr U$ factors through $U_{i}$, so $\phi_{i}$ defines a function $\phi_{i}\colon F T_{i} \arr GT_{i}$, and analogously $\phi_{ij}$ defines functions $\phi_{ij}\colon F T_{ij} \arr G T_{ij}$. There is commutative diagram of sets with rows that are equalizers
   \[
   \xymatrix{
   F T \ar[r] \ar@{-->}[d]
   & {}\prod_{i} F T_{i} \ar@<3pt>[r]\ar@<-3pt>[r] \ar[d]^{\prod_{i}\phi_{i}}
   & {}\prod_{ij} F T_{ij} \ar[d]^{\prod_{ij}\phi_{ij}}
   \\
   G T \ar[r] 
   & {}\prod_{i} G T_{i} \ar@<3pt>[r]\ar@<-3pt>[r]
   & {}\prod_{ij} G T_{ij}\hsmash{.}
   }
   \]
There is a unique function $\phi_{T}\colon F T \arr G T$ that one can insert in the diagram while keeping it commutative. This proves uniqueness. Also, it is easy to check that the collection of the $\phi_{T}$ defines a natural transformation $\phi_{U}\colon F_{U} \arr G_{U}$, whose restriction $F_{U_{i}} \arr G_{U_{i}}$ is $\phi_{i}$.

Now let us show that every object with descent data $(\{F_{i}\}, \{\phi_{ij}\})$ is effective. Here $F_{i}$ is a sheaf on $(\cC/U_{i})$, and $\phi_{ij}$ is an isomorphism of sheaves on $(\cC/U_{ij})$ between the restrictions $(F_{j})_{U_{ij}} \simeq (F_{i})_{U_{ij}}$.

For each object $T \arr U$ of $(\cC/U)$, set $T_{i} = U_{i}\times_{U} T$ as before, and define $F T$ to be the subset of $\prod_{i}F_{i} T_{i}$ consisting of objects $(s_{i}) \in \prod_{i} F T_{i}$, with the property that $\phi_{ij}$ carries the restriction $(s_{j})_{T_{ij}}$ to $(s_{i})_{T_{ij}}$. In other words, $F T$ is the equalizer of two functions $\prod_{i}F_{i} T_{i} \arr \prod_{ij}F_{i} T_{ij}$, where the first sends $(s_{i})$ to the collections of restrictions $\bigl((s_{i})_{T_{ij}}\bigr)$, and the second sends it to $\bigl(\phi_{ij}(s_{j})_{T_{ij}}\bigr)$.

For any arrow $T' \arr T$ in $(\cC/U)$, it is easy to see that the product of the restriction functions $\prod_{i}F_{i} T_{i} \arr \prod_{i}F_{i} T'_{i}$ carries $F T$ to $F T'$; this gives $F$ the structure of a functor $(\cC/U)\op \arr \catset$. We leave it to the reader to check that $F$ is a sheaf.

Now we have to show that the image of $F$ into $\catsh{\cC}(\{U_{i} \to U\})$ is isomorphic to $(\{F_{i}\}, \{\phi_{ij}\})$.

For each index $k$ let us construct an isomorphism of the restriction $F_{U_{k}}$ with $F_{k}$ as sheaves on $(\cC/U_{k})$. Let $T \to U_{k}$ be an object of $(\cC/U_{k})$, $s$ an element of $F_{k} T$. Each $T_{i}$ maps into $U_{ik}$, so we produce an element $\bigl(\phi_{ik}(s_{T_{i}})\bigr) \in \prod F_{i}T_{i}$; the cocycle condition ensures that this is an element of $F T$. This defines a natural transformation $F_{k} \arr F_{U_{k}}$.

In the other direction, let $T \to U_{k}$ be an object of $(\cC/U_{k})$, $(s_{i})$ an element of $F T \subseteq \prod_{i}F_{i} T_{i}$. Factor each $T_{i} \arr U_{i}$ through the projection $U_{ki} \arr U_{i}$. Then $\phi_{ki}(s_{i})$ is an element of $F_k T_{i}$, and the cocycle condition implies that the restrictions of $\phi_{ki}(s_{i})$ and $\phi_{kj}(s_{j})$ to $F_{k}T_{ij}$ coincide. Hence there a unique element of $F_{k} T$ that restricts to $\phi_{ki}(s_{i}) \in F_{k} T_{i}$ for each $i$. This construction defines a function $F T \arr F_{k} T$ for each $T$, which is easily seen to give a natural transformation $F_{U_{k}} \arr F_{k}$.

We leave it to the reader to check that these two natural transformations are inverse to each other, so they define an isomorphism of sheaves $F_{U_{k}} \simeq F_{k}$; and that this collection of isomorphisms constitutes an isomorphism in $\catsh{S}(\{U_{i} \arr U\})$ between the object associated with $F$ and the given object $(\{F_{i}\}, \{\phi_{ij}\})$, which is therefore effective.
\end{example}

\subsection{The functorial behavior of descent data}\label{subsec:func-descent-data}
\index{functoriality!of descent data}
\index{descent data!functoriality of}

Descent data have three kinds of functorial properties: they are functorial for morphisms of fibered categories, functorial on the objects, and functorial under refinement.

Let $F \colon \cF \arr \cG$ be a morphism of categories fibered over $\cC$. For any covering $\cU = \{U_{i} \to U\}$ we get a functor $F_{\cU} \colon \cF(\cU) \arr \cG(\cU)$ defined at the level of objects by the obvious rule
   \[
   F_{\cU}(\{\xi_{i}\}, \{\phi_{ij}\}) = (\{F\xi_{i}\}, \{F\phi_{ij}\})
   \]
and at the level of arrows by the equally obvious rule
   \[
   F_{\cU}\{\alpha_{i}\} = \{F\alpha_{i}\}.
   \]
Furthermore, if $\rho \colon F \arr G$ is a base-preserving natural transformation of morphisms, there is an induced natural transformation of functors $\rho_{\cU} \colon F_{\cU} \arr G_{\cU}$, defined by
   \[
   (\rho_{\cU})_{(\{\xi_{i}\}, \{\phi_{ij}\})} = \{\rho_{\xi_{i}}\}.
   \]
Therefore, if $F$ is an equivalence of fibered categories, $F_{\cU}$ is also an equivalence.

We leave it to the reader to check that the diagram 
   \[
   \xymatrix{
   {}\cF(U) \ar[r]\ar[d] & {}\cF(\cU)\ar[d] \\
   {}\cG(U) \ar[r]       & {}\cG(\cU)
   }
   \]
commutes, in the sense that the two composites $\cF(U) \arr \cG(\cU)$ are isomorphic. From this we obtain the following useful fact.

\begin{proposition}\label{prop:equivtostack->stack}\hfil

\begin{enumeratea}

\item If $F$ is an equivalence of fibered categories and $\cF(U) \arr \cF(\cU)$ is an equivalence of categories, then $\cG(U) \arr \cG(\cU)$ is also an equivalence of categories.

\item If two fibered categories over a site are equivalent, and one of them is a stack, or a prestack, the other is also a stack, or a prestack.

\end{enumeratea}
\end{proposition}

All this can be restated more elegantly using sieves. The morphism $F \colon \cF \arr \cG$ induces a functor $F_{*} \colon \hom(\h_{\cU}, \cF) \arr \hom(\h_{\cU}, \cG)$, that is the composite with $F$ at the level of objects. In this case the diagram becomes
   \[
   \xymatrix{
   {}\hom(\h_{U}, \cF) \ar[r]\ar[d] & {}\hom(\h_{\cU}, \cF)\ar[d] \\
   {}\hom(\h_{U}, \cG) \ar[r]       & {}\hom(\h_{\cU}, \cG)
   }
   \]
which strictly commutative, that is, the two composites are equal, not simply isomorphic. We leave the easy details to the reader.

We are not going to need the functoriality of descent data for the objects, so we will only sketch the idea: if $\{U_{i} \to U\}$ is a covering and $V \to U$ is an arrow, then there is a functor $\cF(\{U_{i}\to U\}) \arr 
\cF(\{V \times_{U} U_{i} \to V\})$. If $(\{\xi_{i}\}, \{\phi_{ij}\})$ is an object of $\cF(\{U_{i}\to U\})$, its image in $\cF(\{V \times_{U} U_{i} \to V\})$ is obtained by pulling back the $\xi_{i}$ and the $\phi_{ij}$ along the projection $V \times_{U} U_{i} \arr U_{i}$.

Now suppose that $\cF$ is a category fibered over $\cC$, $\cU = \{U_{i} \to U\}_{i \in I}$ a covering, $\cV = \{V_{i'} \to U\}_{i' \in I'}$ a refinement of $\cU$. For each index $i'$ choose a factorization $V_{i'} \xrightarrow{f_{i'}} U_{\mu(i')} \arr U$ for a certain $\mu(i') \in I$; this defines a function $\mu\colon I' \arr I$.

This induces a functor $\cF\cU \arr \cF\cV$, as follows. An object $(\{\xi_{i}\}, \{\phi_{ij}\})$ is sent to $(\{f_{i'}^{*}\xi_{\mu(i')}\}, \{f_{i'}^{*}\phi_{\mu(i')\mu(j')}\})$; we leave to the reader to check that this is also an object with descent data. An arrow 
   \[
   \{\alpha_{i}\}\colon (\{\xi_{i}\}, \{\phi_{ij}\}) \arr
                        (\{\eta_{i}\}, \{\psi_{ij}\})
   \]
is a collection of arrows $\alpha_{i}\colon \xi_{i}\arr \eta_{i}$ in $\cF(U_{i})$, and these can be pulled back to arrows $f_{i'}^{*}\alpha_{\mu(i')}\colon f_{i'}^{*}\xi_{\mu(i')} \arr f_{i'}^{*}\eta_{\mu(i')}$. We leave it the reader to verify that the collection $\{f_{i'}^{*}\alpha_{\mu(i')}\}$ yields an arrow
   \[
   (\{f_{i'}^{*}\xi_{\mu(i')}, \{f_{i'}^{*}\phi_{\mu(i')\mu(j')}\})\arr
   (\{f_{i'}^{*}\eta{\mu(i')}, \{f_{i'}^{*}\psi_{\mu(i')\mu(j')}\}),
   \]
and that this defines a functor.

This functor is essentially independent of the function $\mu\colon I' \arr I$, that is, if we change the function we get isomorphic functors. This is seen as follows. Suppose that $\nu\colon I' \arr I$ is another function, and that there are factorizations $V_{i'} \xrightarrow{g_{i'}} U_{\nu(i')} \arr U$. The two arrows $f_{i'}$ and $g_{i'}$ induce arrows $(f_{i'}, g_{i'})\colon U_{\mu(i') \nu(i')}$. We define an isomorphism $f_{i'}^{*}\xi_{\mu(i')} \simeq g_{i'}^{*}\xi_{\nu(i')}$ by composing the isomorphisms in the following diagram
   \[
   \xymatrix@C+40pt{
   f_{i'}^{*}\xi_{\mu(i')}\ar@{-->}[r]\ar[d]^{\simeq} &
      g_{i'}^{*}\xi_{\nu(i')}
      \\
   (f_{i'}, g_{i'})^{*}\pr_{1}^{*}\xi_{\mu(i')}
         \ar[r]^-{(f_{i'}, g_{i'})^{*}\phi_{ij}} & 
      (f_{i'}, g_{i'})^{*}\pr_{2}^{*}\xi_{\nu(i')}\ar[u]_{\simeq}
   }
   \]
The cocycle condition ensures that this gives an isomorphism of descent data between $(\{f_{i'}^{*}\xi_{\mu(i')}\}, \{f_{i'}^{*}\phi_{\mu(i')\mu(j')}\})$ and $(\{g_{i'}^{*}\xi_{\nu(i')}\}, \{g_{i'}^{*}\phi_{\nu(i')\nu(j')}\})$ (we leave the details to the reader). It is easy to check that this gives an isomorphism of functors.

Also, if $\cW = \{W_{i''} \to U\}_{i'' \in I''}$ is a refinement of $\cV$, it is also a refinement of $\cU$. After choosing functions $\mu\colon I' \arr I$, $\nu\colon I'' \arr I'$ and $\rho\colon I'' \arr I$, and factorizations $V_{i'} \xrightarrow{f_{i'}} U_{\mu(i')} \arr U$, $W_{i''} \xrightarrow{g_{i''}} U_{\nu(i'')} \arr U$ and $W_{i''} \xrightarrow{h_{i''}} U_{\rho(i'')} \arr U$ we get functors $\cF\cU \arr \cF\cV$, $\cF\cV \arr \cF\cW$ and $\cF\cU \arr \cF\cW$. I claim that the composite of the first two is isomorphic to the third.

To check this, we may change the factorizations $W_{i''} \xrightarrow{h_{i''}} U_{\rho(i'')} \arr U$, because, as we have just seen, this does not change the isomorphism class of the functor $\cF\cU \arr \cF\cW$; hence we may assume that $\rho = \mu \circ \nu \colon I'' \arr I$, and that $h_{i''}\colon W_{i''} \arr U_{\rho(i'')}$ equals the composite $f_{\nu(i'')}\circ g_{i''} \arr I$. Given an object $(\{\xi_{i}\}, \{\phi_{ij}\})$ of $\cF\cU$, its image in $\cF\cW$ under the functor $\cF\cU \arr \cF\cW$ is the object
   \[
   \bigl(\{(f_{\nu(i'')} g_{i''})^{*}\xi_{i''}\},
      \{(f_{\nu(i'')} g_{i''})^{*}
      \phi_{\nu\mu(i'')\,\nu\mu(j'')}\}\bigr)
   \]
of $\cF\cW$, while its image under the composite $\cF\cU \arr \cF\cV \arr \cF\cV \arr \cF\cW$ is the object
   \[
   \bigl(\{g_{i''}^{*}f_{\nu(i'')}^{*}\xi_{i''}\},
      \{g_{i''}^{*}f_{\nu(i'')}^{*}
      \phi_{\nu\mu(i'')\,\nu\mu(j'')}\}\bigr).
   \]
The canonical isomorphisms $(f_{\nu(i'')} g_{i''})^{*}\xi_{i''} \simeq g_{i''}^{*}f_{\nu(i'')}^{*}\xi_{i''}$ give an isomorphism of the two objects of $\cF\cW$, and this defines the desired isomorphism of functors.

Once again, in the language of sieves everything is much easier: if $\cV$ is a refinement of $\cU$, then $\h_{\cV}$ is a subfunctor of $\h_{\cU}$, and the embedding $\h_{\cV} \into \h_{\cU}$ induces a functor
   \[
   \hom_{\cC}(\h_{\cU}, \cF) \arr \hom_{\cC}(\h_{\cV}, \cF)
   \]
with no choice required. 

Also, in this language the composite 
   \[
   \hom_{\cC}(\h_{\cU}, \cF) \arr \hom_{\cC}(\h_{\cV}, \cF)
                             \arr \hom_{\cC}(\h_{\cW}, \cF)
   \]
equals the functor $\hom_{\cC}(\h_{\cU}, \cF) \arr \hom_{\cC}(\h_{\cW}, \cF)$ on the nose.

\subsection{Stacks and sieves}

Using the description of the category of objects with descent data in Proposition~\ref{prop:funnychar-descent-data} we can give the following very elegant characterization of stacks, which generalizes the characterization of sheaves given in Corollary~\ref{cor:funnychar-sheaves}.

\begin{corollary}\label{cor:funnychar-stacks}
A fibered category $\cF \arr \cC$ is a stack if and only if for any covering $\cU$ of an object $U$ of $\cC$ the functor
   \[
   \hom_\cC(\h_{U},\cF)
   \arr
   \hom_\cC(\h_{\cU},\cF)
   \]
induced by the embedding $\h_{\cU} \subseteq \h_{U}$ is an equivalence.
\end{corollary}

This can sharpened, as in Proposition~\ref{prop:funnychar-sheaves}.

\begin{proposition}\label{prop:funnychar-stacks}
\index{stack!characterization via sieves}
A fibered category $\cF \arr \cC$ is a stack if and only if for any object $\cU$ of $\cC$ and sieve $S$ on $U$ belonging to $\cT$, the functor
   \[
   \hom_\cC(\h_{U},\cF)
   \arr
   \hom_\cC(S,\cF)
   \]
induced by the embedding $\h_{U} \subseteq S$ is an equivalence.

Furthermore, $\cF$ is a prestack if and only if the functor above is fully faithful for all $U$ and $S$.
\end{proposition}

\begin{proof}
The fact that if the functor is an equivalence then $\cF$ is a stack follows from Corollary ~\ref{cor:funnychar-stacks}, so we only need to prove the converse (and similarly for the second statement).

Let $S$ be a sieve belonging to $\cT$ on an object $U$ of $\cC$. Choose a covering $\cU = \{U_{i} \to U\}$ of $U$ such that $\h_{\cU} \subseteq S$: the restriction functor $\hom_\cC(\h_{U},\cF) \arr \hom_\cC(\h_{\cU},\cF)$ is an equivalence, and it factors as
   \[
   \hom_\cC(\h_{U},\cF) \arr \hom_\cC(S,\cF) \arr \hom_\cC(\h_{\cU},\cF).
   \]

Again by Corollary ~\ref{cor:funnychar-stacks}, $\hom_\cC(\h_{U},\cF) \arr \hom_\cC(\h_{\cU},\cF)$ is fully faithful, and it is an equivalence when $\cF$ is a stack: hence $ \hom_\cC(S,\cF) \arr \hom_\cC(\h_{\cU},\cF)$ is full, and it essentially surjective whenever $\cF$ is a stack. So we see that the following lemma suffices.

\begin{lemma}\label{lem:prestack->faithful}
Let $\cF$ be a prestack over a site $\cC$, $S$ and $S'$ be sieves belonging to the topology of $\cC$ with $S' \subseteq S$. Then the induced restriction functor
   \[
   \hom_\cC(S,\cF) \arr \hom_\cC(S',\cF)
   \]
is faithful.
\end{lemma}

This is a generalization of Lemma~\ref{lem:separated->injective}.

\begin{proof}
The proof is very similar to that of Lemma~\ref{lem:separated->injective}. Let $F$ and $G$ be two morphisms $S \arr \cF$, $\phi$ and $\psi$ two base-preserving natural transformations, inducing the same natural transformations from the restriction of $F$ to $\h_{\cU}$ to that of $G$. Let $T \arr U$ be an arrow in $S$; we need to prove that
   \[
   \phi_{T \to U} = \psi_{T \to U} \colon F(T \to U) \arr G(T \to U).
   \]
Consider the fibered products $T \times_{U} U_{i}$, with the first projections $p_{i} \colon T \times_{U} U_{i}$. Since $F$ and $G$ are morphisms, and $S$ is a functor, so that every arrow in $S$ is cartesian, the arrows
   \[
   F(p_{i}) \colon F(T \times_{U} U_{i} \to U) \arr F(T \to U)
   \]
and
   \[
   G(p_{i}) \colon G(T \times_{U} U_{i} \to U) \arr G(T \to U)
   \]
are cartesian. Consider the covering $\{T \times_{U} U_{i} \to T\}$: since the composite $T \times_{U} U_{i} \to T \to U$ is in $\h_{\cU}T$, we have
   \[
   \phi_{T \times_{U} U_{i} \to U} = \psi_{T \times_{U} U_{i} \to U} \colon 
      F(T \times_{U} U_{i} \to U) \arr G(T \times_{U} U_{i} \to U).
   \]
Hence the commutativity of the diagrams
   \[
   \xymatrix{
   F(T \times_{U} U_{i} \to U) \ar[rr]^-{\phi_{T \times_{U} U_{i} \to U}}
         \ar[d]^{F(p_{i}}
      && G(T \times_{U} U_{i} \to U) \ar[d]^{G(p_{i})}\\
   F(T \to U) \ar[rr]^-{\phi_{T \to U}} &&
      G(T \to U)
   }
   \]
can be interpreted as saying that the pullbacks of $\phi_{T\to U}$ and $\psi_{T\to U}$ to $T \times_{U}U_{i}$ are the same. Since the functors of arrows of $\cF$ form a sheaf, since $\cF$ is a prestack, this implies that $\phi_{T\to U}$ and $\psi_{T\to U}$ are equal, as claimed.
\end{proof}

This ends the proof of Proposition~\ref{prop:funnychar-stacks}.
\end{proof}

Since two equivalent topologies on the same category have the same sieves, we obtain the following generalization of Proposition~\ref{prop:same-sheaves}.

\begin{proposition}\label{prop:same-stacks}
\index{topology!equivalent!has the same stacks}
\index{equivalent topologies!has the same stacks}
Let $\cC$ a category, $\cT$ and $\cT'$ two topologies on $\cC$, $\cF \arr \cC$ a fibered category. Suppose that $\cT'$ is subordinate to $\cT$. If $\cF$ is a prestack, or a stack, relative to $\cT$, then it is also a prestack, or a stack, relative to $\cT'$.
\end{proposition}

In particular, if $\cT$ and $\cT'$ are equivalent, then $\cF$ is a stack relative to $\cT$ if and only if it is also a stack relative to $\cT'$.

For later use, we note the following consequence of Lemma~\ref{lem:prestack->faithful}.

\begin{lemma}\label{lem:pass-to-refinement}
If $\cF$ is a prestack on a site, $\cU$ and $\cV$ two coverings of an object $U$ of $\cC$, such that $\cV$ is a refinement of $\cU$, and $\cF(U) \arr \cF(\cV)$ is an equivalence, then $\cF(U) \arr \cF(\cU)$ is also an equivalence.
\end{lemma}

\subsection{Substacks}

\begin{definition}
Let $\cC$ be a site, $\cF \arr \cC$ a stack. A \emph{substack}\index{substack} of $\cF$ is a fibered subcategory that is a stack.
\end{definition}

\begin{example}
Let $\cC$ be a site, $\cF \arr \cC$ a stack, $\cG$ a full subcategory of $\cF$ satisfying the following two conditions.

\begin{enumeratei}

\item Any cartesian arrow in $\cF$ whose target is in $\cG$  is also in $\cG$.

\item Let $\{U_i \to U\}$ be a covering in $\cC$, $\xi$ an object of $\cF(U)$, $\xi_i$ pullbacks of $\xi$ to $U_i$. If $\xi_i$ is in $\cG$ for all $i$, then $\xi$ is in $\cG$.

\end{enumeratei}

Then $\cG$ is a substack.
\end{example}

There are many examples of the situation above: for example, as we shall see (Theorem~\ref{thm:main}) the fibered category $\catqc{S}$ is a stack over $\catsch{S}$ with the fpqc topology. Then the full subcategory of $\catqc{S}$ consisting of locally free sheaves of finite rank satisfies the two conditions, hence it is a substack.

\begin{proposition}\call{prop:criterion-cart-stack}

Let $\cC$ be a site, $\cF \arr \cC$ a fibered category. Recall that $\cF\cart$ is the associated category fibered in groupoids (Definition~\ref{def:ass-groupoids}).

\begin{enumeratei}

\itemref{1} If $\cF$ is a stack, so is $\cF\cart$.

\itemref{2} If $\cF$ is a prestack and $\cF\cart$ is a stack, then $\cF$ is also a stack.

\end{enumeratei}
\end{proposition}

\begin{proof}
The isomorphisms in $\cF$ are all cartesian; hence, given a covering $\{U_i \to U\}$ in $\cC$, the categories $\cF(\{U_i \to U\})$ and $\cF\cart(\{U_i \to U\})$ have the same objects, and the effective objects with descent data are the same. So it is enough to prove that if $\cF$ is a prestack then $\cF\cart$ is a prestack.

Let $\xi$ and $\eta$ be two objects in some $\cF(U)$. Let $\{U_i \to U\}$ be a covering, $\xi_i$ and $\eta_i$ pullbacks of $\xi$ and $\eta$ to $U_i$, $\xi_{ij}$ and $\eta_{ij}$ pullbacks to $U_{ij}$, $\alpha_i \colon \xi_i \arr \eta_i$ arrows in $\cF\cart(U_i)$, such that $\pr_1^* \alpha_i = \pr_2^* \alpha_j \colon \xi_{ij} \arr \eta_{ij}$. Then there is unique arrow $\alpha \colon \xi \arr \eta$ that restricts to $\alpha_i$ for each $i$; and it is enough to show that $\alpha$ is cartesian. But the cartesian arrows in $\cF(U)$ and in each $\cF(U_i)$ are the isomorphisms; hence the $\alpha_i$ are isomorphisms, and the arrow $\alpha_i^{-1} \colon \eta_i \arr \xi$ comes from a unique arrow $\beta \colon \eta \arr \xi$. The composites $\beta \circ \alpha$ and $\alpha \circ \beta$ pull back to identities in each $\cF(U_i)$, and so they must be identities in $\cF(U)$. This shows that $\alpha$ is in $\cF\cart(U)$, and completes the proof.
\end{proof}





\section{Descent theory for quasi-coherent sheaves}

\subsection{Descent for modules over commutative rings}
\label{subsec:descent-rings}

Here we develop an affine version of the descent theory for \qc sheaves.  It is only needed to prove Theorem~\ref{thm:main} below, so it may be a good idea to postpone reading it until after reading the next section on descent for \qc sheaves.

If $A$ is a commutative ring, we will denote by $\catmod{A}$ the category of modules over $A$. 

Consider a ring homomorphism $f \colon A \arr B$. If $M$ is an $A$-module, we denote by $\iota_M \colon M \otimes_A B \simeq B \otimes_A M$ the usual isomorphism of $A$-modules defined by $\iota_M(n \otimes b) = b \otimes n$. Furthermore, we denote by $\alpha_M \colon M \arr B \otimes_{A} M$ the homomorphism defined by $\alpha_M(m) = 1 \otimes m$.

For each $r \ge 0$ set
   \[
   B^{\otimes r} = \overbrace{B \otimes_A B
   \otimes_A \dots \otimes _A B}^{\text{$r$ times}}.
   \]
A $B$-module $N$ becomes a module over $B^{\otimes 2}$ in two different ways, as $N \otimes_A B$ and $B \otimes_A N$; in both cases the multiplication is defined by the formula $(b_1 \otimes b_2)(x_1 \otimes x_2) = b_1x_1 \otimes b_2x_2$. Analogously, $N$ becomes a module over $B^{ \otimes 3}$ as $N \otimes_A B \otimes_A B$, $B \otimes_A N \otimes_A B$ and $B \otimes_A B \otimes_A N$ (more generally, $N$ becomes a module over $B^{\otimes r}$ in $r$ different ways; but we will not need this).

Let us assume that we have a homomorphism of $B^{\otimes 2}$-modules $\psi \colon N \otimes_A B \arr B \otimes_A N$. Then there are three associated homomorphism of $B^{\otimes 3}$-modules
   \begin{align*}
   \psi_1 &\colon B \otimes_A N \otimes_A B \arr
   B \otimes_A B \otimes_A N, \\
   \psi_2 &\colon N \otimes_A B \otimes_A B \arr
   B \otimes_A B \otimes_A N, \\
   \psi_3 &\colon N \otimes_A B \otimes_A B \arr
   B \otimes_A N \otimes_{A} B
   \end{align*}
by inserting the identity in the first, second and third position, respectively. More explicitly, we have $\psi_1 = \id_B \otimes \psi$, $\psi_3 = \psi \otimes \id_B$, while we have $\psi_2(x_1 \otimes x_2 \otimes x_3) = \sum_iy_i \otimes x_2 \otimes z_i$ if $\psi(x_1 \otimes x_3) = \sum_iy_i \otimes z_i \in B \otimes_A N$. Alternatively, $\psi_2 = (\id_B \otimes \iota_N) \circ (\psi\otimes\id_{B}) \circ (\id_N \otimes \iota_B)$.

Let us define a category $\catmod{A \arr B}$ as follows. Its objects are pairs $(N, \psi)$, where $N$ is a $B$-module and $\psi \colon N \otimes_A B \simeq B \otimes_A N$ is an isomorphism of $B^{\otimes 2}$-modules such that
   \[
   \psi_2 = \psi_1\circ \psi_3 \colon
   N \otimes_A B \otimes_A B \arr
   B \otimes_A B \otimes_A N.
   \]
An arrow $\beta \colon (N, \psi) \arr (N', \psi')$ is a homomorphism of $B$-modules $\beta \colon N \arr N'$, making the diagram
   \[
   \xymatrix{
   N \otimes_A B \ar[r]^{\psi} \ar[d]^{\beta \otimes \id_B} &
   B \otimes_A N \ar[d]^{\id_B \otimes \beta} \\
   N' \otimes_A B \ar[r]^{\psi'} &
   B \otimes_A N'
   }
   \]
commutative.

We have a functor $F \colon \catmod{A} \arr \catmod{A \arr B}$, sending an $A$-module $M$ to the pair $(B \otimes_A M, \psi_M)$, where
   \[
   \psi_M \colon (B \otimes_A M) \otimes_A B \arr
   B \otimes_A (B \otimes_A M)
   \]
is defined by the rule
   \[
   \psi_M(b \otimes m \otimes b') = b \otimes b' \otimes m.
   \]
In other words, $\psi_M = \id_B \otimes \iota_M$.

It is easily checked that $\psi_M$ is an isomorphism of $B^{\otimes 2}$-modules, and that $(M \otimes_A B, \psi_M)$ is in fact an object of $\catmod{A \arr B}$.

If $\alpha \colon M \arr M'$ is a homomorphism of $A$-modules, one sees immediately that $\id_B \otimes \alpha \colon B \otimes_A M
\arr B \otimes_A M'$ is an arrow in $\catmod{A \arr B}$. This defines the desired functor $F$.

\begin{theorem}\label{thm:main-affine}
\index{theorem!descent!for modules}
\index{descent!for modules}
If $B$ is faithfully flat over $A$, the functor
   \[
   F \colon \catmod{A} \arr \catmod{A \arr B}
   \]
defined above is an equivalence of categories.
\end{theorem}

\begin{proof}
Let us define a functor $G \colon \catmod{A \arr B} \arr \catmod{A}$. We send an object $(N, \psi)$ to the $A$-submodule $GN \subseteq N$ consisting of elements $n \in N$ such that $1 \otimes n = \psi(n
\otimes 1)$.

Given an arrow $\beta \colon (N, \psi) \arr (N', \psi')$ in $\catmod{A \arr B}$, it follows from the definition of an arrow that $\beta$ takes $GN$ to $GN'$; this defines the functor $G$.

We need to check that the composites $GF$ and $FG$ are isomorphic to the identity. For this we need the following generalization of Lemma~\ref{lem:exact-affine}. Recall that we have defined the two homomorphisms of $A$-algebras
   \[
   e_1, e_2 \colon B \arr B \otimes_A B
   \]
by
$e_1(b) = b \otimes 1$ and $e_2(b) = 1 \otimes b$.

\begin{lemma}\label{lem:exact2-affine}
Let $M$ be an $A$-module.
Then the sequence
   \[
   0 \arr M \overset{\alpha_M}
   \larr B \otimes_A  M
   \xrightarrow{(e_1-e_2) \otimes\id_M} B^{\otimes 2}\otimes M
   \]
is exact.
\end{lemma}

The proof is a simple variant of the proof of Lemma~\ref{lem:exact-affine}.

Now notice that
   \begin{align*}
   \bigl((e_1-e_2) \otimes\id_M\bigr)(b \otimes m)
   &{} = b \otimes 1 \otimes m - 1 \otimes b \otimes m\\
   &{} = \psi_M(b \otimes m \otimes 1) - 1 \otimes b \otimes m
   \end{align*}
for all $m$ and $b$; and this implies that 
   \[
   \bigl((e_1-e_2) \otimes\id_M\bigr)(x) = \psi_M(x \otimes 1) - 1 \otimes x
   \]
for all $x \in B \otimes_A M$. Hence $G(B \otimes_A M, \psi_M)$ is the kernel of $(e_1-e_2) \otimes\id_M$, and the homomorphism $M \arr B \otimes M$ establishes a natural isomorphism between $M$ and $G(M \otimes_A B) = GF(M)$, showing that $GF$ is isomorphic to the identity.

Now take an object $(N, \psi)$ of $\catmod{A \arr B}$, and set $M = G(N, \psi)$. The fact that $M$ is an $A$-submodule of the $B$-module $N$ induces a homomorphism of $B$-modules $\theta \colon B \otimes_A M \arr N$ with the usual rule $\theta(b \otimes m) = bm$.  Let us check that $\theta$ is an arrow in $\catmod{A \arr B}$, that is, that the diagram
   \[
   \xymatrix@C+10pt{
   B \otimes_A B \otimes_A M
      \ar[r]^-{\id_B \otimes \theta} \ar[d]^{\id_B \otimes \iota_M}&
   B \otimes_A N
      \ar[d]^\psi\\
   B \otimes_A M \otimes_A B 
      \ar[r]^-{\theta \otimes \id_B}&
   N \otimes_A B
   }
   \]
commutes. The calculation is as follows:
   \begin{alignat*}3
   \psi(\id_B \otimes \theta)(b \otimes b' \otimes m) &=
      \psi(b \otimes b'm)\\
   &= \psi\bigl((b \otimes b') (1 \otimes m)\bigr)\\
   &= (b \otimes b') \psi(1 \otimes m)\\
   &= (b \otimes b') (1 \otimes m)
   \qquad \text{(because $m \in M$)}\\
   &= (bm \otimes b')\\
   &= (\theta\otimes\id_B)(b \otimes m \otimes b')\\
   &= (\theta \otimes \id_B) (\id_B \otimes \iota_M)
      (b \otimes b' \otimes m).
   \end{alignat*}

So this $\theta$ defines a natural transformation $\id \arr FG$. We have to check that $\theta$ is an isomorphism.

Consider the homomorphisms $\alpha, \beta\colon N \arr B \otimes N$ defined by $\alpha (n) = 1 \otimes n$ and $\beta(n) = \psi(n \otimes 1)$; by definition, $M$ is the kernel of $\alpha - \beta$. There is a diagram with exact rows
   \[
   \xymatrix@C+20pt{0 \ar[r] &
   M \otimes B
   \ar[r]^{i \otimes \id_B} \ar[d]^{\theta \circ \iota_M} &
   N \otimes_A B \ar[d]^\psi
   \ar[r]^-{(\alpha - \beta)\otimes_A\id_B} &
   B \otimes_A N \otimes_A B \ar[d]^{\psi_1} \\
   0 \ar[r] &
   N \ar[r]^{\alpha_M} &
   B \otimes_A N \ar[r]^-{(e_2 - e_1)\otimes\id_N}&
   B \otimes_A B \otimes_A N \\
   }
   \]
where $i \colon M \into N$ denotes the inclusion. Let us show that it is commutative. For the first square, we have
   \[
   \alpha_M \theta \iota_M(m \otimes b) = 1 \otimes bm
   \]
while
   \begin{align*}
   \psi (i \otimes \id_B)(m \otimes b) &= \psi(m \otimes b) \\
   &{}= \psi\bigl((1 \otimes b)(m \otimes 1)\bigr) \\
   & = (1 \otimes b)\psi(m \otimes 1) \\
   & = (1 \otimes b) (1 \otimes m) \\
   & = 1 \otimes bm.
   \end{align*}
For the second square, it is immediate to check that $\psi_1 \circ (\alpha \otimes\id_B) = (e_2\otimes\id_N) \circ \psi$. On the other hand
   \begin{align*}
   \psi_1 (\beta \otimes \id_B)(n \otimes b)
   & = \psi_1\bigl(\psi(n \otimes 1) \otimes b\bigr) \\
   & = \psi_1 \psi_3 (n \otimes 1 \otimes b) \\
   & = \psi_2(n \otimes 1 \otimes b) \\
   & = (e_1 \otimes\id_N)\psi(n \otimes b).
   \end{align*}
   Both $\psi$ and $\psi_1$ are isomorphisms; hence $\theta \circ \iota_N$ is an isomorphism, so $\theta$ is an isomorphism, as desired.

This finishes the proof of Theorem~\ref{thm:main-affine}.
\end{proof}

\subsection{Descent for \qc sheaves}

Here is the main result of descent theory for \qc sheaves. It states that \qc sheaves satisfy descent with respect to the fpqc topology; in other words, they form a stack with respect to either topology. This is quite remarkable, because \qc sheaves are sheaves in that Zariski topology, which is much coarser, so a priori one would not expect this to happen.

Given a scheme $S$, recall that in \S\ref{subsec:fibered-quasi-coherent} we have constructed the fibered category $\catqc{S}$ of \qc sheaves, whose fiber of a scheme $U$ over $S$ is the category $\qcoh U$ of \qc sheaves on $U$.

\begin{theorem}\label{thm:main}
\index{theorem!descent!for quasicoherent sheaves}
\index{descent!for quasicoherent sheaves}
Let $S$ be a scheme. The fibered category $\catqc{S}$ over $\catsch{S}$ is stack with respect to the fpqc topology.
\end{theorem}

\begin{remark}\label{rmk:really-really-need-finiteness}
\index{wild flat topology}\index{topology!wild flat}
This would fail in the ``wild'' flat topology of Remark~\ref{rmk:really-need-finiteness}: in this topology $\catqc{S}$ is not even a prestack.

Take the covering $\{V_p \arr U\}$ defined there, and the \qc sheaf $\oplus_p \cO_p$, the direct sum of the structure sheaves of all the closed points. The restriction of $\oplus_p \cO_p$ to each $V_p$ is the structure sheaf $\cO_p$ of the closed point, since pullbacks commute with direct sums, and the restriction of each $\cO_q$ to $V_p$ is zero for $p \neq q$. For each $p$ consider the projection $\pi_p \colon \cO_{V_p} \arr \cO_p$; it easy to see that $\pr_1^* \pi_p = \pr_2^* \pi_q \colon \cO_{V_p \times_U V_q} \arr f_{p,q}^* (\oplus_{p} \cO_p)$, where $f_{p,q} \colon V_p \times_U V_q \arr U$ is the obvious morphism.

On the other hand there is no homomorphism $\cO_U \arr \oplus_{p} \cO_p$ that pulls back to $\pi_p$ for each $p$. In fact, such a homomorphism would correspond to a section of $\oplus_{p} \cO_p$ that is $1$ at each closed point, and this does not exist because of the definition of direct sum.
\end{remark}

For the proof of the theorem we will use the following criterion, a generalization of that of Lemma~\ref{lem:criterion-sheaf}.

\begin{lemma}\label{lem:criterion-stack}
\index{stack!characterization of fpqc}
\index{fpqc!characterization of fpqc stacks}

Let $S$ be a scheme, $\cF$ be a fibered category over the category $\catsch{S}$. Suppose that the following conditions are satisfied.

\begin{enumeratei}

\item $\cF$ is a stack with respect to the Zariski topology.

\item Whenever $V \arr U$ is a flat surjective morphism of affine $S$-schemes, the functor
   \[
   \cF(U) \arr \cF(V \to U)
   \]
is an equivalence of categories.

\end{enumeratei}

Then $\cF$ is a stack with respect to the the fpqc topology.
\end{lemma}

\begin{proof} The proof is a little long and complicates, so for clarity we will divide it in several steps. According to Theorem~\ref{thm:equivalent-split} and Proposition~\ref{prop:equivtostack->stack} we may assume that $\cF$ is split: this will only be used in the last two steps.

\steps

\step[: $\cF$ is a prestack] Given an $S$-scheme $T \arr S$ and two objects $\xi$ and $\eta$ of $\cF(T)$, consider the functor
   \[
   \underhom_T(\xi, \eta) \colon \catsch{T}\op
   \arr \catset.
   \]
We see immediately that the two conditions of Lemma~\ref{lem:criterion-sheaf} are satisfied, so the functor $\underhom_T(\xi, \eta)$ is a sheaf, and $\cF$ is a prestack in the fpqc topology.

Now we have to check that every object with descent data is effective.

\step[: reduction to the case of a single morphism]

We start by analyzing the sections of $\cF$ over the empty scheme $\emptyset$.

\begin{lemma}
The category $\cF(\emptyset)$ is equivalent to a category with one object and one morphism.
\end{lemma}

Equivalently, between any two objects of $\cF(\emptyset)$ there is a unique arrow.

\begin{proof}
The scheme $\emptyset$ has the empty Zariski covering $\cU = \emptyset$. By this I really mean the empty set, consisting of no morphisms at all, and not the set consisting of the embedding of $\emptyset \subseteq \emptyset$. There is only one object with descent data $(\emptyset, \emptyset)$ in $\cF(\cU)$, and one morphism from $(\emptyset, \emptyset)$ to itself. Hence $\cF(\cU)$ is equivalent to the category with one object and one morphism; but $\cF(\emptyset)$ is equivalent to $\cF(\cU)$, because $\cF$ is a stack in the Zariski topology.
\end{proof}

\begin{lemma}\label{lem:disjoint-product}
If a scheme $U$ is a disjoint union of open subschemes $\{U_i\}_{i \in I}$, then the functor $\cF(U) \arr \prod_{i \in I} \cF(U_i)$ obtained from the various restriction functors $\cF(U) \arr \cF(U_i)$ is an equivalence of categories.
\end{lemma}

\begin{proof}
Let $\xi$ and $\eta$ be objects of $\cF(U)$; denote by $\xi_i$ and $\eta_i$ their restrictions to $U_i$. The fact that $\underhom_U(\xi, \eta) \colon \catsch{T}\op \arr \catset$ is a sheaf ensures that the function
   \[
   \hom_{\cF(U)}(\xi, \eta) \arr
   \prod_i \hom_{\cF(U_i)}(\xi_i, \eta_i)
   \]
is a bijection; but this means precisely that the functor is fully faithful.

To check that it is essentially surjective, take an object $(\xi_i)$ in $\prod_{i \in I}\cF(U_i)$. We have $U_{ij} = \emptyset$ when $i \neq j$, and $U_{ij} = U_i$ when $i = j$; we can define transition isomorphisms $\phi_{ij} \colon \pr_2^* \xi_j \simeq \pr_1^* \xi_i$ as the identity when $i = j$, and as the only arrow from $\pr_2^* \xi_j$ to $\pr_1^* \xi_i$ in $\cF(U_{ij}) = \cF(\emptyset)$  when $i \neq j$. These satisfy the cocycle condition; hence there is an object $\xi$ of $\cF(U)$ whose restriction to each $U_i$ is isomorphic to $\xi_i$. Then the image of $\xi$ into $\prod_{i \in I} \cF(U_i)$ is isomorphic to $(\xi_i)$, and the functor is essentially surjective.
\end{proof}

Given an arbitrary covering $\{U_i \to U\}$, set $V = \coprod_i U_i$, and denote by $f \colon V \arr U$ the induced morphism. I claim that the functor $\cF(U) \arr \cF(V \to U)$ is an equivalence if and only if $\cF(U) \arr \cF(\{U_i \to U\})$ is. In fact, we will show that there is an equivalence of categories
   \[
   \cF(V \to U) \arr \cF(\{U_i \to U\})
   \]
such that the composite
   \[
   \cF(U) \arr \cF(V \to U)
   \arr \cF(\{U_i \to U\})
   \]
is isomorphic to the functor
   \[
    \cF(U) \arr \cF(\{U_i \to U\}).
   \]
This is obtained as follows. We have a natural isomorphism of $U$-schemes
   \[
   V \times_U V \simeq \coprod_{i,j} U_i \times_U U_j,
   \]
so Lemma~\ref{lem:disjoint-product} gives us equivalences of categories
   \begin{equation}\label{eq:first-equivalence}
   \cF(V) \arr \prod_i \cF(U_i)
   \end{equation}
and
   \begin{equation}\label{eq:second-equivalence}
    \cF(V \times_U V) \arr
   \prod_{i,j}\cF(U_i \times_U U_j).
   \end{equation}

An object of $\cF(V \to U)$ is a pair $(\eta, \phi)$, where $\eta$ is an object of $\cF(V)$ and $\phi \colon \pr_2^* \eta \simeq \pr_1^* \eta$ in $\cF(V \times_U V)$ satisfying the cocycle condition. If $\eta_i$ denotes the restriction of $\eta$ to $U_i$ for all $i$ and $\phi_{ij} \colon \pr_2^* \eta \simeq \pr_1^* \eta$ the arrow pulled back from $\phi$, the image of $\phi$ in $\prod_{i,j}\cF(U_i \times_U U_j)$ is precisely the collection $(\phi_{ij})$; it is immediate to see that $(\phi_{ij})$ satisfies the cocycle condition. 

In this way we associate with each object  $(\eta, \phi)$ of $\cF(V \to U)$ an object $(\{\eta_i\}, \{\phi_{ij}\})$ of $\cF(\{U_i \to U\})$. An arrow $\alpha \colon (\eta, \phi) \arr (\eta', \phi')$ is an arrow $\alpha \colon \eta \arr \eta'$ in $\cF(V)$ such that
   \[
   \pr_1^* \alpha \circ \phi = \phi \circ \pr_2^* \alpha \colon 
   \pr_2^* \eta \arr \pr_1^* \eta';
   \]
then one checks immediately that the collection of restrictions $\{\alpha_i \colon \eta_i \arr \eta'_i\}$ gives an arrow $\{\alpha_i\} \colon (\{\eta_i\}, \{\phi_{ij}\}) \arr (\{\eta'_i\}, \{\phi'_{ij}\})$.

Conversely, one can use the inverses of the functors (\ref{eq:first-equivalence}) and (\ref{eq:second-equivalence}) to define the inverse of the functor constructed above, thus showing that it is an equivalence (we leave the details to the reader). This equivalence has the desired properties.

This means that to check that descent data in $\cF$ are effective we can restrict consideration to coverings consisting of one arrow.

\step[: the case of a quasi-compact morphism with affine target]

Consider the case that $V \arr U$ is a flat surjective morphism of $S$-schemes, with $U$ affine and $V$ quasi-compact. Let $\{V_{i}\}$ be a finite covering of $V$ by open affine subschemes, $V'$ the disjoint union of the $V_{i}$. Then $\cF(U) \arr \cF(V' \to U)$ is an equivalence, by hypothesis, so by Lemma~\ref{lem:pass-to-refinement} $\cF(U) \arr \cF(V \to U)$ is also an equivalence.

\step[: the case of a morphism with affine target]

Now $U$ is affine and $V \arr U$ is an arbitrary fpqc morphism. By hypothesis, there is an open covering $\{V_{i}\}$ of $V$ by quasi-compact open subschemes, all of which surject onto $U$. We will use the fact that $\cF$ has a splitting. We need to show that $\cF(U) \arr \cF(V \to U)$ is essentially surjective.

Choose an index $i$; $V_{i} \arr U$ is also an fpqc cover, with $V_{i}$ quasi-compact. We have a strictly commutative diagram of functors
   \[
   \xymatrix@C-30pt{
   {}\cF(U) \ar[rr] \ar[rd] && {}\cF(V \to U) \ar[ld] \\
   &{}\cF(V_{i} \to U)
   }
   \]
in which $\cF(U) \arr \cF(V_{i} \to U)$ is an equivalence, because of the previous step, and $\cF(U) \arr \cF(V \to U)$ is fully faithful. From this we see that to show that $\cF(U) \arr \cF(V \to U)$ is essentially surjective it is enough to prove that $\cF(V \to U) \to \cF(V_{i} \to U)$ is fully faithful.

It is clear that $\cF(V \to U) \to \cF(V_{i} \to U)$ is full (because $\cF(U) \arr \cF(V_{i} \to U)$ is), so it is enough to show that it is faithful. Let $\alpha$ and $\beta$ be two arrows in $\cF(V \to U)$ with the same image in $\cF(V_{i} \to U)$. For any other index $j$ we have that the restriction functor $\cF(V_{i} \cup V_{j} \to U) \to \cF(V_{i} \to U)$ is an equivalence, because in the strictly commutative diagram
   \[
   \xymatrix@C-30pt{
   {}\cF(U) \ar[rr] \ar[rd] && {}\cF(V_{i}\cup V_{j} \to U) \ar[ld] \\
   &{}\cF(V_{i} \to U)
   }
   \]
the top and left arrows are equivalences, so the restrictions of $\alpha$ and $\beta$ to $V_{i} \cup V_{j}$ are the same. Hence the restrictions of $\alpha$ and $\beta$ to each $V_{j}$ are the same: since $\cF$ is a prestack we can conclude that $\alpha = \beta$.

\step[: the general case] Now consider a general fpqc morphism $f\colon V \arr U$, with no restrictions on $U$ or $V$. Take an open covering $\{U_{i}\}$ of $U$ by affine subschemes, and let $V_{i}$ be the inverse image of $U_{i}$ in $V$. We need to show that any object $(\eta, \phi)$ of $\cF(V \to U)$ comes from an object of $\cF(U)$. 

For each open subset $U' \subseteq U$, denote by $\Phi_{U'}\colon \cF(U') \arr \cF(f^{-1}U' \to U')$ the functor that sends objects to objects with descent data, defined via the splitting. We will use the obvious fact that if $U''$ is an open subset of $U'$, the diagram
   \[
   \xymatrix{
   {}\cF(f^{-1}U') \ar[r]\ar[d]^{\Phi_{U'}}    
      & {}\cF(f^{-1}U'')\ar[d]^{\Phi_{U''}} \\
   {}\cF(f^{-1}U'\to U') \ar[r] & {}\cF(f^{-1}U''\to U'')
   },
   \]
where the rows are given by restrictions, is strictly commutative.

For each index $i$, let $(\eta_{i}, \phi_{i})$ be the restriction of $(\eta, \phi)$ to $V_{i}$. This is an object of $\cF(V_{i} \to U_{i})$, hence by the previous step it there exists an object $\xi_{i}$ of $\cF(U_{i})$ with an isomorphism $\alpha_{i}\colon \Phi_{U_{i}}\xi_{i} \simeq (\eta_{i}, \phi_{i})$ in $\cF(\{V_{i} \to U_{i}\})$. Now we want to glue together the $\xi_{i}$ to a global object $\xi$ of $\cF(U)$; for this we need Zariski descent data $\phi_{ij}\colon \xi_{j}\mid_{U_{ij}} \simeq \xi_{i}\mid_{U_{ij}}$.

For each pair of indices $i$ and $j$, set $V_{ij} = V_{i} \cap V_{j}$, so that $V_{ij}$ is the inverse image of $U_{ij}$ in $V$. By restricting the isomorphisms $\alpha_{i}\colon \Phi_{U_{i}}\xi_{i} \simeq (\eta_{i}, \phi_{i})$ to $V_{ij}$ we get isomorphisms
   \[
   \Phi_{U_{ij}}(\xi_{i}\mid_{U_{ij}}) = (\Phi_{U_{i}}\xi_{i})\mid_{V_{ij}}
      \stackrel{\alpha_{i}\mid_{V_{ij}}}{\simeq}
      (\eta\mid_{V_{ij}}, \phi\mid_{V_{ij}})
   \]
and from these isomorphisms
   \[
   \alpha_{i}^{-1}\alpha_{j} \colon \Phi_{U_{ij}}(\xi_{j}\mid_{U_{ij}})
      \simeq \Phi_{U_{ij}}(\xi_{i}\mid_{U_{ij}}).
   \]
Since the functor $\Phi_{U_{ij}}$ is an equivalence of categories, there exists a unique isomorphism $\phi_{ij}\colon \xi_{j}\mid_{V_{ij}} \simeq \xi_{i}\mid_{V_{ij}}$ such that $\Phi_{U_{ij}}\phi_{ij} = \alpha_{i}^{-1}\alpha_{j}$.

By applying $\Phi_{U_{ijk}}$ we see easily that the cocycle condition $\phi_{ik} = \phi_{ij}\phi_{jk}$ is satisfied; hence there exists an object $\xi$ of $\cF(U)$, with isomorphisms $\xi\mid_{U_{i}} \simeq \xi_{i}$. If we denote by $f_{i}\colon V_{i} \arr U$ the restriction of $\colon V \arr U$, we obtain an isomorphism of $f_{i}^{*}\xi = f^{*}\xi\mid_{V_{i}}$ with the pullback of $\xi_{i}$ to $V_{i}$. We also have an isomorphism of the pullback of $\xi_{i}$ to $V_{i}$ with $\eta_{i}$, obtained by pulling back $\alpha_{i}$ along $V_{i}\arr U_{i}$. These isomorphisms $f^{*}\xi\mid_{V_{i}} \simeq \eta_{i}$ coincide when pulled back to $V_{ij}$, so they glue together to give an isomorphism $f^{*}\xi \simeq \eta$; and this is the desired isomorphism of $\Phi_{U}\xi$  with $(\eta, \phi)$.

This completes the proof of Lemma~\ref{lem:criterion-stack}.
\end{proof}

It is a standard fact that $\catqc S$ is a stack in the Zariski topology; so we only need to check that the second condition of Lemma~\ref{lem:criterion-stack} is satisfied; for this, we use the theory of \S\ref{subsec:descent-rings}. Take a flat surjective morphism $V \arr U$, corresponding to a faithfully flat ring homomorphism $f \colon A \arr B$. We have the standard equivalence of categories $\qcoh{U} \simeq \catmod{A}$; I claim that there is also an equivalence of categories $\qcoh{V \arr U} \simeq \catmod{A \arr B}$. A \qc sheaf $\cM$ on $U$ corresponds to an $A$-module $M$. The inverse images $\pr_1^* \cM$ and  $\pr_2^* \cM$ in $V \times_U B = \spec B \otimes_A B$ correspond to the modules $M \otimes_A B$ and $B \otimes_A M$, respectively; hence an isomorphism $\phi\colon \pr_2^* \cM \simeq \pr_1^* \cM$ corresponds to an isomorphism $\psi \colon M \otimes_A B \simeq B \otimes_A M$. It is easy to see that $\phi$ satisfies the cocycle condition, so that $(\cM, \phi)$ is an object of $\qcoh{V\arr U}$, if and only if $\psi$ satisfies the condition $\psi_1\psi_3 = \psi_2$; this gives us the equivalence $\qcoh{V \arr U} \simeq \catmod{A \arr B}$. The functor $\qcoh U \arr \qcoh{V \arr U}$ corresponds to the functor $\catmod{A} \arr \catmod{A \arr B}$ defined in \S\ref{subsec:descent-rings}, in the sense that the composites
   \[
   \qcoh U \arr \qcoh{V \arr U} \simeq \catmod{A \arr B}
   \]
and
   \[
   \qcoh U \simeq \catmod{A} \simeq \catmod{A \arr B}
   \]
are isomorphic. Since $\catmod{A} \arr \catmod{A \arr B}$ is an equivalence, this finishes the proof of Theorem~\ref{thm:main}.

Here is an interesting question. Let us call a morphism of schemes $V \arr U$ a  \emph{descent morphism}\index{descent morphism} if the functor $\qcoh{U} \arr \qcoh{V \to U}$ is an equivalence.

Suppose that a morphism of schemes $V \arr U$ has local sections in the fpqc topology, that is, there exists an fpqc covering $\{U_{i} \to U\}$ with sections $U_{i} \to V$. This is equivalent to saying that $V \to U$ is a covering in the saturation of the fpqc topology, so, by Theorem~\ref{thm:main} and Proposition~\ref{prop:same-stacks}, it is a descent morphism.

\begin{namedr}{Open question}
Do all descent morphisms have local sections in the fpqc topology? If not, is there an interesting characterization of descent morphisms?
\end{namedr}

\subsection{Descent for sheaves of commutative algebras}

There are many variants of Theorem~\ref{thm:main}. The general principle is that one has descent in the fpqc topology for \qc sheaves with an additional structure, as long as this structure is defined by homomorphisms of sheaves, satisfying conditions that are expressed by the commutativity of certain diagrams.

Here is a typical example. If $U$ is a scheme, we may consider \qc sheaves of commutative algebras on $U$, that is, sheaves of commutative $\cO_U$-algebras that are \qc as sheaves of $\cO_U$-modules. The \qc sheaves of commutative algebras on a scheme $U$ form a category, denoted by $\qcohalg U$.

We get a pseudo-functor on the category $\catsch{S}$ by sending each $U \arr S$ to the category $\qcohalg U$; we denote the resulting fibered category on $\catsch{S}$ by $\catqcalg{S}$.

\begin{theorem}\label{thm:algebras-stack}
\index{theorem!descent!for sheaves of commutative quasi-coherent algebras}
\index{descent!for sheaves of commutative quasi-coherent algebras}
$\catqcalg{S}$ is a stack over $\catsch{S}$.
\end{theorem}

Here is the key fact.

\begin{lemma} \call{lem:local-algebrahom}
Let $\{\sigma_i \colon U_i \arr U\}$ be an fpqc covering of schemes.
\begin{enumeratei}

\itemref{1} If $\cA$ and $\cB$ are \qc sheaves of algebras over $U$, $\phi \colon \cA \arr \cB$ is a homomorphism of \qc sheaves, such that each pullback $\sigma_i^* \phi \colon \sigma_i^* \cA \arr \sigma_i^* \cB$ is a homomorphism of algebras for all $i$, then $\phi$ is a homomorphism of algebras.

\itemref{2} Let $\cA$ be a \qc sheaf on $U$. Assume that each pullback $\sigma_i^* \cA$ has a structure of sheaf of commutative algebras, and that the canonical isomorphism of \qc sheaves $\pr_2^* \sigma_j^* \cA \simeq \pr_1^* \sigma_i^* \cA$ is an isomorphism of sheaves of algebras for each $i$ and $j$. Then there exists a unique structure of sheaf of commutative algebras on $\cA$ inducing the given structure on each $\sigma_i^*\cA$.

\end{enumeratei}
\end{lemma}

\begin{proof}
For part~\refpart{lem:local-algebrahom}{1}, we need to check that the two composites
   \[
   \cA\otimes_{\cO_U} \cA
   \overset{\mu_{\cA}}\larr
   \cA
   \arr
   \cB
   \]
and
   \[
   \cA\otimes_{\cO_U} \cA
   \xrightarrow{\phi \otimes \phi}
   \cB\otimes_{\cO_U} \cB
   \overset{\mu_{\cB}}\larr
   \cB
   \]
coincide. However, the composites
   \[
   \sigma_i^*\cA \otimes_{\cO_{U_i}}
   \sigma_i^*\cA
   \simeq
   \sigma_i^*(\cA\otimes_{\cO_U} \cA)
   \xrightarrow{\sigma_i^* \mu_{\cA}}
   \sigma_i^*\cA
   \xrightarrow{\sigma_i^* \phi}
   \sigma_i^*\cB
   \]
and 
   \[
   \sigma_i^*\cA \otimes_{\cO_{U_i}} \sigma_i^*\cA
   \simeq
   \sigma_i^*(\cA\otimes_{\cO_U} \cA)
   \xrightarrow{\sigma_i^* \phi\otimes\phi}
   \sigma_i^*(\cB\otimes_{\cO_U} \cB)
   \xrightarrow{\sigma_i^* \mu_{\cB}}
   \sigma_i^*\cB   
   \]
coincide, because $\sigma_i^*\phi$ is a homomorphism of sheaves of algebras; and we know that two homomorphisms of \qc sheaves on a scheme that are locally equal in the fpqc topology, are in fact equal, because $\catqc{S}$ is a stack over $\catsch{S}$.

Let us prove part~\refpart{lem:local-algebrahom}{2}. From the algebra structure on each $\sigma_i^*\cA$ we get homomorphisms of \qc sheaves
   \[
   \mu_i \colon \sigma_i^*(\cA \otimes_{\cO_U} \cA)
   \simeq
   \sigma_i^* \cA \otimes_{\cO_{U_i}} \sigma_i^*\cA
   \arr \sigma_i^*\cA.
   \]
We need to show that these homomorphisms are pulled back from a homomorphism $\cA \otimes_{\cO_U} \cA \arr \cA$. Denote by $\sigma_{ij} \colon U_{ij} \arr U$ the obvious morphism; since $\catqc{S}$ is a stack, and in particular the functors of arrows are sheaves, this is equivalent to proving that the pullbacks $\sigma_{ij}^*(\cA \otimes_{\cO_{U}} \cA) \arr \cA$ of $\mu_i$ and $\mu_j$ coincide for all $i$ and $j$. This is most easily checked at the level of sections.

Similarly, the homomorphisms $\cO_{U_{i}} \arr \cA_{i}$ corresponding to the identity come from a homomorphism $\cO_{U} \arr \cA$.

We also need to show that the resulting homomorphism $\cA \otimes_{\cO_U} \cA \arr \cA$ gives $\cA$ the structure of a sheaf of commutative algebras, that is, we need to prove that the product is associative and commutative, and that the homomorphism $\cO_{U} \arr \cA$ gives an identity. Once again, this is easily done by looking at sections, and is left to the reader.
\end{proof}

From this it is easy to deduce Theorem~\ref{thm:algebras-stack}. If $\{\sigma_{i}\colon U_{i} \arr U\}$ is an fpqc covering, and  $(\{\cA_{i}\}, \{\phi_{ij}\})$ is an object of $\catqcalg{S}(\{\sigma_{i}\colon U_{i} \to U\})$, we can forget the algebra structure on the $\cA_{i}$, and simply consider it as an object of $\catqc{S}(\{\sigma_{i}\colon U_{i} \to U\})$; then it will come from a \qc sheaf $\cA$ on $U$. However, each $\sigma_{i}^{*}\cA$ is isomorphic to $\cA_{i}$, thus it inherits a commutative algebra structure: and Lemma~\ref{lem:local-algebrahom} implies that this comes from a structure of sheaf of commutative algebras on $\cA$. This finishes the proof of Theorem~\ref{thm:algebras-stack}.

Exactly in the same way one can defined fibered categories of sheaves of (not necessarily commutative) associative algebras, sheaves of Lie algebras, and so on, and prove that all these structures give stacks.

\section{Descent for morphisms of schemes}
\label{sec:descent-morphisms}

Consider a site $\cC$, a stable class $\cP$ of arrows, and the associate fibered category  $\cP \arr \cC$, as in Example~\ref{ex:restrictedarrows}.

The following fact is often useful.

\begin{proposition}\label{prop:stable->prestack}
Let $\cC$ be a subcanonical site, $\cP$ a stable class of arrows. Then $\cP \arr \cC$ is a prestack.
\end{proposition}

Recall (Definition~\ref{def:subcanonical}) that a site is subcanonical when every representable functor is a sheaf. The site $\catsch{S}$ with the fpqc topology is subcanonical (Theorem~\ref{thm:rep-fppf}).

\begin{proof}
Let $\{U_i \to U\}$ be a covering, $X \arr U$ and $Y \arr U$ two arrows in $\cP$. The arrows in $\cP(U)$ are the arrows $X \arr Y$ in $\cC$ that commute with the projections to $U$. Set $X_i = U_i \times_U X$ and $X_{ij} = U_{ij}\times_U X = X_i \times_X X_j$, and analogously for $Y_i$ and $Y_{ij}$. Suppose that we have arrows $f_i \colon X_i \arr Y_i$ in $\cP(U_i)$, such that the arrows $X_{ij} \arr Y_{ij}$ induced by $f_i$ and $f_j$ coincide; we need to show that there is a unique arrow $f \colon X \arr Y$ in $\cP(U)$ whose restriction $X_i \arr Y_i$ coincides with $f_i$ for each $i$.

The composites $X_i \xrightarrow{f_i} Y_i \arr  Y$ give sections $g_i \in \h_Y(X_i)$, such that the pullbacks of $g_i$ and $g_j$ to $X_{ij}$ coincide. Since $\h_Y$ is a sheaf, $\{X_i \arr X\}$ is a covering, and $X_{ij} = X_i \times_X X_j$ for any $i$ and $j$, there is a unique arrow $f \colon X \arr Y$ in $\cC$, such the composite $X_i \arr X \xrightarrow{f} Y$ is $g_i$, so that the diagram
   \[
   \xymatrix{
   X_i \ar[r]^{f_i} \ar[d] &
   Y_i \ar[d] \\
   X \ar[r]^f &
   Y
   }
   \]
commutes for all $i$. It is also clear that the arrows $X \arr U$ and $X \xrightarrow{f} Y \arr U$ coincide, since they coincide when composed with $U_i \arr U$ for all $i$, and since $\h_U$ is a sheaf, and in particular a  separated functor. Hence the diagram
   \[
   \xymatrix{
   X \ar[r]^f \ar[d] & Y \ar[d] \\
   U \ar@{=}[r]      & U
   }
   \]
commutes, and $f$ is the only arrow in $\cP(U)$ whose restriction to each $U_i$ coincides with $f_i$.
\end{proof}

However, in general $\cP$ will not be a stack. It is easy to see that $\cP$ cannot be a stack unless it satisfies the following condition.

\begin{definition}\label{def:local-class}
A class of arrows $\cP$ in $\cC$ is \emph{local}\index{class of arrows!local}\index{local class of arrows} if it is stable (Definition~\ref{def:stable-class}), and the following condition holds. Suppose that you are given a covering $\{U_i \to U\}$ in $\cC$ and an arrow $X \arr U$. Then, if the projections $U_i \times_U X \arr U_i$ are in $\cP$ for all $i$, $X \arr U$ is also in $\cP$.
\end{definition}

Still, a local class of arrow does not form a stack in general, effectiveness of descent data is not guaranteed, not even when $\cP$ is the class of all arrows. Consider the following example. Take $\cC$ to be the class of all schemes locally of finite type over a field $k$, of bounded dimension, with the arrows being morphisms of schemes over $k$. Let us equip it with the Zariski topology, and let $\cP$ be the class of all arrows. Call $U = U_1 \coprod U_2 \coprod U_3 \coprod U_4 \coprod \ldots$ the union of countably many copies of $U_1 = \spec k$. The collection of inclusions $\{U_i \into U\}$ forms a covering. Over each $U_i$ consider the scheme $\AA^i_k$. Obviously $U_{ij} = \emptyset$ if $i \neq j$, and $U_{ij} = \spec k$ if $i = j$, so we define transition isomorphisms in the only possible way as the identity $\phi_{ii} = \id_{\AA^i_k} \colon \AA^i_k \arr \AA^i_k$, and as the identity $\id_{\emptyset} \colon \emptyset \arr \emptyset$ when $i \neq j$. These obviously satisfy the cocycle condition, being all identities. On the other hand there cannot be a scheme of bounded dimension over $U$, whose pullback to each $U_i$ is $\AA^i_{k}$.

This is an artificial example; obviously if we want to glue together infinitely many algebraic varieties, we shouldn't ask for the dimension to be bounded. And in fact, morphisms of schemes form a stack in the Zariski topology, and therefore a local category of arrows also forms a stack in the Zariski topology. 

On the other hand, most of the interesting properties of morphisms of schemes are local in the fpqc topology on the codomain, such as for example being flat, being of finite presentation, being quasi-compact, being proper, being smooth, being affine, and so on (Proposition~\ref{prop:local-fpqc}). For each of these properties we get a prestack of morphisms of schemes over $\catsch{S}$, and we can ask if this is a stack in the fpqc topology. 

The issue of effectiveness of descent data is rather delicate, however. We will give an example to show that it can fail even for proper and smooth morphisms, in the \'etale topology (see~\ref{subsec:failure}). In this section we will prove some positive results.

\subsection{Descent for affine morphisms}
\label{subsec:affine-descent}

Let $\cP$ be the class of affine arrows in $\catsch{S}$, and denote by $\cataff{S} \arr \catsch{S}$ the resulting fibered category. The objects of $\cataff{S}$ are affine morphisms $X \arr U$, where $U \arr S$ is an $S$-scheme.

\begin{theorem}\label{thm:affine-stack}
\index{theorem!descent!for affine morphisms}
\index{descent!for affine morphisms}
The fibered category $\cataff{S}$ is a stack over $\catsch{S}$ in the fpqc topology.
\end{theorem}

First of all, $\cataff{S}$ is a prestack, because of Proposition~\ref{prop:stable->prestack}, so the only issue is effectiveness of descent data. By Proposition~\refall{prop:criterion-cart-stack}{2} it is enough to check that $\cataff{S}\cart$ is a stack.

Let $\cA$ be a \qc sheaf of algebras on a scheme $U$. Then we denote by $\cursspec_X \cA$ the relative spectrum of $\cA$; this is an affine scheme over $U$, and if $V \subseteq U$ is an open affine subscheme of $U$, the inverse image of $V$ in $\cursspec_U \cA$ is the spectrum of the ring $\cA(V)$.

A homomorphism of sheaves of commutative rings $\cA \arr \cB$ induces a homomorphism of $U$-schemes $\cursspec_U \cB \arr \cursspec_U \cA$; this is a contravariant functor from $\qcohalg U$ to the category $\aff U$ of affine schemes over $U$, which is well-known to be an equivalence of categories $\qcohalg{U}\op \simeq \aff{U}$. The inverse functor sends an affine morphism $h \colon X \arr U$ to the \qc sheaf of commutative algebras $h_* \cO_X$.

There is a morphism of fibered categories
   \[
   \catqcalg{S}\cart \arr \cataff{S}\cart
   \]
that sends a an object $\cA$ of $\qcohalg U$ to the affine morphism $\cursspec \cA \arr U$. Let $(f, \alpha) \colon (U, \cA) \arr (V, \cB)$ be an arrow in $\catqcalg{S}\cart$; $f \colon U \arr V$ is a morphism of $S$-schemes, $\alpha \colon \cA \simeq f^{*}\cB$ an isomorphism of sheaves of $\cO_{U}$-modules. Then $\alpha^{-1}$ gives an isomorphism $\cursspec_U \cA \simeq \cursspec_U f^* \cB = U \times_V \cursspec_V \cB$ of schemes over $U$, and the composite of this isomorphism with the projection $U \times_V \cursspec_V \cB \arr \cursspec_V \cB$ gives an arrow from $\cursspec_U \cA \arr U$ to $\cursspec_V \cB \arr V$ in $\cataff S\cart$.

If we restrict the morphism to a functor
   \[
   \catqcalg{S}\cart (U) \arr \cataff{S}\cart(U)
   \]
for some $S$-scheme $U$ we obtain an equivalence of categories; hence this morphism is an equivalence of fibered categories over $\catsch{S}$, by Proposition~\ref{prop:char-equivalence}. Since $\catqc S\cart$ is a stack, by Theorem~\ref{thm:algebras-stack} and Proposition~\refall{prop:criterion-cart-stack}{1}, we see from Proposition~\ref{prop:equivtostack->stack} that $\cataff{S}\cart$ is also a stack, and this concludes the proof of Theorem~\ref{thm:affine-stack}.

The following corollary will be used in \S\ref{subsec:descent-ample}.

\begin{corollary}\label{cor:descent-embeddings}
Let $P \arr U$ be a morphism of schemes, $\{U_{i} \to U\}$ an fpqc cover. For each $i$ set $P_{i} \eqdef U_{i}\times_{U} P$ and $P_{ij} \eqdef U_{ij}\times_{U} P$. Suppose that for each $i$ we have a closed subscheme $X_{i}$ of $P_{i}$, with the property that for each pair of indices $i$ and $j$ the inverse images of $X_{i}$ and $X_{j}$ in $P_{ij}$, through the first and second projection respectively, coincide. Then there is a unique closed subscheme of $P$ whose inverse image in $P_{i}$ coincides with $X_{i}$ for each $i$.
\end{corollary}

\begin{proof}
We have that $\{P_{i} \to P\}$ is an fpqc cover, and $P_{ij} = P_{i}\times_{P}P_{j}$. The pullbacks $\pr_{2}^{*}X_{j}$ and $\pr_{1}^{*}X_{i}$ to $P_{ij}$ coincide as subschemes of $P_{ij}$, and this yields a canonical isomorphism $\phi_{ij} \pr_{2}^{*}X_{j} \simeq \pr_{1}^{*}X_{i}$. The cocycle condition is automatically satisfied, because any two morphisms of $P_{ijk}$-schemes that are embedded in $P_{ijk}$ automatically coincide. Hence there is an affine morphism $X \arr P$ that pulls back to $X_{i} \to P_{i}$ for each $i$; and this morphism is a closed embedding, because of Proposition~\ref{prop:local-fpqc}.

Uniqueness is clear, because two closed subschemes of $P$ that are isomorphic as $P$-schemes are in fact equal.
\end{proof}

\subsection{The base change theorem}\label{subsec:base-change}

For the next result we are going to need a particular case of the base change theorem for \qc sheaves.

Suppose that we have a commutative diagram of schemes
   \begin{equation}\label{eq:cart-diagram}
   \xymatrix{
   X \ar[r]^f \ar[d]^\xi & Y\ar[d]^\eta\\
   U \ar[r]^\phi         & V
   }
   \end{equation}
and a sheaf of $\cO_Y$-modules $\cL$. Then there exists a natural base change homomorphism of $\cO_U$-modules
   \[
   \beta_{f, \phi}(\cL) \colon \phi^*\eta_*\cL
   \arr \xi_* f ^* \cL
   \]
that is defined as follows. First of all, start from the natural adjunction homomorphism $\cL \arr f_*f^* \cL$ (this is the homomorphism that corresponds to $\id_{f^*\cL}$ in the natural adjunction isomorphism $\hom_Y(\cL, f_*\phi^* \cL) \simeq \hom_X (f^*\cL, f^*\cL)$). This gives a homomorphism of $\cO_V$-modules
   \[
   \eta_*\cL \arr \eta_*f_*f^* \cL =
   \phi_*\xi_* f^* \cL.
   \]
Then $\beta_{f, \phi}(\cL)$ corresponds to this homomorphism under the adjunction isomorphism\index{base-change!map}\index{base-change!homomorphism}
   \[
   \hom_U(\phi^*\eta_* \cL, \xi_*f^* \cL) \simeq
   \hom_V(\eta_* \cL, \phi_*\xi_*f^* \cL).
   \]

The homomorphism $\beta_{f, \phi}(\cL)$ has the following useful characterization at the level of sections. If $V_{1}$ is an open subset of $V$, and $s \in \cL(\eta^{-1}V_{1}) = \eta_{*}\cL(V_{1})$, then there is a pullback section $\phi^{*}s \in \phi^{*}\eta_{*}\cL(\phi^{-1}V_{1})$. The sections of this form generate $\phi^{*}\eta_{*}\cL$ as an $\cO_{U}$-module. The section
   \[
   f^{*}s \in \cL(f^{-1}\eta^{-1}V_{1}) = \cL(\xi^{-1}\phi^{-1}s)
   \]
can be considered as an element of $\xi_{*}f^{*}\cL(\phi^{-1}V_{1})$; and then $\beta_{f,\phi}(\cL)$ is characterized as the only $\cO_{U}$-linear homomorphism of sheaves such that
   \[
   \beta_{f,\phi}(\cL)(\phi^{*}s) = f^{*}s \in
      \xi_{*}f^{*}\cL(\phi^{-1}V_{1})
   \]
for all $s$ as above.

The base change homomorphism is functorial in $\cL$. That is, there are two functors $\phi^*\eta_*$ and $\xi_* f ^*$ from $\catqc Y$ to $\catqc U$, and $\beta_{f, \phi}$ gives a natural transformation $\phi^*\eta_* \arr \xi_* f^*$. 

The base change homomorphism also satisfies a compatibility condition.

\begin{proposition}\label{prop:compatibility-base-change}
\index{base-change!homomorphism!compatibility}
Let
   \[
   \xymatrix{
   X \ar[r]^f \ar[d]^\xi & Y \ar[r]^g \ar[d]^\eta & Z \ar[d]^\zeta \\
   U \ar[r]^\phi         & V \ar[r]^\psi           & W
   }
   \]
be a commutative diagram of schemes, $\cL$ a sheaf of $\cO_{Z}$-modules. Then the diagram of $\cO_U$-modules
   \[
   \xymatrix{
   {}\phi^*\psi^*\zeta_*\cL
   \ar[rr]^-{\isoass_{\phi, \psi}(\zeta_*\cL)} 
   \ar[d]_{\phi^* \beta_{g, \psi}(\cL)}&&
   {} (\psi\phi)^*\zeta_* \cL
   \ar[d]^{\beta_{gf, \psi\phi}(\cL)}\\
   {}\phi^*\eta_*g^* \cL \ar[rd]_{\beta_{f, \phi}(g^*\cL)}&&
   {}\xi_*(gf)^* \cL\ar[ld]^{\xi_*\isoass_{f,g}(\cL)}\\
   & {}\xi_*f^*g^* \cL
   }
   \]
commutes.
\end{proposition}

\begin{proof}
This is immediately proved by taking an open subset $W_{1}$ of $W$, a section $s \in \zeta_{*}\cL(W_{1}) = \cL(\zeta^{-1}W_{1})$, and following $\phi^{*}\psi^{*}s$ in the diagram above.
\end{proof}

Since in Proposition~\ref{prop:compatibility-base-change} the homomorphisms $\isoass_{\phi, \psi}(\zeta_*\cL)$ and $\xi_*\isoass_{f,g}(\cL)$ are always isomorphism, we get the following corollary.

\begin{corollary}\label{cor:basechange->basechange}
In the situation of Proposition~\ref{prop:compatibility-base-change}, assume that the base change homomorphism $\beta_{g, \psi}(\cL)$ is an isomorphism. Then $\beta_{gf, \psi\phi}(\cL)$ is an isomorphism if and only if $\beta_{f, \phi}(g^*\cL)$ is an isomorphism
\end{corollary}

Here is the base change theorem, in the form in which we are going to need it. This is completely standard in the noetherian case; the proof reduces to this case with reduction techniques that are also standard.

\begin{proposition}
\index{base-change!theorem}
\index{theorem!base-change}
Suppose that the diagram (\ref{eq:cart-diagram}) is cartesian, that $\eta$ is proper and of finite presentation, and that $\cL$ is \qc and of finite presentation, and flat over $V$. For any point $v \in V$ denote by $Y_v$ the fiber of $\eta$ over $v$, and by $\cL_v$ the restriction of $\cL$ to $Y_v$.

If $\H^1(Y_v, \cL_v) = 0$ for all $v \in V$, then $\eta_*\cL$ is locally free over $V$, and the base change homomorphism $\beta_{f, \phi}(\cL) \colon \phi^*\eta_*\cL \arr \xi_* f ^* \cL$ is an isomorphism.
\end{proposition}

\begin{proof}
If $V$ is noetherian, then the result follows from \cite[7.7]{ega3-2} (see also \cite[III 12]{hartshorne}).

In the general case, the base change homomorphism is easily seen to localize in the Zariski topology on $U$; hence we may assume that $V$ is affine. Set $V = \spec A$. According to \cite[Proposition 8.9.1, Th\'eor\`eme 8.10.5 and Th\'eor\`eme~11.2.6]{ega4-3} there exists a subring $A_0 \subseteq A$ that is of finite type over $\ZZ$, hence noetherian, a scheme $Y'_0$ that is proper over $\spec A_0$, and a coherent sheaf $\cL'_0$ on $Y'_0$ that is flat over $A_0$, together with an isomorphism of $(Y, \cL)$ with the pullback of $(Y'_0, \cL'_0)$ to $\spec A$.

By semicontinuity (\cite[Th\'eor\`eme~7.6.9]{ega3-2}), the set of points $v_0 \in \spec A_0$ such that the restriction of $\cL'_0$ to the fiber of $Y'_0$ over $v_0$ has nontrivial $\H^1$ is closed in $\spec A_0$; obviously, it does not contain the image of $\spec A$. Denote by $V_0$ the open subscheme that is the complement of this closed subset, by $Y_0$ and $\cL_0$ the restrictions of $Y'_0$ and $\cL'_0$ to $V_0$. Then $\spec A$ maps into $V_0$, and $(Y, \cL)$ is isomorphic to the pullback of $(Y_0, \cL_0)$ to $\spec A$; hence the result follows from Corollary~\ref{cor:basechange->basechange} and from the noetherian case.
\end{proof}

\subsection{Descent via ample invertible sheaves}\label{subsec:descent-ample}

Descent for affine morphism can be very useful, but is obviously limited in scope. One is more easily interested in projective morphisms, rather than in affine ones. Descent works in this case, as long as the projective morphisms are equipped with ample invertible sheaves, and these also come with descent data.

\begin{theorem}\label{thm:descent-ample}
\index{theorem!descent!via ample invertible sheaves}
\index{descent!via ample invertible sheaves}
Let $S$ be a scheme, $\cF$ be a class of flat proper morphisms of finite presentation in $\catsch S$ that is local in the fpqc topology (Definition~\ref{def:local-class}). Suppose that for each object $\xi\colon X \arr U$ of $\cF$ one has given an invertible sheaf $\cL_\xi$ on $X$ that is ample relative to the morphism $X \arr U$, and for each cartesian diagram
   \[
   \xymatrix{
   X \ar[r]^f \ar[d]^\xi & Y\ar[d]^\eta\\
   U \ar[r]^\phi         & V
   }
   \]
an isomorphism $\rho_{f, \phi} \colon f^*\cL_\eta \simeq \cL_\xi$ of invertible  sheaves on $X$. These isomorphisms are required to satisfy the following condition: whenever we have a cartesian diagram of schemes
   \[
   \xymatrix{
   X \ar[r]^f \ar[d]^\xi & Y \ar[r]^g \ar[d]^\eta & Z \ar[d]^\zeta \\
   U \ar[r]^\phi         & V \ar[r]^\psi           & W
   }
   \]
whose columns are in $\cF$, then the diagram
   \[
   \xymatrix@C+15pt{
   f^*g^*\cL_\zeta \ar[r]^-{\isoass_{f,g}(\cL_\zeta)}
      \ar[d]^{f^* \rho_{g, \psi}}
      & (gf)^*\cL_\zeta \ar[d]^{\rho_{gf, \psi\phi}}
   \\
   f^* \cL_\eta  \ar[r]^-{\rho_{f, \phi}}&
      {}\cL_\xi 
   }
   \]
of \qc sheaves on $X$ commutes. Here $\isoass_{f,g}(\cL_\zeta)$ is the canonical isomorphism of \S\ref{subsec:fibered-quasi-coherent}.

Then $\cF$ is a stack in the fpqc topology.
\end{theorem}

Another less cumbersome way to state the compatibility condition is using the formalism of fibered categories, which will be freely used in the proof: since $(gf)^{*}\cL_{\zeta}$ is a pullback of $g^{*}\cL_{\zeta}$ to $Y$, we can consider the pullback $f^{*}\rho_{g, \psi}\colon (gf)^{*}\cL_{\zeta} \arr f^{*}\cL_{\eta}$, and then the condition is simply the equality
   \[
   \rho_{gf, \psi\phi} = \rho_{f, \phi} \circ f^* \rho_{g, \psi}
      \colon (gf)^*\cL_\zeta \arr \cL_{\xi}.
   \]

\begin{example}\label{ex:stacks-curves}
For any fixed base scheme $S$ and any non-negative integer $g$ we can consider the class $\cF_{g, S}$ of proper smooth morphisms, whose geometric fibers are connected curves of genus $g$. These morphisms form a local class in $\catsch{S}$.

If $g \neq 1$ then the theorem applies. For $g \ge 2$ we can take $\cL_{X \arr U}$ to be the relative cotangent sheaf $\Omega^{1}_{X/U}$, or one of its powers, while for $g = 0$ we can take its dual. So $\cF_{g, S}$ is a stack. The stack $(\cF_{g, S})\cart$ (Definition~\ref{def:ass-groupoids}) is usually denoted by $\cM_{g, S}$, and plays an important role in algebraic geometry.

There is no natural ample sheaf on families of curves of genus $1$, so this theorem does not apply. In fact, $\cF_{1, S}$, as we have defined it here, is not a stack: this follows from the counterexample of Raynaud in \cite[XIII~3.2]{raynaudample}.

See Remark~\ref{rmk:algebraic-spaces} for further discussion.
\end{example}

\begin{proof}[Proof of Theorem~\ref{thm:descent-ample}]

The fact that $\cF$ is a prestack in the fpqc topology follows from Proposition~\ref{prop:stable->prestack}. It is also easy to check that $\cF$ is a stack in the Zariski topology.

For each object $\xi \colon X \arr U$ of $\cF$ we define a \qc finitely presented sheaf $\cM_{\xi}$ on $U$ as $\cM_{\xi} \eqdef \xi_{*}\cL_{X}$. Given a cartesian square
   \[
   \xymatrix{
   X \ar[r]^f \ar[d]^\xi & Y\ar[d]^\eta\\
   U \ar[r]^\phi         & V
   }
   \]
we get a homomorphism of \qc sheaves
   \[
   \sigma_{f,\phi}\colon \phi^{*}\cM_{\eta} = \phi^{*}\eta_{*}\cL_{\eta}
   \arr \xi_{*}\cL_{X} = \cM_{\xi}
   \]
by composing the base change homomorphism
   \[
   \beta_{f,\phi} \colon \phi^{*}\eta_{*}\cL_{\eta} \arr 
      \xi_{*}f^{*}\cL_{\eta}
   \]
with the isomorphism
   \[
   \xi_{*}\rho_{f,\phi} \colon \xi_{*}f^{*}\cL_{\eta}
      \simeq \xi_{*}\cL_{\xi}.
   \]
This homomorphism satisfies the following compatibility condition.

\begin{proposition}\label{prop:compatibility-sigma}
Given a cartesian diagram of schemes
   \[
   \xymatrix{
   X \ar[r]^f \ar[d]^\xi & Y \ar[r]^g \ar[d]^\eta & Z \ar[d]^\zeta \\
   U \ar[r]^\phi         & V \ar[r]^\psi           & W
   }
   \]
whose columns are in $\cF$, the composite
   \[
   (\psi\phi)^{*}\cM_{\zeta} \xarr{\phi^{*}\sigma_{g,\psi}}
   \phi^{*}\cM_{\eta} \xarr{\sigma_{f,\phi}}
   \cM_{\xi}
   \]
equals $\sigma_{gf,\psi\phi}$.
\end{proposition}

\begin{proof}
Take a section $s$ of $\cM_{\zeta}$ on some open subset of $W$, and let us check that 
   \[
   \sigma_{f,\phi} \circ \phi^{*}\sigma_{g,\psi} \bigl((\psi\phi)^{*}s\bigr) =
      \sigma_{gf,\psi\phi}\bigl((\psi\phi)^{*}s\bigr)
   \]
(since the sections of the form $(gf)^{*}s$ generate $(\psi\phi)^{*}\cM_{\zeta}$ as a sheaf of $\cO_{U}$-modules, this is enough). The section $s$ is a section of $\cL_{\zeta}$ on some open subset of $Z$, and we have
\begin{align*}
   \phi^{*}\sigma_{g,\psi} \bigl((\psi\phi)^{*}s\bigr) &=
   \phi^{*}(\eta_{*}\rho_{g,\psi} \circ
      \beta_{g,\psi})\bigl((\psi\phi)^{*}s\bigr)\\
   &= \phi^{*}\bigl(\rho_{g,\psi} \circ 
      \beta_{g,\psi}(\psi^{*}s)\bigr)\\
   &= \phi^{*}\bigl(\rho_{g,\psi} (g^{*}s)\bigr);
\end{align*}
hence
   \begin{align*}
   \sigma_{f,\phi} \circ \phi^{*}\sigma_{g,\psi} \bigl((\psi\phi)^{*}s\bigr)
   &= (\xi_{*}\rho_{f,\phi} \circ \beta_{f,\phi})
      \phi^{*}\bigl(\rho_{g,\psi} (g^{*}s)\bigr)\\
   &= \rho_{f,\phi}\bigl(f^{*}\rho_{g,\psi} (g^{*}s)\bigr)\\
   &= \phi_{gf, \psi\phi} \bigl((gf)^{*}s\bigr)\\
   &= \xi_{*}\phi_{gf, \psi\phi} \circ \beta_{gf, \psi\phi} 
      \bigl((\psi\phi)^{*}s\bigr)\\
   &= \sigma_{gf,\psi\phi}\bigl((\psi\phi)^{*}s\bigr).\qedhere
   \end{align*}
\end{proof}

We use once again the criterion of Lemma~\ref{lem:criterion-stack}. Consider a flat surjective morphism $V \arr U$ of affine $S$-schemes and object $\eta \colon Y \arr V$ of $\cF(V)$. Notice that, given a positive integer $N$, the isomorphisms $\rho_{f, \phi} \colon f^*\cL_\eta \simeq \cL_\xi$ as in the statement of the theorem induce isomorphisms $\rho_{f, \phi}^{\otimes N} \colon  f^*(\cL_\eta^{\otimes N}) \simeq \cL_\xi^{\otimes N}$. These also satisfy the conditions of the theorem:  whenever we have a diagram of schemes
   \[
   \xymatrix{
   X \ar[r]^f \ar[d]^\xi & Y \ar[r]^g \ar[d]^\eta & Z \ar[d]^\zeta \\
   U \ar[r]^\phi         & V \ar[r]^\psi           & W
   }
   \]
whose columns are in $\cF$, then the diagram
   \[
   \xymatrix@C+30pt{
   f^*g^*\cL_\zeta^{\otimes N} \ar[r]^-{\isoass_{f,g}(\cL_\zeta^{\otimes N})}
      \ar[d]^{f^* \rho_{g, \psi}^{\otimes N}}
      & (gf)^*\cL_\zeta^{\otimes N} \ar[d]^{\rho_{gf, \psi\phi}^{\otimes N}}
   \\
   f^* \cL_\eta  \ar[r]^-{\rho_{f, \phi}^{\otimes N}}&
      {}\cL_\xi ^{\otimes N}
   }
   \]
commutes (this is easily checked by following the action of the arrows on sections of the form $f^{*}g^{*}s^{\otimes N}$, where $s$ is a section of $\cL_{\zeta}$ over some open subset of $W$: since those generate $f^{*}g^{*}\cL_{\zeta}^{\otimes N}$, if the two composites $f^{*}g^{*}\cL_{\zeta} \arr \cL_{\xi}$ agree on them, they must be equal).

By substituting $\cL_{\eta}$ with $\cL_{\eta}^{\otimes N}$ for a sufficiently large integer $N$, we may assume that $\cL_\eta$ is very ample on $Y$, and for any point $v \in V$ we have $\H^1(Y_v, \cL_\eta \mid_{Y_v}) = 0$, where $Y_v$ is the fiber of $Y$ over $v$. This will have the consequence that all the base change homomorphisms that intervene in the following discussion are isomorphisms.

 The two diagrams
   \[
   \xymatrix{
   Y \times_U V \ar[r]^-{\pr_1} \ar[d]_{\eta \times \id_V}
   & Y \ar[d]^\eta   \\
   V \times_U V \ar[r]^-{\pr_1}    & V
   }
   \qquad \text{and} \qquad
   \xymatrix{
   V \times_U Y\ar[r]^-{\pr_2} \ar[d]_{\id_V \times \eta}
   & Y \ar[d]^\eta  \\
   V \times_U V \ar[r]^-{\pr_2} & V
   }
   \]
are cartesian; therefore we can take $\eta \times \id_V \colon Y \times_U V \arr V \times_U V$ as the pullback $\pr_1^* \eta \in \cF(V \times_U V)$, and analogously, $\id_V \times f \colon V \times_U V \arr V \times_U V$ as $\pr_2^* \eta \in \cF(V \times_U V)$.

Analogously, the pullbacks of $f$ along the three projections $V \times_U V \times_U V \arr V$ are
   \begin{align*}
   \eta \times \id_V \times \id_V &\colon Y \times_U V \times_U V 
   \arr V \times_U V \times_U V, \\
   \id_V \times \eta \times \id_V &\colon
   V \times_U Y \times_U V \arr V \times_U V \times_U V
   \quad\text{and} \quad\\
   \id_V \times \id_V \times \eta &\colon
   V \times_U V \times_U Y \arr V \times_U V \times_U V.
   \end{align*}

Suppose that we are given an object $\eta\colon Y \arr V$ in $\cF(V)$ with descent data $\phi \colon V \times_{U} Y \simeq Y \times_{U} V$, that consists of an isomorphism of schemes over $V \times_{U} V$ satisfying the cocycle condition, that is the commutativity of the diagram
   \[
   \xymatrix{
   V \times_{U} V \times_{U} Y
      \ar[rr]^-{\pr_{23}^{*}\phi \ =\ \id_{V}\times\phi}
      \ar[rd] _-{\pr_{13}^{*}\phi}
   && V \times_{U} Y \times_{U} V
      \ar[ld]^-{\hskip10pt\pr_{12}^{*}\phi\ =\ \phi\times\id_{V}}
   \\
   & Y\times_{U} V \times_{U} V
   }
   \]
We will use $\phi$ to construct descent data for the \qc sheaf of finite presentation $\cM_{\eta}$ on $V$.
From the two cartesian diagrams
   \[
   \xymatrix{
   V \times_U Y\ar[r]^-{\pr_2} \ar[d]_{\id_V \times \eta}
   & Y \ar[d]^\eta
   \\
   V \times_U V \ar[r]^-{\pr_2} & V
   }
   \qquad \text{and} \qquad
   \xymatrix{
   Y \times_U V \ar[r]^-{\pr_1} \ar[d]_{\eta \times \id_V}
   & Y \ar[d]^\eta
   \\
   V \times_U V \ar[r]^-{\pr_1}    & V
   }
   \]
we get isomorphisms
   \[
   \sigma_{\pr_{2},\pr_{2}} \colon
      \pr_{2}^{*}\cM_{\eta} \simeq \cM_{\id_V \times \eta}
   \quad
   \text{and}
   \quad
   \sigma_{\pr_{1},\pr_{1}} \colon
      \pr_{1}^{*}\cM_{\eta} \simeq \cM_{\eta\times\id_V};
   \]
and from the cartesian diagram
   \[
   \xymatrix{
   V\times_{U} Y \ar[r]^-{\phi} \ar[d]^{\id_{V}\times\eta}
   & Y\times_{U} V \ar[d]^{\eta\times \id_{V}}
   \\
   V\times_{U} V \ar@{=}[r]
   &V\times_{U} V
   }
   \]
another isomorphism
   \[
   \sigma_{\phi,\id_{V\times_{U}V}}\colon \cM_{\eta\times\id_V} 
      \simeq \cM_{\id_V \times \eta}.
   \]
With these we define an isomorphism
   \[
   \psi \eqdef \sigma_{\pr_{1},\pr_{1}}^{-1} \circ
      \sigma_{\phi,\id_{V\times_{U}V}}^{-1} \circ
      \sigma_{\pr_{2},\pr_{2}}\colon 
      \pr_{2}^{*}\cM_{\eta} \simeq \pr_{1}^{*}\cM_{\eta}
   \]
of \qc sheaves on $V \times_{U} V$.

Let us check that $\psi$ satisfies the cocycle condition. We will use our customary notation $\pr_{1}$ and $\pr_{2}$ to denote the projection onto the first and second factor of a product $X_{1} \times_{Y} X_{2}$, and $p_{1}$, $p_{2}$ and $p_{3}$ for the first, second and third projection onto the factors of the triple product $X_{1} \times_{Y} X_{2} \times_{Y} X_{3}$. (Previously we have also denoted these by $\pr_{1}$, $\pr_{2}$ and $\pr_{3}$, but here the risk of confusion seems more real.)

Consider the cartesian diagram
   \[
   \xymatrix@C+20pt{
   V\times_{U} Y\times_{U} V \ar[r]^{\pr_{12}} \ar@/^1.5pc/[rr]^{p_{2}}   
      \ar[d]^{\id_{V}\times \eta\times \id_{V}}
   & V\times_{U} Y \ar[r]^{\pr_{2}} \ar[d]^{\id_{V}\times \eta}
   & Y \ar[d]^{\eta}
   \\
   V\times_{U} V\times_{U} V \ar[r]^{\pr_{12}} \ar@/_1.5pc/[rr]^{p_{2}}
   &V\times_{U} V \ar[r]^{\pr_{2}}
   &V
   };
   \]
according to Proposition~\ref{prop:compatibility-sigma}, we have that $\sigma_{p_{2},p_{2}}\colon p_{2}^{*}\cM_{\eta} \arr \cM_{\id_{V}\times \eta\times \id_{V}}$ is the composite
   \[
   p_{2}^{*}\cM_{\eta} \xarr{\pr_{12}^{*}\sigma_{\pr_{2},\pr_{2}}}
   \pr_{12}^{*}\cM_{\id_{V}\times\eta} \xarr{\sigma_{\pr_{12}. \pr_{12}}}
   \cM_{\id_{V}\times \eta\times \id_{V}},
   \]
so we have the equality
   \[
   \pr_{12}^{*}\sigma_{\pr_{2},\pr_{2}} =
      \sigma_{\pr_{12},\pr_{12}}^{-1}\circ \sigma_{p_{2},p_{2}}\colon 
   p_{2}^{*}\cM_{\eta} \arr \cM_{\id_{V}\times \eta\times \id_{V}}.
   \]
In a completely analogous fashion we get the equalities
   \begin{align*}
   \pr_{12}^{*}\sigma_{\pr_{1},\pr_{1}} &=
   \sigma_{\pr_{12},\pr_{12}}^{-1} \circ 
      \sigma_{p_{1},p_{1}},\\
   \pr_{23}^{*}\sigma_{\pr_{1},\pr_{1}} &=
   \sigma_{\pr_{23},\pr_{23}}^{-1} \circ 
      \sigma_{p_{2},p_{2}},\\
   \pr_{23}^{*}\sigma_{\pr_{2},\pr_{2}} &=
   \sigma_{\pr_{23},\pr_{23}}^{-1} \circ 
      \sigma_{p_{3},p_{3}},\\
   \pr_{13}^{*}\sigma_{\pr_{1},\pr_{1}} &=
   \sigma_{\pr_{13},\pr_{13}}^{-1} \circ 
      \sigma_{p_{1},p_{1}}\quad\text{and}\\
   \pr_{13}^{*}\sigma_{\pr_{2},\pr_{2}} &=
   \sigma_{\pr_{13},\pr_{13}}^{-1} \circ 
      \sigma_{p_{3},p_{3}}.\\
   \end{align*}
We need to prove the equality $\pr_{12}^{*}\psi \circ \pr_{23}^{*}\psi = \pr_{13}^{*}\psi$; using definition of $\psi$ and the identities above, and doing some simplifications, the reader can check that this equality is equivalent to the equality
   \begin{align*}
      (\sigma_{\pr_{12}, \pr_{12}} \circ \pr_{12}^{*}\sigma_{\phi}^{-1}
      \circ \sigma_{\pr_{12}, \pr_{12}}^{-1})
   &\circ
   (\sigma_{\pr_{23}, \pr_{23}} \circ \pr_{23}^{*}\sigma_{\phi}^{-1}
      \circ \sigma_{\pr_{23}, \pr_{23}}^{-1}) \\
   &=
   (\sigma_{\pr_{23}, \pr_{23}} \circ \pr_{23}^{*}\sigma_{\phi}^{-1}
      \circ \sigma_{\pr_{23}, \pr_{23}}^{-1})
   \end{align*}
that we are going to prove as follows.

From the cartesian diagram
   \[
   \xymatrix{
   V\times_{U} Y\times_{U} V \ar[r]^-{\pr_{12}}
      \ar[d]^{\id_{V}\times \eta\times \id_{V}}
   & V\times_{U} Y \ar[r]^-{\phi} \ar[d]^{\id_{V}\times \eta}
   & Y\times_{U}V \ar[d]^{\eta\times \id_{V}}
   \\
   V\times_{U} V\times_{U} V \ar[r]^-{\pr_{12}}
   &V\times_{U} V \ar@{=}[r]
   &V\times_{U} V
   }
   \]
we get the equality
   \[
   \pr_{12}^{*}\sigma_{\phi}^{-1} \circ \sigma_{\pr_{12}, \pr_{12}}^{-1}
      = \sigma_{\phi\circ\pr_{12},\pr_{12}}^{-1};
   \]
from the other cartesian diagram
   \[
   \xymatrix@C+20pt{
   V\times_{U}Y\times_{U}V \ar[r]^{\phi\times\id_{V}}
      \ar[d]^{\id_{V}\times\eta\times\id_{V}}
      \ar@/^1.5pc/[rr]^{\phi\circ p_{12}}   
   & Y\times_{U}V\times_{U}V \ar[r]^-{\pr_{12}}
      \ar[d]^{\eta\times\id_{V}\times\id_{V}}
   & Y \times_{U} V\ar[d]^{\eta\times\id_{V}}
   \\
   V\times_{U}V\times_{U}V \ar@{=}[r]
   & V\times_{U}V\times_{U}V \ar[r]^-{\pr_{12}}
   & V\times_{U}V
   }
   \]
we obtain that
   \[
   \sigma_{\pr_{12},\pr_{12}} \circ \sigma_{\phi\circ\pr_{12},\pr_{12}}^{-1}
      = \sigma_{\phi\times\id_{V},\id_{V\times_{U}V\times_{U}V}}^{-1}
      = \sigma_{\pr_{12}^{*}\phi,\id_{V\times_{U}V\times_{U}V}}^{-1}.
   \]
With analogous arguments we get the equalities
   \[
   \sigma_{\pr_{23},\pr_{23}} \circ \sigma_{\phi\circ\pr_{23},\pr_{23}}^{-1}
      = \sigma_{\pr_{23}^{*}\phi,\id_{V\times_{U}V\times_{U}V}}^{-1}
   \]
and
   \[
   \sigma_{\pr_{13},\pr_{13}} \circ \sigma_{\phi\circ\pr_{13},\pr_{13}}^{-1}
      = \sigma_{\pr_{13}^{*}\phi,\id_{V\times_{U}V\times_{U}V}}^{-1},
   \]
from which we see that the cocycle condition for $\psi$ is equivalent to the equality
   \[
   \sigma_{\pr_{23}^{*}\phi,\id_{V\times_{U}V\times_{U}V}}\circ
   \sigma_{\pr_{12}^{*}\phi,\id_{V\times_{U}V\times_{U}V}} =
   \sigma_{\pr_{13}^{*}\phi,\id_{V\times_{U}V\times_{U}V}}.
   \]
But this follows immediately, once again thanks to Proposition~\ref{prop:compatibility-sigma}, from the cocycle condition on $\phi$.

So $(\cM_{\eta},\psi)$ is a \qc sheaf with descent data, hence it will come from some \qc sheaf of finite presentation $\cM$ on $U$.

Now let us go back to the general case. Given an arrow $\xi\colon X \arr U$ in $\cF$, we have an adjunction homomorphism
   \[
   \tau_{\xi}\colon \xi^{*}\cM_{\xi} = \xi^{*}\xi_{*}\cL_{\xi} \arr \cL_{\xi}
   \]
that is characterized at the level of sections by the equality $\tau_{\xi}(\xi^{*}s) = s$ for any section $s$ of $\cL_{\xi}$ over the inverse image of an open subset of $U$.

\begin{proposition}\label{prop:compatibility-adjunction}
Given a cartesian square
   \[
   \xymatrix{
   X \ar[r]^{f} \ar[d]^{\xi} 
   &Y\ar[d]^{\eta}
   \\
   U\ar[r]^{\phi}
   & V
   }
   \]
the two composites
   \[
   (\phi\xi)^{*}\cM_{\eta} = (\eta f)^{*}\cM_{\eta} \xarr{f^{*}\tau_{\eta}}
      f^{*}\cL_{\eta} \xarr{\rho_{f,\phi}}
      \cL_{\xi}
   \]
and
   \[
   (\phi\xi)^{*}\cM_{\eta} \xarr{\xi^{*}\sigma_{f,\phi}}
      \xi^{*}\cM_{\xi} \xarr{\tau_{\xi}}
      \cL_{\xi}
   \]
coincide.
\end{proposition}

\begin{proof}
Let $s$ be a section of $\cM_{\eta}$ over an open subset of $V$, that is, a section of $\cL_{\eta}$ over the inverse image in $Y$ of an open subset of $V$. Then both composites are characterized by the property of sending $(\phi\xi)^{*}s = (\eta f)^{*}s$ to $\rho_{f,\phi}(f^{*}s)$.
\end{proof}

In the situation of the Proposition above, assume that $\cL_{\eta}$ is very ample on $Y$ relative to $\eta$, and that the base change homomorphism $\beta_{f,\phi}\colon \phi^{*}\eta_{*}\cL_{\eta} \arr \xi_{*}f^{*}\cL_{\eta}$ is an isomorphism. Then $\sigma_{f,\phi}\colon \phi^{*}\cM_{\eta} \arr \cM_{\xi}$ is an isomorphism, and this induces an isomorphism
   \[
   U\times_{V} \PP(\cM_{\eta}) = \PP(\phi^{*}\cM_{\eta}) \simeq \PP(\cM_{\xi})
   \]
of schemes over $U$, hence a cartesian square
   \[
   \xymatrix{
   {}\PP(\cM_{\xi}) \ar[r]\ar[d]
   & U\ar[d]^{\phi}
   \\
   {}\PP(\cM_{\eta}) \ar[r]
   & V
   }
   \]
Also, since $\cL_{\xi}$ and $\cL_{\eta}$ are very ample, the base change homomorphisms $\tau_{\xi}\colon \xi^{*}\cM_{\xi} \arr \cL_{\xi}$ and $\tau_{\eta}\colon \eta^{*}\cM_{\eta} \arr \cL_{\eta}$ are surjective, and the corresponding morphisms of schemes $X \arr \PP(\cM_{\xi})$ and $Y \arr \PP(\cM_{\eta})$ are closed embeddings. Proposition~\ref{prop:compatibility-adjunction} implies that the diagram
   \[
   \xymatrix{
   X\ar[d]^{f} \ar@{^(->}[r]
   & {}\PP(\cM_{\xi}) \ar[r]\ar[d]
   & U\ar[d]^{\phi}
   \\
   Y \ar@{^(->}[r]
   & {}\PP(\cM_{\eta}) \ar[r]
   & V
   }
   \]
commutes, and is cartesian.

Going back to our covering $V \arr U$, we have that the two diagrams
   \[
   \xymatrix{
   V\times_{U}Y\ar[d]^{\pr_{1}} \ar@{^(->}[r]
   & {}\PP(\cM_{\id_{V}\times \eta}) \ar[r]\ar[d]
   & V\times_{U}V \ar[d]^{\pr_{1}}
   \\
   Y \ar@{^(->}[r]
   & {}\PP(\cM_{\eta}) \ar[r]
   & V
   }
   \]
and
   \[
   \xymatrix{
   Y\times_{U}V\ar[d]^{\pr_{2}} \ar@{^(->}[r]
   & {}\PP(\cM_{\eta\times\id_{V}}) \ar[r]\ar[d]
   & V\times_{U}V \ar[d]^{\pr_{2}}
   \\
   Y \ar@{^(->}[r]
   & {}\PP(\cM_{\eta}) \ar[r]
   & V
   }
   \]
are cartesian; these, together with the diagram
   \[
   \xymatrix{
   V\times_{U}Y\ar[d]^{\phi} \ar@{^(->}[r]
   & {}\PP(\cM_{\id_{V}\times \eta}) \ar[r]\ar[d]
   & V\times_{U}V \ar@{=}[d]
   \\
   Y\times_{U} V \ar@{^(->}[r]
   & {}\PP(\cM_{\eta}\times\id_{V}) \ar[r]
   & V\times_{U}V
   }
   \]
and the definition of $\psi$, show that the two inverse images of $Y \subseteq \PP(\cM_{\eta})$ in $\PP(\pr_{1}^{*}\cM_{\eta})$ and in $\PP(\pr_{2}^{*}\cM_{\eta})$ coincide. On the other hand, the quasi-coherent sheaf with descent data $(\cM_{\eta},\psi)$ is isomorphic to the pullback of the quasi-coherent sheaf $\cM$. If $f\colon V \arr U$ is the given morphism, $F\colon V \times_{U} V \arr U$ the composite with the projections $V\times_{U}V \arr V$, this implies that the two pullbacks of $Y \subseteq \PP(\cM_{\eta}) \simeq \PP(f^{*}\cM)$ coincide. This implies that there is a unique closed subscheme $X \subseteq \PP(\cM)$ that pulls back to $Y \subseteq \PP(\cM_{\eta})$; and, since $\cF$ is a local class in the fpqc topology, the morphism $X \arr U$ is in $\cF$. The scheme with descent data associated with $X \arr U$ is precisely $(Y \to V, \phi)$; and this completes the proof of Theorem~\ref{thm:descent-ample}.
\end{proof}

\section{Descent along torsors}

One of the most interesting examples of descent is descent for \qc sheaves along fpqc torsors. This can be considered as a vast generalization of the well know equivalence between the category of real vector spaces and the category of complex vector spaces with an anti-linear involution. Torsors are generalizations of principal fiber bundles in topology; and I always find it striking that among the simplest examples of torsors are Galois field extensions (see Example~\ref{ex:galois-extensions}).

Here we only introduce the bare minimum of material that allows us to state and prove the main theorem. For a fuller treatment, see \cite{demazure-gabriel}.

In this section we will work with a subcanonical site $\cC$ with fibered products, and a group object $G$ in $\cC$. We will assume that $\cC$ has a terminal object $\pt$.

The examples that we have in mind are $\cC = \cattop$, endowed with the global classical topology, where $G$ is any topological group, and $\cC = \catsch{S}$, with the fpqc topology, where $G \arr S$ is a group scheme.

\subsection{Torsors}

Torsors are what in other fields of mathematics are called \emph{principal bundles}. Suppose that we have an object $X$ of $\cC$, with a left action $\alpha \colon G \times X \arr X$ of $G$. An arrow $f\colon X \arr Y$ is called \emph{invariant}\index{invariant arrow}\index{arrow!invariant} if for each object $U$ of $\cC$ the induced function $X(U) \arr Y(U)$ is invariant with respect to the action of $G(U)$ on $X(U)$. Another way of saying this is that the composites of $f\colon X \arr Y$ with the two arrows $\alpha$ and $\pr_2$ from $G \times X$ to $X$ are equal (the equivalence with the definition above follows from Yoneda's lemma).

Yet another equivalent definition is that the arrow $f$ is $G$-equivariant, when $Y$ is given the trivial $G$-action $\pr_2 \colon G \times Y \arr Y$.

If $\pi \colon X \arr Y$ is an invariant arrow and $f \colon Y' \arr Y$ is an arrow, there is an induced action of $G$ on $Y' \times_{Y} X$; this is the unique action that makes the first projection $\pr_1 \colon Y' \times_{Y} X \arr Y'$ invariant, and the second projection $\pr_2 \colon Y' \times_{Y} X \arr X$ $G$-equivariant. In functorial terms, if $U$ is an object of $\cC$, $g \in G(U)$, $x \in X(U)$ and $y' \in Y'(U)$ are elements with the same image in $Y(U)$, we have $g \cdot (y', x) = (y', g \cdot x)$.

The first example of a torsor is the \emph{trivial torsor}. For each object $Y$ of $\cC$, consider the product $G \times Y$. This has an action of $G$, defined by the obvious formula
   \[
   g \cdot (h, y) = (gh,y)
   \]
for all objects $U$ of $\cC$, all $g$ and $h$ in $G(U)$ and all $y$ in $Y(U)$.

More generally, a \emph{trivial torsor}\index{torsor!trivial}\index{trivial torsor} consists of an object $X$ of $\cC$ with a left action of $G$, together with an invariant arrow $f \colon X \arr Y$, such that there is a $G$-equivariant isomorphism $\phi \colon G \times Y \simeq X$ making the diagram
   \[
   \xymatrix@-15pt{
   G \times Y\ar[rr]^{\phi} \ar[rd]_-{\pr_2} &&
   X \ar[ld]^f \\
   & Y
   }
   \]
commute. (The isomorphism itself is not part of the data, only its existence is required.)

A $G$-torsor is an object $X$ of $\cC$ with an action of $G$ and an invariant arrow\ $\pi\colon X \arr Y$ that locally on $Y$ is a trivial torsor. Here is the precise definition.

\begin{definition}
A \emph{$G$-torsor}\index{torsor} in $\cC$ consists of an object $X$ of $\cC$ with an action of $G$ and an invariant arrow\ $\pi \colon X \arr Y$, such that there exists a covering $\{Y_i \arr Y\}$ of $Y$ with the property that for each $i$ the arrow $\pr_1 \colon Y_i \times_Y X \arr Y_i$ is a trivial torsor.
\end{definition}

Here is an important characterization of torsors. Notice that every time we have an action $\alpha \colon G \times X \arr X$ of $G$ on an object $X$ and an invariant arrow $f \colon X \arr Y$, we get an arrow $\delta_\alpha \colon G \times X \arr X \times_Y X$, defined as a natural transformation by the formula $(g, x) \mapsto (gx, x)$ for any object $U$ of $\cC$ and any $g \in G(U)$ and $x \in X(U)$.

\begin{proposition}\call{prop:char-torsor}
\index{torsor!characterization of}
Let $X$ be an object of $\cC$ with an action of $G$. An invariant arrow $\pi \colon X \arr Y$ is a $G$-torsor if and only if

\begin{enumeratei}

\itemref{1} There exists a covering $\{Y_i \arr Y\}$ such that every arrow $Y_i \arr Y$ factors through $\pi \colon X \arr Y$, and

\itemref{2} the arrow $\delta_\alpha \colon G \times X \arr X \times_Y X$ is an isomorphism.

\end{enumeratei}
\end{proposition}

Notice that part~\refpart{prop:char-torsor}{1} says that $X \arr Y$ is a covering in the saturation of the topology of $\cC$ (Definition~\ref{def:saturation}).

\begin{proof}
Assume that the two conditions are satisfied. The arrow $\delta_\alpha$ is immediately checked to be $G$-equivariant; hence the pullback $X \times_Y X \arr X$ of $\pi$ through the covering $\pi \colon X \arr Y$ is a trivial torsor, and therefore $\pi \colon X \arr Y$ is a $G$-torsor.

Conversely, take a torsor $\pi \colon X \arr Y$. First of all, assume that $\pi \colon X \arr Y$ is a trivial torsor, and fix a $G$-equivariant isomorphism $\phi \colon G \times Y \simeq X$ over $Y$. There is a section $Y \arr X$ of $\pi\colon X \arr Y$, so condition~\refpart{prop:char-torsor}{1} is satisfied for the covering $\{Y = Y\}$.

To verify condition~\refpart{prop:char-torsor}{2}, notice that $\delta_\alpha$ can be written as the composite of isomorphisms
   \[
   G \times X \xrightarrow{\id_G \times\phi^{-1}}
   G \times G \times Y \simeq
   (G \times Y) \times_Y (G \times Y) \xrightarrow{\phi \times \phi}
   X \times_Y X,
   \]
where the isomorphism is in the middle is defined as a natural transformation by the rule $(g, h, y) \mapsto \bigl((g,y), (h,y)\bigr)$ for any object $U$ of $\cC$ and any $g, h \in G(U)$ and $y \in Y(U)$.

In the general case, when $\pi \colon X \arr Y$ is not necessarily trivial, the result follows from the previous case and the following lemma.

\begin{lemma}\label{lem:local-isom}
Let $\cC$ be a subcanonical site, 
   \[
   \xymatrix@-10pt{
   X \ar[rr]^f \ar[rd] &&
   Y \ar[ld]\\
   &S
   }
   \]
a commutative diagram in $\cC$. Suppose that there is a covering $\{S_i \arr S\}$ such that the induced arrows
   \[
   \id_{S_i} \times f \colon S_i \times_S X
   \arr S_i \times_S Y
   \]
are isomorphisms. Then $f$ is also an isomorphism.
\end{lemma}

\begin{proof}
The site $(\cC/S)$ is subcanonical (Proposition~\ref{prop:comma-subcanonical}): this means that we can substitute $(\cC/S)$ for $\cC$, and suppose that $S$ is a terminal object of $\cC$.

By Yoneda's lemma, it is enough to show that for any object $U$ of $\cC$ the function $f_U \colon X(U) \arr Y(U)$ induced by $f$ is a bijection. First of all, assume that the arrow $U \arr S$ factors through some $S_i$. By hypothesis $\id_{S_i} \times f \colon S_i \times X \arr S_i \times Y$ is an isomorphism, hence $\id_{S_i(U)} \times f_U \colon S_i(U) \times X(U) \arr S_i(U) \times Y(U)$ is a bijection. If $S_i(U) \neq \emptyset$ it follows that $f_U$ is a bijection.

For the general case we use the hypothesis that $\cC$ is subcanonical. If $U$ is arbitrary, and we set $U_i \eqdef S_i \times U$, then $\{U_i \to U\}$ is a covering. Hence we have a diagram of sets
   \[
   \xymatrix{
   X(U)\ar[r]\ar[d]^{f_{U}} &
   {}\prod_{i} X(U_i) \ar[d]^{\prod_i f_{U_i}} \ar@<3pt>[r]\ar@<-3pt>[r] &
   {}\prod_{ij} Y(U_{ij}) \ar[d]^{\prod_{ij} f_{U_{ij}}} \\
   Y(U)\ar[r] & {}\prod_{i} X(U_i) \ar@<3pt>[r]\ar@<-3pt>[r] &
   {}\prod_{ij} Y(U_{ij}) 
   }
   \]
in which the rows are equalizers, because $\cC$ is subcanonical. On the other hand each arrow $U_i \arr S$ and $U_{ij} \arr S$ factors through $S_i$, so $f_{U_i}$ and $f_{U_{ij}}$ are bijections. It follows that $f_U$ is a bijection, as required.
\end{proof}

Consider the arrow $\delta_\alpha \colon G \times X \arr X \times_Y X$, and choose a covering $\{Y_i \arr Y\}$ such that for each $i$ the pullbacks $X_i \eqdef Y_i \times_Y X$ are trivial as torsors over $Y_i$. Denote by $\alpha_i \colon G \times X_i \arr X_i$ the induced action; then $\delta_{\alpha_i} \colon G \times X_i \arr X_i \times_{Y_i} X_i$ is an isomorphism. On the other hand there are standard isomorphisms $(G \times X) \times_{Y} Y_i \simeq G \times X_i$ and $(X \times_Y X) \times_Y Y_i \simeq X_i \times_{Y_i} X_i$, and the diagram
   \[
   \xymatrix@C+15pt{
   (G \times X)\times_{Y}Y_i \ar[r]^-{\delta_\alpha\times\id_{Y_i}}
   \ar[d]^{\pr_1} &
   (X \times_Y X)\times_Y Y_i \ar[d]^{\pr_1} \\
   G \times X_i \ar[r]^{\delta_{\alpha_i}} &
   X_i \times_{Y_i} X_i
   }
   \]
commutes. Hence $\delta_\alpha\times\id_{Y_i}$ is an isomorphism for all $i$, and it follows that $\delta_\alpha$ is an isomorphism.
\end{proof}

\begin{example}\label{ex:galois-extensions}
Let $K \subseteq L$ be a finite Galois extensions, with Galois group $G$. Denote by $G_{K}$ the discrete group scheme $G\times \spec K \arr \spec K$ associated with $G$, as in \S\ref{subsec:discrete-groups}. The action of $G$ on $L$ defines an action $\alpha \colon G_{K}\times~_{\spec K} \spec L = G \times \spec L \arr \spec L$ of $G_{K}$ on $\spec L$ (Proposition~\ref{prop:action-discrete-group}), which leaves the morphism $\spec L \arr \spec K$ invariant. (For convenience we will write the action of $G$ on $L$ on the right, so that the resulting action of $G$ on $\spec L$ is naturally written as a left action.)

By the primitive element theorem, $L$ is generated as an extension of $K$ by a unique  element $u$; denote by $f \in K[x]$ its minimal polynomial. Then $L = K[x]/\bigl(f(x)\bigr)$. The group $G$ acts on the roots of $f$ simply transitively, so $f(x) = \prod_{g \in G}(x - ug) \in L[x]$.

The morphism
   \[
   \delta_\alpha \colon G_{K} \times \spec L
   \arr \spec L \times_{\spec K} \spec L
      = \spec(L\otimes_{K} L)
   \]
corresponds to the homomorphism of $K$-algebras $L \otimes_{K}L \arr L^{G}$ defined as
   \[
   a \otimes b \mapsto \bigl((ag)b\bigr)_{g\in G},
   \]
where by $L^{G}$ we mean the product of copies of $L$ indexed by $G$. We have an isomorphism
   \[
   L \otimes_{K} L = K[x]/\bigl( \prod_{g \in G}(x - ug)\bigr)\otimes_{K} L 
   \simeq L[x]/\bigl( \prod_{g \in G}(x - ug)\bigr);
   \]
by the Chinese remainder theorem, the projection
   \[
   L[x]/\bigl( \prod_{g \in G}(x - ug)\bigr) \arr \prod_{g \in G}L[x]/(x-ug)
   \simeq L^{G}
   \]
is an isomorphism. Thus we get an isomorphism $L\otimes_{K} L \simeq L^{G}$, that is easily seen to coincide with the homomorphism corresponding to $\delta_{\alpha}$. Thus $\delta_{\alpha}$ is an isomorphism; and since $\spec L \arr \spec K$ is \'etale, this shows that $\spec L$ is $G_{K}$-torsor over $\spec K$.
\end{example}

Here is our main result.

\begin{theorem}
\index{theorem!descent!along torsors}
\index{descent!along torsors}
Let $X \arr Y$ be a $G$-torsor, and $\cF \arr \cC$ a stack. Then there exists a canonical equivalence of categories between $\cF(Y)$ and the category of $G$-equivariant objects $\cF^{G}(X)$ defined in \S\ref{sec:fibered-actions}.
\end{theorem}

\begin{proof}
Because of Proposition~\ref{prop:same-stacks}, we have an equivalence of $\cF(Y)$ with $\cF(X \to Y)$, so it is enough to produce an equivalence between $\cF(X \to Y)$ and $\cF^{G}(X)$.

For this we need the isomorphism $\delta_{\alpha} \colon G \times X \simeq X \times_{Y} X$ defined above, and also the one defined in the next Lemma.

\begin{lemma}
If $X \arr Y$ is a $G$-torsor, the arrow
   \[
   \delta_{\alpha}' \colon G \times G \times X 
      \arr X \times_{Y} X \times_{Y} X
   \]
defined in functorial terms by the rule
   \[
    \delta_{\alpha}'(g, h, x) = (ghx, hx, x)
   \]
is an isomorphism.
\end{lemma}

Once again, one reduces to the case of a trivial torsor using Lemma~\ref{lem:local-isom}. We leave the proof of this case to the reader.

Since the category $\cF(X \to Y)$ does not depend on the choice of the fibered products $X \times_{Y} X$ and $X \times_{Y} X \times_{Y} X$, we can make the choice $X \times_{Y} X = G\times X$ and $X \times_{Y} X \times_{Y} X = G \times G \times X$, in such a way that $\delta_{\alpha}$ and $\delta_{\alpha}'$ become the identity. Then we have
   \begin{align*}
   \pr_{1} &= \alpha \colon G \times X \arr X,\\
   \pr_{13}&= \mul_{G} \times \id_{X} \colon
      G \times G \times X \arr G \times X,\\
   \pr_{23}&= \id_{G} \times \alpha \colon
      G \times G \times X \arr G \times X,
   \end{align*}
while $\pr_{23}= X \times_{Y} X \times_{Y} X \arr X \times_{Y} X$ coincides with the projection $G \times G \times X \arr G \times X$ on the second and third factor.

Then an object $(\rho, \phi)$ of $\cF(X \to Y)$ is an object $\rho$ of $\cF(X)$, together with an isomorphism $\phi\colon \pr_{2}^{*}\rho \simeq \alpha^{*}\rho$ satisfying the cocycle condition: and the cocycle condition is precisely the condition for $\phi$ to define a $G$-equivariant structure on $\rho$, according to Proposition~\ref{prop:definition-equivariant}. Hence the category $\cF(X\to Y)$ is canonically isomorphic to $\cF^{G}(X)$ (for this we need to check what happens to arrows, but this is easy and left to the reader), and this concludes the proof of the theorem.
\end{proof}

\subsection{Failure of descent for morphisms of schemes}
\label{subsec:failure}

Now we construct an example to show how descent can fail for proper smooth
morphisms of proper schemes of finite type over a field.

The starting point is a variant on Hironaka's famous example of a nonprojective threefold\index{example!Hironaka} (\cite{hironaka62}, \cite[Chapter~3, \S~3]{git}), \cite[Appendix~B, Example~3.4.1]{hartshorne}); this has been already been used to give examples of a smooth three-dimensional algebraic space over a field that is not a scheme (\cite[p.~14]{knutson}).

Fix an algebraically closed field $\kappa$. Then one constructs a smooth proper connected three dimensional scheme $M$ over $\kappa$, with an action of a cyclic group of order two $\mathrm{C}_2 = \{1, \sigma\}$, containing two copies $L_1$ of $L_2$ of $\PP^1$ that are interchanged by $\sigma$, with the property that the 1-cycle $L_1 + L_2$ is algebraically equivalent to $0$. This implies that there is no open affine subscheme $U$ of $M$ that intersects $L_1$ and $L_2$ simultaneously: if not, the complement $S$ of $U$ would be a surface in $M$ that intersects both $L_1$ and $L_2$ in a finite number of points. But since $S \cdot (L_1 + L_2) = 0$ this finite number of points would have to be zero, and this would mean that $L_1$ and $L_2$ are entirely contained in $U$. This is impossible, because $U$ is affine.

Now take a $\mathrm{C}_2$-torsor $V \arr U$ (a Galois \'etale cover with group $\mathrm{C}_2$) with $V$ irreducible, and set $Y = M \times_\kappa V$. The projection $\pi \colon Y \arr V$ is smooth and proper. We need descent data for the covering $V \arr U$; these are given by the diagonal action of $\mathrm{C}_2$ on $Y$, obtained from the two actions on $M$ and $V$. More precisely, the action $\rC_{2}\times Y \arr Y$ gives a cartesian diagram
   \[
   \xymatrix{
   {}\rC_{2}\times Y \ar[r]\ar[d] & Y\ar[d] \\
   {}\rC_{2}\times V \ar[r]       & V
   }
   \]
yielding an isomorphism of $\rC_{2} \times Y$ with the pullback of $Y \arr V$ to $\rC_{2}\times V$, and defines an object with descent data on the covering $V \arr V$ (keeping in mind that $V \times_{U} V = \rC_{2}\times V$, since $V \arr U$ is a $\rC_{2}$-torsor).

I claim that these descent data are not effective. Suppose that it is not so: then there is a cartesian diagram of schemes
   \[
   \xymatrix{
   Y \ar[r]^f \ar[d]^\pi & X\ar[d] \\
   V \ar[r]              & U
   }
   \]
such that $f$ is invariant under the action of $\mathrm{C}_2$ on $Y$. Take an open affine subscheme $W \subseteq X$ that intersects $f(L_1 \times V) = f(L_2 \times V)$; then its inverse image $f ^{-1}W \subseteq Y$ is affine, and if $p$ is a generic closed point of $V$ the intersection $\pi^{-1} p \cap f ^{-1}W \subseteq M$ is an affine open subscheme of $M$ that intersects both $L_1$ and $L_2$. As we have seen, this is impossible.

Another counterexample, already mentioned in Example~\ref{ex:stacks-curves}, is in \cite[XIII 3.2]{raynaudample}.

\begin{remark}\label{rmk:algebraic-spaces}
There is an extension of the theory of schemes, the theory of \emph{algebraic spaces}\index{algebraic space}, due to Michael Artin (see \cite{artin-algspaces}, \cite{artinrep} and \cite{knutson}). An algebraic space over a scheme $S$ is an \'etale sheaf $\catsch{S}\op \arr \catset$, that is, in some sense, \'etale locally a scheme. The category of algebraic stacks contains the category of schemes over $S$ with quasi-compact diagonal; furthermore, by a remarkable result of Artin, it is a stack in the fppf topology (it is probably also a stack in the fpqc topology, but I do not know this for sure: however, for most applications fppf descent is what is needed). Also, most of the concepts and techniques that apply to schemes extend to algebraic spaces. This is obvious for properties of schemes, and morphisms of schemes, such us being Cohen--Macaulay, smooth, or flat, that are local in the \'etale topology (on the domain). Global properties, such as properness, require more work.

So, in many contexts, when some descent data in the fppf topology fail to define a scheme, an algebraic space appears as a result. Also, algebraic spaces can be used to define stacks in situations when descent for schemes fails. For example, if we redefine that stack $\cF_{1,S}$ of Example~\ref{ex:stacks-curves} so that the objects are proper smooth morphisms $X \arr U$ whose fibers are curves of genus~$1$, where $U$ is an $S$-scheme and $X$ is an algebraic space, then $\cF_{1,S}$ is a stack in the fppf topology.
\end{remark}

\nocite{moerdijk02}
\nocite{knus-ojanguren-descent}
\nocite{benabou-roubaud}
\nocite{giraud}
\nocite{giraud-descent}

\bibliographystyle{amsalpha}
\bibliography{mrabbrev,VistoliRefs}

\printindex

\end{document}